\documentclass[a4paper,11pt,reqno]{article}

\usepackage{authblk}
\usepackage[T1]{fontenc}
\usepackage[latin1]{inputenc}
\usepackage{empheq,amsthm,amssymb}

\usepackage{mathrsfs}

\usepackage{xcolor}
\usepackage{graphicx,caption,subcaption}
\graphicspath{{FIGS/}}

\usepackage{hyperref}

\usepackage[margin=2.55cm]{geometry}

\numberwithin{equation}{section}

\newtheorem{theorem}{Theorem}[section]

\newtheorem{proposition}[theorem]{Proposition}

\theoremstyle{definition}

\newtheorem{assumption}{Assumption}

\newcommand{\calA}{\mathcal{A}}
\newcommand{\calB}{\mathcal{B}}
\newcommand{\calC}{\mathcal{C}}

\newcommand{\calE}{\mathcal{E}}

\newcommand{\calG}{\mathcal{G}}

\newcommand{\calM}{\mathcal{M}}

\newcommand{\calO}{\mathcal{O}}
\newcommand{\calP}{\mathcal{P}}

\newcommand{\calV}{\mathcal{V}}
\newcommand{\calW}{\mathcal{W}}
\newcommand{\calX}{\mathcal{X}}
\newcommand{\calY}{\mathcal{Y}}

\newcommand{\sfN}{\mathsf{N}}

\newcommand{\bbE}{\mathbb{E}}

\newcommand{\bbL}{\mathbb{L}}
\newcommand{\bbM}{\mathbb{M}}
\newcommand{\bbN}{\mathbb{N}}

\newcommand{\bbP}{\mathbb{P}}

\newcommand{\bbR}{\mathbb{R}}
\newcommand{\bbS}{\mathbb{S}}

\newcommand{\scrB}{\mathscr{B}}

\newcommand{\scrV}{\mathscr{V}}
\newcommand{\scrW}{\mathscr{W}}
\newcommand{\scrX}{\mathscr{X}}
\newcommand{\scrY}{\mathscr{Y}}

\newcommand{\vertiii}[1]{{\left\vert\kern-0.25ex\left\vert\kern-0.25ex\left\vert #1
    \right\vert\kern-0.25ex\right\vert\kern-0.25ex\right\vert}}

\begin{document}

\title{\vspace{-.75cm}On the Mathematical Theory of Ensemble (Linear-Gaussian) Kalman-Bucy Filtering}

\author[1]{Adrian N. Bishop*}
\author[2]{Pierre Del Moral}
\affil[1]{{\small CSIRO; and University of Technology Sydney (UTS), Australia}}
\affil[2]{{\small INRIA, Bordeaux Research Center, France}}
\date{}

\maketitle

\vspace{-1.cm}
\begin{abstract}
	The purpose of this review is to present a comprehensive overview of the theory of ensemble Kalman-Bucy filtering for continuous-time, linear-Gaussian signal and observation models. We present a system of equations that describe the flow of individual particles and the flow of the sample covariance and the sample mean in continuous-time ensemble filtering. We consider these equations and their characteristics in a number of popular ensemble Kalman filtering variants. Given these equations, we study their asymptotic convergence to the optimal Bayesian filter. We also study in detail some non-asymptotic time-uniform fluctuation, stability, and contraction results on the sample covariance and sample mean (or sample error track). We focus on testable signal/observation model conditions, and we accommodate fully unstable (latent) signal models. We discuss the relevance and importance of these results in characterising the filter's behaviour, e.g. it's signal tracking performance, and we contrast these results with those in classical studies of stability in Kalman-Bucy filtering.We also provide a novel (and negative) result proving that the bootstrap particle filter cannot track even the most basic unstable latent signal, in contrast with the ensemble Kalman filter (and the optimal filter). We provide intuition for how the main results extend to nonlinear signal models and comment on their consequence on some typical filter behaviours seen in practice, e.g. catastrophic divergence.
\end{abstract}

{\scriptsize
\setcounter{tocdepth}{2}
\tableofcontents
}

\section{Introduction}

Consider a time-invariant, continuous-time, signal and observation model of the form,
\begin{equation}\label{nonlinear-diffusion-filtering-intro}
\begin{split}
d\scrX_t~&=~a(\scrX_t)\,dt \,+\, R^{1/2}\,d\scrV_t \\
d\scrY_t~&=~h(\scrX_t)\,dt \,+\, R_1^{1/2}\,d\scrW_t
\end{split}
\end{equation}
where $\scrX_t$ is the underlying signal (latent) process, $\scrY_t$ is the observation signal, $a(\cdot)$ and $h(\cdot)$ are the signal and sensor model functions, and $\scrV_t$ and $\scrW_t$ are continuous-time Brownian motion (noise) signals. The filtering problem \cite{anderson79,bainCrisan} is concerned with estimating some statistic(s) of the signal $\scrX_t$ conditioned on the observations $\scrY_s$, $0\leq s\leq t$. For example, one may want to characterise fully the distribution of $\scrX_t$ given $\scrY_t$, or one may seek some moments of this distribution. The conditional distribution of $\scrX_t$ given $\scrY_s$, $0\leq s\leq t$ is called the (optimal, Bayesian) filtering distribution. When the model functions $a(\cdot)$, $h(\cdot)$ are linear, the exact (optimal, Bayesian) solution to this problem is completely characterised by the first two moments of the filtering distribution and these moments are given by the celebrated Kalman-Bucy filter \cite{kalman61,anderson79,Bishop/DelMoral:2016}. 

Apart from the most special of nonlinear models, there is in general no finite dimensional optimal filter \cite{Benes1981,bainCrisan}. In practice, some filter approximations are needed. For example, one may consider a type of ``extended'' Kalman filter \cite{anderson79} based on linearisation of the nonlinear model and application of the classical Kalman-Bucy filter. This method works well in suitably regular, and sufficiently close to linear problems. This method does not handle well multiple modes in the true filtering distribution. So-called Gaussian-sum filters are another Kalman-filter-type/based approximation designed to handle in some sense multiple modes in the filtering distribution \cite{anderson79}. More recently there has been some focus on Monte Carlo integration methods for approximating the optimal Bayesian filter \cite{mf-dm-04,bainCrisan}. Such methods, termed particle filters or sequential Monte Carlo filters/methods \cite{gordon1993,doucetSMCpaper2000,doucetSMCbook2001}, have the advantage of not being subject to the assumption of linearity or Gaussianity in the model. These particle filters are consistent in the number of Monte Carlo samples, i.e. with infinite computational power these methods converge to the optimal nonlinear filter. However, typical particle filtering algorithms exhibit high computational costs with approximation errors that grow (with a fixed sample size) with the signal/observation dimensions \cite{mf-dm-04,Rebeschini2015}. These methods are not scalable to the high-dimensional filtering or state estimation problems found in the geosciences and other areas \cite{evensen2009book,kalnay2003atmospheric,law2015data,Reich2015book}.

The ensemble Kalman-Bucy filter (generally abbreviated {\tt EnKF}) \cite{evensen03,evensen2009book} is a type of Monte Carlo sample approximation of a class of linear (in the observations) filter in the spirit of the Kalman filter. The {\tt EnKF} is a recursive algorithm for propagating and updating the sample mean and sample covariance of an approximated Bayesian filter \cite{evensen2009book}. The filter works via the evolution of a collection (i.e. an ensemble) of samples (i.e. ensemble members, or particles) that each satisfies a type of Kalman-Bucy update equation, linear in the observations. In classical Kalman-Bucy filtering \cite{kalman61,anderson79,Bishop/DelMoral:2016}, a gain function, that depends on the filter error covariance, is used to weight a predicted state estimate with the signal observations, see \cite{evensen2009book,Bishop/DelMoral:2016}. In the {\tt EnKF}, the error covariance in the gain function is replaced by a type of sample covariance. The result is a system of interacting particles in the spirit of a mean-field approximation of a certain McKean-Vlasov-type diffusion equation \cite{mckean1966class,DelMoral/Tugaut:2016}. We may refine this discussion by giving the relevant equations for a most basic form of {\tt EnKF}. Let $(\calV^i_t,\calW^i_t,\calX_0^i)$ with ${1\leq i\leq \sfN+1}$ be $(\sfN+1)$ independent copies of $(\scrV_t,\scrW_t,\scrX_0)$. The most basic ensemble Kalman filter, originally due to Evensen \cite{burgers1998analysis,evensen03,evensen2009book}, is defined by,
\begin{align}\label{EnKF-sampled-diffusions-nonlinear-intro}
d\calX_t^i~=&~a(\calX_t^i)\,dt~+~R^{1/2}\,d\calV^i_t+\widehat{P}^{\,h}_{t}\,R_1^{-1}\left[d\scrY_t-\left(h(\calX_t^i)\,dt+R_1^{1/2}\,
d\calW_{t}^i\right)\right]
\end{align}
with $1\leq i\leq \sfN+1$ and the (particle) sample mean and the sample cross-covariance defining the so-called Kalman gain matrix given by,
\begin{equation}\label{fv1-3-2-intro}
\begin{array}{l}
\displaystyle
\widehat{X}_{t} :=\frac{1}{\sfN+1}\sum_{i=1}^{\sfN+1}\calX_t^i \quad\mathrm{and}\quad \displaystyle  \widehat{P}^{\,h}_t := \frac{1}{\sfN}\,\sum_{i=1}^{\sfN+1}\left[\calX_t^i-\widehat{X}_{t}\right]\left[h(\calX_t^i)-\frac{1}{\sfN+1}\sum_{i=1}^{\sfN+1}h(\calX_t^i)\right]^{\prime}
\end{array}\end{equation}
and we may also write the standard sample covariance,
\begin{equation}\label{fv1-3-2-nonlinear-intro}
\widehat{P}_{t} :=\,\frac{1}{\sfN}\,\sum_{i=1}^{\sfN+1}\left[\calX_t^i-\widehat{X}_{t}\right]\left[\calX_t^i-\widehat{X}_{t}\right]^{\prime}
\end{equation}
In this work we study this most basic ensemble Kalman filter as described above, and also more sophisticated variants, including the method of Sakov and Oke \cite{sakov2008deterministic}, that exhibit less fluctuation due to sampling noise. Readers familiar with the Kalman filter will recognise immediately some structural similarities as discussed above. However, there is no evolution equation given above for the covariance as in the Kalman filter (e.g. no Riccati-type matrix flow equation). Instead, we replace the relevant covariance matrices with their sample-based counterparts. 

Importantly, if the underlying model is linear and Gaussian, then the filtering distribution is Gaussian, and the {\tt EnKF} propagates exactly the sample mean and covariance of the optimal Bayesian filter, and is provably consistent. If the model is nonlinear and/or non-Gaussian, then a standard implementation of the {\tt EnKF} propagates a sample-based estimate of the filtering mean and covariance (but not the true posterior sample mean or covariance and with no results on consistency). In the context of estimation theory, we may contrast the notion of a state estimator (or observer) with the notion of a Bayesian filter. The goal of the former is to design an observer that tracks in some suitable (typically point-wise) sense the underlying signal and perhaps provides some usable measure of uncertainty on this estimate. The goal of the latter is to compute or approximate the true (Bayesian) filtering distribution (or some related statistics). In the nonlinear setting, even with infinite computational power, the {\tt EnKF} methods do not converge to the optimal nonlinear filter; and indeed their limiting objects are not well understood in this setting. As discussed more technically later, ensemble Kalman filters are probably best viewed in practice as a type of (random) sample-based state estimator for nonlinear signal/observation models. However, in the special case of linear signal and observation models they are indeed provably consistent approximations of the optimal Bayesian filter. 

In practice, the ensemble Kalman filtering methodology is applied in high-dimensional, nonlinear state-space models, e.g. see \cite{evensen03,evensen2009book} and the application references listed later in this introduction. Empirically, this method has shown good tracking performance in these applications, see \cite{evensen2009book} and the application references listed later. This tracking behaviour of the {\tt EnKF} when applied to practical models may be explainable by viewing the {\tt EnKF} as a dynamic state estimator. The fluctuation, stability, and contraction properties of the {\tt EnKF} studied in this article (albeit mainly for linear-Gaussian models) may be viewed in this context also, and provide some insight into the state estimate tracking behaviour seen in practice.

\subsection{Purpose}

The purpose of this review is to present a comprehensive overview of the theory of ensemble Kalman-Bucy filtering with an emphasis on rigorous results and behavioural characterisations for linear-Gaussian signal/observation models. We present a system of equations that describe the flow of individual particles, the flow of the sample covariance, and the flow of the sample mean in continuous-time ensemble filtering. We consider these equations and their characteristics in a number of popular {\tt EnKF} varieties. Given these equations, \emph{we study in detail some fluctuation, stability, and error contraction results for the various ensemble Kalman filtering equations}. We discuss the relevance and importance of these results in terms of characterising the {\tt EnKF} behaviour, and we contrast these results with those considered in classical studies of stability in Kalman-Bucy filtering.

Classical studies of stability in (traditional, non-ensemble-type) Kalman-Bucy filtering are important because they rigorously establish the type of ``tracking'' properties desired in a filtering or estimation problem; and they establish intuitive, testable, model-based conditions (e.g. model observability) for achieving these convergence properties. Classical results in Kalman-Bucy filtering also establish the (exponential) convergence of the error covariance to a fixed steady-state value computable from the model parameters. See the review \cite{Bishop/DelMoral:2016,bd-CARE} for detailed results in the classical context and historical remarks. \emph{The results in this work seek to characterise in an analogous manner the practical performance and behaviour of ensemble Kalman filtering, and these results then provide guidance and intuition on the tracking, approximation error, and other properties of these practical methods. Notably, the stochastic fluctuation properties of ensemble Kalman methods also need to be established; and counterparts of this latter analysis do not arise at all in classical Kalman filtering analyses.} Our results are presented under testable, model-based assumptions. In particular, we rely on the standard controllability assumption from classical Kalman filtering theory; and, typically, a more restrictive (but testable) observability-type assumption (i.e. linear fully-observed processes, which imply classical observability).

\subsection{Overview of the Main Topics and Literature}
 
In this subsection we touch on the main topics and related literature as it pertains to the {\tt EnKF}. These topics include the fluctuation, stability, and contractive properties of the relevant {\tt EnKF} stochastic equations. Later, toward the end of this article, we discuss some of these topics in the context of filtering and state estimation more broadly, and we touch on other related but somehow distinct results as they pertain to the {\tt EnKF} more specifically.

The {\tt EnKF} is a key numerical method for solving high-dimensional forecasting and data assimilation problems; see, e.g., \cite{evensen03,evensen2009book}. In particular, applications have been motivated by inference problems in ocean and atmosphere sciences~\cite{Lisaeter2003,majda2012,kalnay2003atmospheric,ott2004}, weather forecasting~\cite{anderson2001ensemble,anderson2003local,burgers1998analysis,houtekamer1998data}, environmental and ecological statistics~\cite{allen2003,johns2008}, as well as in oil reservoir simulations~\cite{evensen2007reservoir,nydal2003,seiler2010}, and many others. This list is by no means exhaustive, nor the cited articles fully representative of the respective applications. We refer to (some of) the seminal methodology papers in \cite{evensen1994orig,houtekamer1998data,burgers1998analysis,andersonanderson1999,houtekamer2001seq,hamill01,CHBishop2001,anderson2001ensemble,Whitaker2002,anderson2003local,tippetsqrt2003,sakov2008deterministic,Reich2013,Yang2016}. This long list is not exhaustive; see also the books \cite{evensen2009book,kalnay2003atmospheric,law2015data,Reich2015book} for more background, and the detailed chronological list of references in Evensen's text \cite{evensen2009book}. 

In continuous-time, we may broadly break down the class of {\tt EnKF} methods into three distinct types; distinguished by the level of fluctuation added via sampling noise needed to ensure that the {\tt EnKF} sample mean and covariance are consistent in the linear-Gaussian setting. The original form of the {\tt EnKF} is the so-called `vanilla' {\tt EnKF} of Evensen \cite{burgers1998analysis,evensen2009book}, see also \cite{Lange2021a}; and this method exhibits the most fluctuations due to sampling of both signal and observation noises. The next class is the so-called  `deterministic' {\tt EnKF} of Sakov and Oke \cite{sakov2008deterministic}, see also \cite{Bergemann2012,Reich2013}; which exhibits (considerably) less fluctuation. In the continuous-time linear-Gaussian setting this class is representative of the so-called square-root {\tt EnKF} methods \cite{tippetsqrt2003,Lange2021b} (which differ somewhat in discrete-time, e.g. contrast \cite{sakov2008deterministic} with \cite{tippetsqrt2003}, see also \cite{Lange2021b}). Finally, there has been recent interest in so-called transport-inspired {\tt EnKF} methods \cite{Reich2013,Taghvaei2016ACC}; which apart from initialisation noise/randomisation are completely deterministic and whose analysis in the linear model setting follows closely that of the classical Kalman-Bucy filter, cf. \cite{Bishop/DelMoral:2016}. These classes do not distinguish the totality of {\tt EnKF} methodology (especially in nonlinear or non-Gaussian models); which may further consist of so-called covariance regularisation methods \cite{andersonanderson1999,houtekamer2001seq,hamill01,Mitchell2002,evensen03}, etc. However, in the linear-Gaussian case, these three classes broadly capture the fundamentals.

As discussed later, the fully deterministic, transport-inspired {\tt EnKF} method, see \cite{Reich2013,Taghvaei2016ACC}, is a rather special case in the linear-Gaussian setting and is not studied in detail in this article where linear-Gaussian models are the focus. Nevertheless, we point to \cite{deWiljes2018,deWiljes2019} for certain mean-field consistency results, non-asymptotic fluctuation (e.g. finite sample size) results, and the long-time behaviour of this particular method in the case of a nonlinear signal model and linear observations. We also touch on this method briefly throughout; but when we refer to the general {\tt EnKF} we typically mean the so-called `vanilla' \cite{burgers1998analysis,Lange2021a} or `deterministic' \cite{sakov2008deterministic,Bergemann2012} methods (which will become clear as the article progresses).

Convergence to a mean-field limit, and large-sample asymptotics, of the discrete-time {\tt EnKF} was studied in~\cite{legland,mandel2011convergence,Kwiatkowski2015,Lange2020}, in the sense of taking the number of particles to infinity. The discrete-time square root form of the {\tt EnKF} is accommodated in \cite{Kwiatkowski2015,Lange2020}, and nonlinear state-space models are accommodated in \cite{legland}. In the continuous-time, linear-Gaussian, setting, the convergence (in sample size) of the three broad classes of {\tt EnKF} to the true Kalman-Bucy filter is more immediate; and follows from the sample mean and sample covariance evolution equations in \cite{DelMoral/Tugaut:2016,Bishop/DelMoral/multiDimRicc}. In this latter sense, we recover the fact that the {\tt EnKF} is a consistent approximation of the optimal, Bayesian filter (i.e. the classical Kalman-Bucy filter) in the linear-Gaussian setting as discussed earlier. The mean-field limit of various {\tt EnKF} methods in the continuous-time, nonlinear model setting is studied in \cite{law2016deterministic,deWiljes2018,Lange2020}.

We remark in the nonlinear model setting (discrete or continuous-time), see \cite{legland,Reich2014,law2016deterministic,deWiljes2018,Lange2020}, the mean-field limiting equations (and distribution) are not easily related to the optimal filter. Moreover, in practice, one is typically interested in the non-asymptotic (in terms of ensemble size) fluctuation properties as well as the long time/stability behaviour of the particle-type filtering approximations. 

The fluctuation analysis of the {\tt EnKF} is studied in detail in the linear-Gaussian setting in \cite{Bishop/DelMoral/Niclas:2017,2017arXiv171110065B,Bishop/DelMoral/multiDimRicc}. In \cite{Bishop/DelMoral/Niclas:2017} a complete Taylor-type stochastic expansion of the sample covariance is given at any order with bounded remainder terms and estimates. Both non-asymptotic and asymptotic bias and variance estimates for the {\tt EnKF} sample covariance and sample mean are given explicitly in \cite{Bishop/DelMoral/Niclas:2017}. These latter expansions directly imply an almost sure strong form of a central limit-type result on the sample covariance and sample mean at any time. The analysis in \cite{Bishop/DelMoral/Niclas:2017} is considered over the entire path space of the matrix-valued Riccati stochastic differential equation that describes the flow of the sample covariance. However, most of the non-asymptotic time-uniform results in \cite{Bishop/DelMoral/Niclas:2017} hold only when the underlying signal is stable. In \cite{2017arXiv171110065B,Bishop/DelMoral/multiDimRicc} we consider the case in which the underlying signal may be unstable, and we provide time-uniform, non-asymptotic moment estimates and time-uniform control over the fluctuation of the sample covariance and mean about their limiting Riccati and Kalman-Bucy filtering terms.

The emphasis of \emph{time-uniformity} on the moment bounds and on the fluctuation bounds on the sample mean and sample covariance (about the true optimal Bayesian filtering mean and covariance) is important. If these bounds are allowed to grow in time, e.g. typically in this analysis one can easily obtain bounds that grow exponentially in time, then these bounds quickly become useless for any practical numerical application; e.g. an exponent $>200$ may induce an exceedingly pessimistic bound greater than the estimated number of particles of matter in the visible universe. We remark also that our emphasis on {\emph{accommodating unstable (latent) signal models} is important because time-uniform fluctuation results in such cases (which are of real practical importance) are significantly more difficult to obtain under testable and realistic model assumptions (like the classical observability and controllability model assumptions in the control and filtering literature \cite{anderson79,Bishop/DelMoral:2016}).

In \cite{DelMoral/Tugaut:2016}, stability of the {\tt EnKF} in continuous-time linear-Gaussian models is considered under the assumption that the underlying signal model is also stable. This latter assumption is in contrast with classical Kalman-Bucy filter stability results, which hold in the linear-Gaussian setting under the much weaker (and more natural) condition of signal detectability \cite{VanHandel2009,Bishop/DelMoral:2016,bd-CARE}. The classical Kalman-Bucy filter is stable as a result of the closed-loop stabilising properties of the so-called Kalman gain matrix, which is closely connected to the flow of the filter error covariance described by a Riccati differential equation. The {\tt EnKF} analogue, in linear-Gaussian settings, is the sample covariance, and its random fluctuation properties (noted in the preceding paragraph) are the main source of difficulty in establishing the closed-loop filter stability in those models in which the underlying signal itself is unstable.

In \cite{tong2016stability}, the authors analyse the long-time behaviour of the (discrete-time) {\tt EnKF} in a class of nonlinear systems, with finite ensemble size, using Foster-Lyapunov techniques. Applying the results of \cite{tong2016stability} to the basic linear-Gaussian filtering problem, the analysis and assumptions in \cite{tong2016stability} then also require stability of the underlying signal model. In a traditional sense, the conditions needed in \cite{tong2016stability} are hard to check, e.g. as compared to the classical observability or controllability-type model conditions in Kalman filtering analysis; but a range of examples are given in \cite{tong2016stability}. In \cite{Kelly2014}, the long-time behaviour of the {\tt EnKF} is analysed in both discrete and continuous time settings with similar conditions on the model as in \cite{tong2016stability}; and which again if linearised equates to a form of stability on the signal model. 

We emphasise again that the type of analysis in \cite{Kelly2014,tong2016stability,DelMoral/Tugaut:2016} cannot handle unstable, or transient, signal models; i.e. signals with sample paths with at least one coordinate that may grow unbounded. In the context studied in \cite{Kelly2014,tong2016stability,DelMoral/Tugaut:2016} dealing with stable or bounded latent signal processes (e.g. the Lorenz-class of signal models \cite{Kelly2014,tong2016stability}), the important question on the filter stability or filtering error estimates relies on obtaining meaningful quantitative fluctuation constants decreasing with the number of ensemble members to achieve a desired performance. Of course, time uniformity of these bounds follows trivially in this setting from the boundedness properties of the latent signal process.

Covariance inflation is a mechanism used in practical methods to increase the positive-definiteness of the sample covariance matrix and essentially amplify its effect on the stabilisation properties of the Kalman gain matrix. In \cite{Kelly2014} time-uniform {\tt EnKF} error boundedness results follow under a true signal stability condition and given a sufficiently large variance inflation regime. See also~\cite{tong2016inflation,majda2018} for related stability analysis in the presence of adaptive covariance inflation and projection techniques. In \cite{Bishop/DelMoral/multiDimRicc} in the continuous-time linear-Gaussian setting, the mechanism by which covariance inflation acts to stabilise the ensemble filter is exemplified, see also \cite{Bishop/DelMoral/Pathiraja:2017}. Covariance localisation is studied rigorously in \cite{deWiljes2019} in the case of the fully deterministic, transport inspired ensemble filter \cite{Reich2013,Taghvaei2016ACC}. 

In the continuous-time, and linear-Gaussian setting, the first work to relax the assumption of underlying signal stability for the {\tt EnKF} is in \cite{Bishop/DelMoral/Niclas:2017,2017arXiv171110065B,Bishop/DelMoral/multiDimRicc}. In those articles, latent signals with sample paths that may grow unbounded (to infinity exponentially fast) are accommodated. That work is based on both a fluctuation analysis of the sample covariance and the sample mean \cite{Bishop/DelMoral/Niclas:2017,2017arXiv171110065B,Bishop/DelMoral/multiDimRicc}, followed by studies on the long-time behaviour, e.g. stability properties, of both the sample covariance and mean \cite{2017arXiv171110065B,Bishop/DelMoral/multiDimRicc}. Time-uniform fluctuation properties are given under a type of (strong) signal observability condition. In this setting time-uniformity of these results is non-trivial. This assumption is in keeping with classical Kalman-Bucy filtering and Riccati equation results; and does not require any form of underlying signal stability. As the authors of \cite{tong2016stability} note in their stability analysis, they use ``few properties of the forecast [predicted] covariance matrix other than positivity''. As noted in \cite{tong2016stability}, this lends generality to their results, but conversely places the burden back on the signal model assumptions (including those assumptions of true signal stability). Contrast this with the work in \cite{Bishop/DelMoral/Niclas:2017,2017arXiv171110065B,Bishop/DelMoral/multiDimRicc} where emphasis is placed on the fluctuation analysis of the sample covariance, with a primary aim of removing the stability assumptions needed on the underlying signal model. The time-uniform fluctuation and stochastic perturbation contributions in \cite{Bishop/DelMoral/Niclas:2017,2017arXiv171110065B,Bishop/DelMoral/multiDimRicc}, were discussed earlier. Given this fluctuation analysis, the stability of the filter sample mean and sample covariance and their (time) asymptotic properties are studied in \cite{Bishop/DelMoral/Niclas:2017,2017arXiv171110065B,Bishop/DelMoral/multiDimRicc} without stability assumptions on the underlying signal model. These results rigorously establish the type of ``tracking'' properties desired by a filtering or estimation solution.

Although of lesser practical use in applications, strong results in the one-dimensional setting are also derived in \cite{2017arXiv171110065B} that converge, e.g. in the limit with the ensemble size, to those properties of the classical Kalman-Bucy filter. For example, we can recover the optimal exponential contraction and filter stability rates, etc. In the multidimensional setting, the decay rates to equilibrium are not sharp, and the stationary measures are not given in closed form.

\subsection{Aims and Contributions}

The main goal of this article is to: 1). present a novel formulation for ensemble filtering in linear-Gaussian, continuous-time, systems that lends itself naturally to analysis; 2). provide detailed fluctuation analysis of the ensemble Kalman-Bucy flow, the sample mean, and the stochastic Riccati equation describing the sample covariance; 3). study the stability of the resulting stochastic Riccati differential equation that describes the flow of the sample covariance; 4). study the stability of the continuous-time ensemble Kalman-Bucy update equation that is coupled to this stochastic Riccati equation, and which describes the flow of the sample mean (or the sample mean minus the true signal, i.e. the sample error signal). This article is primarily a review of the literature and results in these directions. The prime focal point of this review are the articles \cite{DelMoral/Tugaut:2016,Bishop/DelMoral/Niclas:2017,2017arXiv171110065B,Bishop/DelMoral/STV2018,Bishop/DelMoral/multiDimRicc}, which focus heavily on the linear-Gaussian model setting. In this review, \emph{an emphasis is placed on deriving time-uniform fluctuation, stability, and contraction results under testable model conditions equivalent and/or closely related to the classical observability and controllability-type model assumptions. Importantly, we do not generally assume the true underlying signal is stable in this review. }

Throughout this review we contrast and discuss the presented results with the broader literature on the rigorous mathematical behaviour of ensemble Kalman-type filtering. For example, we find easily that the sample covariance matrix in the broad class of {\tt EnKF} methods considered is always under-biased when compared to the true covariance matrix. This may motivate, from a pure uncertainty quantification viewpoint, some form of covariance regularisation \cite{andersonanderson1999,houtekamer2001seq,hamill01,Mitchell2002,evensen03}. We provide detailed analysis illustrating the effect of inflation regularisation on stability (similarly to \cite{Kelly2014,tong2016inflation,majda2018}). As another example, we provide strong intuition for so-called catastrophic filter divergence (studied previously in \cite{Harlim2010,Gottwald2013,kelly2015}) based on rigorous (heavy-tailed) fluctuation properties inherent to the relevant sample covariance matrices and their invariant distributions. We contrast the so-called `vanilla' {\tt EnKF} of \cite{burgers1998analysis,evensen2009book} with the  `deterministic' {\tt EnKF} of Sakov and Oke \cite{sakov2008deterministic} in terms of their fluctuation and sample noise characteristics and we show how this affects their respective sample behaviour and stability properties. 

 As with classical (non-ensemble) Kalman filtering, the importance of the results reviewed is in rigorously establishing the type of tracking and stability behaviour desired in filtering applications \cite{anderson79,mf-dm-04,bainCrisan,Bishop/DelMoral:2016}. For example, our results imply conditions under which the initial estimation errors are forgotten, and that the flow of the sample mean converges to the true Kalman filtering (conditional mean) state estimate (and thus the signal) in the average. In the case of the {\tt EnKF}, there must be some emphasis placed on the stochastic behaviour of the ensemble (Monte Carlo) mean and covariance in order to establish filter stability. We also provide the analogue of the error covariance fixed point in classical Kalman filtering \cite{anderson79,Bishop/DelMoral:2016}; whereby we state results that ensure the sample covariance matrix converges to an invariant, steady-state, distribution. We characterise the properties of this invariant distribution and relate this to the sample behaviour of the `vanilla' {\tt EnKF} \cite{burgers1998analysis,evensen2009book} and the  `deterministic' {\tt EnKF} \cite{sakov2008deterministic}.

We focus on the linear, continuous-in-time, Gaussian setting in this review and note that in this case the sample mean and sample covariance are consistent approximations of the optimal Bayesian filtering mean and covariance. We emphasise that even in the linear-Gaussian case, the samples themselves are not in general independent. The analysis even in the linear setting is highly technical \cite{DelMoral/Tugaut:2016,Bishop/DelMoral/Niclas:2017,2017arXiv171110065B,Bishop/DelMoral/STV2018,Bishop/DelMoral/multiDimRicc} and the results presented in this case are aimed as a step in the progression to more applied results and intuition in nonlinear model settings. There is some precedent for studying the relative properties, behaviour, or performance of ensemble Kalman filtering firstly with linear-Gaussian signal models \cite{evensen03}. For example, the seminal article \cite{burgers1998analysis} illustrated that a perturbation of the observations in the ensemble Kalman filter was necessary to recover a consistent covariance limit (to the true Kalman filter for linear-Gaussian systems); or to achieve the standard Monte Carlo error rate with a finite set of particles. The analysis (and even derivation) of ensemble square root filters for linear-Gaussian system models is standard \cite{SakovOkeSRF2008,Livings08}, etc. Convergence of the ensemble Kalman filter in inverse problems is studied in \cite{Schillings2017a} in the linear setting. We discuss connections and extensions of the results in this article to the nonlinear model setting toward the end.

We also briefly contrast the approximation capabilities of particle filtering (sequential Monte Carlo) methods \cite{gordon1993,doucetSMCpaper2000} with the {\tt EnKF}. We give a revealing, and perhaps surprising, simple result illustrating the complete failure of the bootstrap particle filter \cite{gordon1993} to track unstable linear-Gaussian latent signals. Compared to the {\tt EnKF}, the fluctuation and stability of various particle filtering methods (e.g. see \cite{dm-g-1999,mf-dm-04,Chopin2004,Oudjane2005,VanHandel2009bb,mf-dm-13,Douc2014}) is a rather mature topic. Nevertheless, time-uniform particle filtering estimates rely on mixing-type, or certain contractive, conditions on the mutation transition which do not hold in general in the case of unstable linear-Gaussian models. We contrast this new (rather negative) particle filtering result with its (positive) counterpart for the {\tt EnKF}.

Note that the analysis and proofs in \cite{DelMoral/Tugaut:2016,Bishop/DelMoral/Niclas:2017,2017arXiv171110065B,Bishop/DelMoral/STV2018,Bishop/DelMoral/multiDimRicc}, while motivated originally by ensemble Kalman-type filtering methods, are largely presented as independent technical results on certain general classes of matrix-valued Riccati diffusion equations and associated linear stochastic differential equations with random coefficients. In this review we emphasise the work in \cite{DelMoral/Tugaut:2016,Bishop/DelMoral/Niclas:2017,2017arXiv171110065B,Bishop/DelMoral/STV2018,Bishop/DelMoral/multiDimRicc} via a series of results directly and solely stated in the context of ensemble Kalman-type filtering. Throughout we relate our results to the broader technical literature on ensemble Kalman filtering and we emphasise the practical significance of these results, e.g. via the tracking property of the filter, its stability, or via their error fluctuation or catastrophic divergence behaviour, among other topics. We also contrast the behaviour of the various classes of continuous-time {\tt EnKF} methods.

\subsection{Notation}

We remark firstly that some care must be taken throughout to keep track of the font stylings; e.g. upright vs. calligraphic vs. script, etc. There is typically a relationship between like symbols appearing with different stylings. 

Hatted terms ~$\widehat{\cdot}$~ should be viewed as being indexed to the ensemble size $\sfN\geq1$, i.e. ~$\widehat{\cdot}:=\cdot^\sfN$. Time is indexed variously by $s,t,u,\tau\in[0,\infty[$. We write $c,c_{n},c_{\tau},c_{n,\tau},c_{n,\tau}(Q),c_{n,\tau}(z,Q)\ldots$ for some positive constants whose values may vary from result to result, and which only depend on the indexed/referenced parameters $n,\tau,z,Q$, etc, as well as implicitly on the model parameters $(A,H,R,R_1)$ introduced later. Importantly, these constants do not depend on the time horizon $t$, nor on the number of ensemble particles $\sfN$.

Let $\bbM_{d}$ be the set of $(d\times d)$ real matrices with $d\geq 1$ and $\bbM_{d_1,d_2}$ the set of $(d_1\times d_2)$ real matrices. Let $\bbS_d\subset \bbM_{d}$ be the subset of symmetric matrices, and $\bbS^0_d$, and $\bbS^+_d$ the subsets of positive semi-definite and definite matrices respectively. We write $A \geq B$ when $A-B\in \bbS^0_d$; and $A > B$ when $A-B\in \bbS^+_d$. We denote by $0$ and $I$ the null and identity matrices, for any $d\geq 1$. Given $R\in \partial \bbS_d^+:= \bbS_d^0-\bbS_d^+$ we denote by $R^{1/2}$ a (non-unique) symmetric square root of $R$. When $R\in\bbS_d^+$ we choose the unique symmetric square root. We write $A^{\prime}$ the transpose of $A$, and $A_{\mathrm{sym}}=(A+A^{\prime})/2$ its symmetric part. We denote by $\mathrm{Absc}(A):=\max{\left\{\mbox{\rm Re}(\lambda)\,:\,\lambda\in \mathrm{Spec}(A)\right\}} $ its spectral abscissa. We also denote by $\mathrm{Tr}(A)$ the trace. When $A\in\bbS_d$ we let $\lambda_1(A)\geq \ldots\geq \lambda_d(A)$ denote the ordered eigenvalues of $A$. We equip $\bbM_{d}$ with the spectral norm $\Vert A \Vert=\Vert A \Vert_2=\sqrt{\lambda_{1}(AA^{\prime})}$ or the Frobenius norm $\Vert A \Vert=\Vert A \Vert_{\mathrm{Frob}}=\sqrt{\mathrm{Tr}(AA^{\prime})}$. 

Let $\mu(A)$ denote a matrix logarithmic ``norm'' (which can be $<0$), see \cite{strom1975logarithmic}. The logarithmic norm is a tool to study the growth of solutions to ordinary differential equations and the error growth in approximation methods. For any square matrix $A\in\bbM_{d}$, the logarithmic norm is the smallest element in the set $\{h\in\mathbb{R}\,:\, \|\exp(At) \|\leq\exp(ht),\,t\geq0\}$ where $\|\cdot\|$ is any matrix norm and the value $\mu(A)$ may be considered to be indexed to the matrix norm employed. For example, the (2-)logarithmic ``norm'', or spectral log-norm, is given by $\mu(A)=\lambda_{1}(A_{\mathrm{sym}})$. We have $\mu(\cdot)\geq\mathrm{Absc}(\cdot)$ in general, but importantly we note that if $\mathrm{Absc}(\cdot)<0$, then there is a matrix norm $\|\cdot\|$ defining a logarithmic norm such that $\mu(\cdot)<0$, see \cite[Theorem 5]{strom1975logarithmic}.

\section{Kalman-Bucy Filtering}

Consider a time-invariant linear-Gaussian filtering model of the following form,
\begin{equation}\label{lin-Gaussian-diffusion-filtering}
\begin{split}
d\scrX_t~&=~A\,\scrX_t\,dt \,+\, R^{1/2}\,d\scrV_t \\
d\scrY_t~&=~H\,\scrX_t\,dt \,+\, R_1^{1/2}\,d\scrW_t
\end{split}
\end{equation}
where $A\in\bbM_{d}$ and $H\in\bbM_{d_y,d}$ are the signal and sensor model matrices respectively, and $R\in\bbS^0_{d}$ and $R_1\in \bbS^+_{d_y}$ are the respective signal and sensor noise covariance matrices. The noise inputs $\scrV_t$ and $\scrW_t$ are $d$ and $d_y$-dimensional Brownian motions, and $\scrX_0$ is an $d$-dimensional Gaussian random variable (independent of $(\scrV_t,\scrW_t)$) with mean $\bbE(\scrX_0)$ and covariance $P_0\in\bbS_d^0$. 

We let $\scrY_0=0$ and $\calY_t=\sigma\left(\scrY_s,~s\leq t\right)$ be the $\sigma$-algebra generated by the observations. The conditional distribution $\eta_t:=\mbox{\rm Law}\left(\scrX_t~|~\calY_t\right)$ of the signal states $\scrX_t$ given $\calY_t$ is Gaussian with a conditional mean and covariance given by
$$
X_t:=\bbE\left(\scrX_t~|~\calY_t\right)\quad \mbox{\rm and}\quad
P_t:=\bbE\left(\left[\scrX_t-X_t\right]\left[\scrX_t-X_t\right]^{\prime}~|~\calY_t\right).
$$
The mean and the covariance obey the Kalman-Bucy and the Riccati equations 
\begin{eqnarray}
dX_t&=&A\,X_t\,dt+P_t\,H^{\prime}R_1^{-1}\left(d\scrY_t-HX_t\,dt\right)\label{nonlinear-KB-mean}\\
\partial_tP_t&=&\mathrm{Ricc}(P_t) \label{nonlinear-KB-Riccati}
\end{eqnarray}
with the Riccati drift function from $\bbS^0_{d}$ into $\bbS_{r}$ defined for any $Q\in \bbS^0_{d}$ by
\begin{equation}\label{riccdrift}
	\mathrm{Ricc}(Q):=AQ+QA^{\prime}-QSQ+R
\end{equation}
and with,
\begin{equation}\label{Sdef}
	 S:=H^{\prime}R_1^{-1}H
\end{equation}
Importantly, the covariance of the conditional distribution $\mbox{\rm Law}(\scrX_t~|~\calY_t)$ in this case does not depend on the observations $\calY_t$. The error $Z_t := (X_t - \scrX_t)$ satisfies
\begin{eqnarray}
	dZ_t 
	&=&  \left(A-P_t\,S\right)Z_t\,dt +P_t\,H^{\prime}R_1^{-1/2}\,d\scrW_t-R^{1/2}\,d\scrV_t \nonumber\\
	&\stackrel{ law}{=}& \left(A-P_t\,S\right)Z_t\,dt + \left(P_t\,S\,P_t + R\right)^{1/2}d\scrB_t \label{kf-error1}
\end{eqnarray}
where $\scrB_t$ is some independent $d$-dimensional Brownian motion. Here we make use of a martingale representation theorem, e.g. \cite[Theorem 4.2]{karatzas}, see also~\cite{doob}.

Let $\phi_t(Q):=P_t$ denote the flow of the matrix differential equation (\ref{nonlinear-KB-Riccati}) with $P_0=Q\in\bbS^0_d$. Let $\psi_t(z,Q):=Z_t$ denote the flow of the stochastic error (\ref{kf-error1}) with $Z_0=z=(x-\scrX_0)\in\bbR^d$ and $P_t= \phi_t(Q)$. Finally, we denote the flow of the Kalman-Bucy update (\ref{nonlinear-KB-mean}) with $X_0=x\in\bbR^d$ by $\chi_t(x,Q):=X_t$. This notation allows us to reference the flows $\psi_t(z,Q)$, $\phi_t(Q)$, $\chi_t(x,Q)$ with respect to their initialisation at $t=0$ which is useful when we compare flows and study stability.

Throughout this section, we assume that $(A,R^{1/2})$ and $(A,H)$ are controllable and observable pairs in the sense that
\begin{equation}\label{def-contr-obs}
\left[R^{1/2},AR^{1/2}\ldots, A^{r-1}R^{1/2}\right]\quad
\mbox{\rm and}\quad
\left[\begin{array}{c}
H\\
HA\\
\vdots\\
HA^{r-1}
\end{array}
\right] \quad\mbox{have rank $d$}.
\end{equation}
Note that if $R\in\bbS^+_d$ is positive definite, which is quite common in filtering problems, it follows that controllability holds trivially. We consider the observability and controllability Gramians $(\calO_{t},\calC_{t}(\calO))$ and $(\calC_{t},\calO_{t}(\calC))$ associated with the triplet $(A,R,S)$ and defined by
\begin{eqnarray}
\calO_{t}~:=~ \int_{0}^{t}\,e^{-A^{\prime}s}\,S\,e^{-As}\,ds ~\quad\mbox{and}\quad~ \calC_{t}(\calO)&:=& \calO_{t}^{-1}\left[\int_0^t\,e^{-(t-s)A^{\prime}}\,\calO_{s}\,R\,\calO_{s}\,e^{-(t-s)A}\,ds\right]\calO_{t}^{-1} \nonumber
\\~\\
\calC_{t} ~:=~ \int_{0}^{t}\,e^{As}\,R\,e^{A^{\prime}s}\,ds ~\quad\mbox{and}\quad~ \calO_{t}(\calC)&:=&\calC_{t}^{-1}\left[\int_0^t\,e^{(t-s)A}\,\calC_{s}\, S\,\calC_{s}\,e^{(t-s)A}\,ds\right]\calC_{t}^{-1}. \nonumber
\end{eqnarray}
Given (\ref{def-contr-obs}), for any finite $\tau>0$, there exists some finite parameters $\varpi^{o,c}_{\pm},\varpi^{c}_{\pm}(\calO),\varpi^{o}_{\pm}(\calC)>0$ such that
\begin{eqnarray}
\varpi_-^{c} \leq \|\calC_{\tau}\|
\leq \varpi_+^{c}
\quad\mbox{\rm
and}
\quad
\varpi_-^{o} \leq \|\calO_{\tau}\| \leq \varpi_+^{o}
\end{eqnarray}
\begin{eqnarray}
 \varpi_-^{c}(\calO) \leq \|\calC_{\tau}(\calO)\| \leq \varpi_+^{c}(\calO) \quad
\mbox{\rm and}\quad \varpi_-^{o}(\calC) \leq \|\calO_{\tau}(\calC)\| \leq \varpi_+^{o}(\calC).~
\end{eqnarray}
The parameter $\tau$ is often called the interval of observability-controllability, see \cite{bucy67}. 

These rank conditions (\ref{def-contr-obs}) ensure the existence and the uniqueness of a positive definite fixed-point matrix $P_\infty$ solving the algebraic Riccati equation
\begin{equation}\label{steady-state-eq}
\mathrm{Ricc}(P_\infty):=AP_\infty+P_\infty A^{\prime}-P_\infty SP_\infty+R=0.
\end{equation}
Indeed, if (\ref{def-contr-obs}) holds, then $P_\infty\in\bbS_d^+$ and $\mathrm{Absc}(A-P_\infty S)<0$. We may relax the controllability assumption to just stabilisability, in which case $P_\infty\in\bbS_d^0$ and $\mathrm{Absc}(A-P_\infty S)<0$; see \cite{kucera72,Molinari77,Lancaster1995} and the convergence results in \cite{Kwakernaak72,callier81}. Under just a detectability condition, it follows that $P_\infty\in\bbS_d^0$ and $\mathrm{Absc}(A-P_\infty S)\leq0$, i.e. $(A-P_\infty S)$, is only marginally stable, and convergence to this solution is given under mild additional conditions in \cite{Poubelle86,callier95,Park97}. In \cite{VanHandel2009}, given only detectability, the time-varying ``closed loop'' matrix $(A-\phi_t(Q)S)$ is shown to be stabilising, even when $(A-P_\infty S)$ is only marginally stable.
  
  In the context of ensemble Kalman-Bucy filtering considered later, we will require the same controllability assumption as considered above, and a more restrictive observability condition (that implies the classical observability/detectability discussed above).

For any $s\leq t$ and $Q\in\bbS_{d}^0$ we define the state-transition matrix,
\begin{equation}\label{KB-semigroup-def-main}
\calE_{s,t}(Q):=\exp{\left(\oint_s^t\left(A-\phi_u(Q)\,S\right)du\right)} ~~~\Longleftrightarrow~~~ \partial_t \calE_{s,t}(Q)=\left(A-\phi_u(Q)S\right)\calE_{s,t}(Q).
\end{equation}
When $s=0$ we often write $\calE_{t}(Q)$ instead of $\calE_{0,t}(Q)$. The matrix $\calE_{t}(Q)$ is the fundamental matrix. We have $\calE_{s,t}(Q)=\calE_{t}(Q)\calE_{s}(Q)^{-1}$. The following convergence estimates follow from~\cite{Bishop/DelMoral:2016,bd-CARE}: For any $Q,Q_1,Q_2\in\bbS^{0}_{d}$ and any $t\geq 0$ we have the local contraction inequalities
\begin{equation}\label{ref-E-1}
\Vert \mathcal{E}_{t}(Q)\Vert \,\leq\, c\,(1+\Vert Q\Vert)\,\Vert \mathcal{E}_{t}({P}_{\infty})\Vert \qquad \mbox{\rm and}\qquad \Vert \mathcal{E}_{t}({P}_{\infty})\Vert\,= \,\Vert e^{t(A-P_\infty S)} \Vert \,\leq\, c\,e^{-\alpha\, t}
\end{equation}
for some finite $\alpha,c>0$ and with ${P}_\infty$ solving (\ref{steady-state-eq}) and
\begin{equation}
\Vert \calE_t(Q_2)-\calE_t(Q_1)\Vert \,\leq\,  c( Q_1, Q_2)~e^{-2\,\alpha\, t}~\Vert Q_2-Q_1\Vert \label{expo-Estable} 
\end{equation}
for some finite constant $c( Q_1, Q_2)>0$. In addition, there exists some parameter $\tau> 0$ such that for any $s\geq 0$ and any $t\geq \tau>0$ we have the uniform estimates,
\begin{equation}\label{ref-krause-inf}
\Vert \mathcal{E}_{s,s+t}(Q)\Vert \,\leq\,c_\tau\,\Vert \calE_{t}({P}_{\infty})\Vert 
\end{equation}
Note it is desirable to relate the decay of $\mathcal{E}_{s,s+t}(Q)$ to the decay at the fixed point $\Vert \mathcal{E}_{t}({P}_{\infty})\Vert= \Vert e^{t(A-P_\infty S)} \Vert \leq c\,e^{-\alpha\, t}$ (since as $t\rightarrow\infty$ it is clear that we cannot do better). See \cite{bd-CARE} for an explicit Floquet-type expression of $\mathcal{E}_{t}(Q)$ in terms of $\mathcal{E}_{t}(P_\infty)$.

The convergence and stability properties of the Kalman-Bucy filter and the associated Riccati equation are directly related to the contraction properties of the state-transition matrix $\calE_{s,t}(Q)$. To get some intuition for this we note, 
\begin{equation}\label{closed-form-OU}
\psi_t(z,Q) \,=\, \mathcal{E}_{s,t}(Q)\,\psi_s(z,Q)+\int_s^t~\mathcal{E}_{u,t}(Q)\left(\phi_{u}(Q)\,S\,\phi_{u}(Q) + R\right)^{1/2}\,d\scrB_u
\end{equation}
and
\begin{equation}\label{closed-form-ricc-diff}
	\phi_t(Q) \,=\, \mathcal{E}_{s,t}(Q)\,\phi_{s}(Q)\,\mathcal{E}_{s,t}(Q)' + \int_{0}^t\,\mathcal{E}_{s,t}(Q)\left(\phi_{s}(Q)\,S\,\phi_{s}(Q) + R\right)\mathcal{E}_{s,t}(Q)'\,ds
\end{equation}
for any $s\leq t$.

From~\cite{Bishop/DelMoral:2016}, for any $t\geq \tau>0$ and any $Q\in \bbS_d^0$ we have the uniform estimates
\begin{equation}\label{after-upsilon}
\left(\calO_{\tau}(\calC)+ \calC_{\tau}^{-1}\right)^{-1}\,\leq\,\phi_{t}(Q)~\leq \,
\calO_{\tau}^{-1}+\calC_{\tau}(\calO).
\end{equation}
We also have
\begin{equation}\label{upper-bound-Qhi}
	0\,\leq\, \phi_{t}(Q) \,\leq\, P_\infty+e^{(A-P_\infty S)t}(Q-P_\infty)e^{(A-P_\infty S)'t}
\end{equation}
The following stability result follows from \cite{Bishop/DelMoral:2016,bd-CARE}: For any $Q_1,Q_2\in\bbS^{0}_{d}$ and for any $t\geq 0$,
\begin{equation}\label{ref-phi-1-stability}
	\Vert \phi_{t}(Q_1) -  \phi_{t}(Q_2)\Vert \,\leq\, c\,(1+\Vert Q_1\Vert^2+\Vert Q_2\Vert^2)\,\Vert \mathcal{E}_{t}({P}_{\infty})\Vert^2 \,\Vert Q_2-Q_1\Vert
\end{equation}
and recall the exponential contraction estimate on $\Vert \mathcal{E}_{t}({P}_{\infty})\Vert$ in (\ref{ref-E-1}). Similarly, using (\ref{ref-krause-inf}), for any $s\geq 0$ and any $t\geq \tau>0$, we have 
\begin{equation}\label{ref-phi-1-stability-upsilon}
	\Vert \phi_{s,s+t}(Q_1) -  \phi_{s,s+t}(Q_2)\Vert \,\leq\, c_\tau\,\Vert \mathcal{E}_{t}({P}_{\infty})\Vert^2 \,\Vert Q_2-Q_1\Vert
\end{equation}
Note that both (\ref{ref-phi-1-stability}) and (\ref{ref-phi-1-stability-upsilon}) imply immediately that $\phi_{t}(Q)\rightarrow_{t\rightarrow\infty}{P}_\infty$ exponentially fast for any $Q\in \bbS^0_d$; e.g. by letting $Q_2=P_\infty$.

Note that the uniform estimates with constants independent of the initial condition stated throughout, involve some arbitrarily small, positive time parameter $\tau$, which can be directly related to the notion of a so-called observability/controllability interval introduced earlier; for further details on this topic we refer to \cite{bucy67,Bishop/DelMoral:2016}. Contrast, for example, the stability results (\ref{ref-phi-1-stability}) and (\ref{ref-phi-1-stability-upsilon}). The symbol $\tau$ is reserved for this arbitrary small time parameter throughout the article.

Results (e.g. bounds and convergence results) on the flow of the inverse of the solution of the Riccati equation are considered in \cite{Bishop/DelMoral:2016} and are relevant for proving results on the flow of the Riccati equation itself; e.g. upper bounds on the flow of the inverse solution help to lower bound solutions of the Riccati flow. The flow of the inverse Riccati solution may also be of interest on its own as it relates to the flow of ``information'' (as the inverse of covariance).  

Given the contraction properties on $\mathcal{E}_{s,t}(Q)$ it is often said the ``deterministic part'' of the filter error $\partial_tZ_t=\left(A-P_t\,S\right)Z_t$ is stable. From \cite{Bishop/DelMoral:2016} we can be more explicit if desired, for example, for any $t\geq \tau$ we have the uniform estimate,
\begin{equation}\label{KB-bias-conv}
\sup_{Q\in\bbS^0_{d}} \left\Vert\,  \bbE\left( \psi_{t}(z,Q) \,\vert\, \scrX_0\right)\, \right\Vert
~\leq~ c\,e^{-\alpha\, t}\,\Vert\,x-\scrX_0 \Vert
\end{equation}
for some rate $\alpha>0$ and some finite constant $c>0$. Moreover, the conditional probability of the following
event
\begin{equation}\label{event-control-KB-bis}
\left\Vert\psi_{t}(z,Q)\right\Vert~ \leq~  c(Q)\left( e^{-\alpha t}~\Vert x-\scrX_0\Vert + \frac{e^2}{\sqrt{2}}\left[\frac{1}{2}+\left(\delta+\sqrt{\delta}\right)\right] \right)
\end{equation}
given the state variable $\scrX_0$ is greater than $1-e^{-\delta}$, for any $\delta\geq 0$. And, for any $t\geq 0$, $z_1,z_2\in \bbR^{d}$, $Q_1,Q_2\in\bbS^0_{d}$ and any $n\geq 1$ we have the almost sure local contraction estimate
\begin{equation}\label{control-KF-error}
\begin{array}{l}
\bbE\left(\Vert\psi_{t}(z_1,Q_1)-\psi_{t}(z_2,Q_2)\Vert^{n}~\vert~\scrX_0\right)^{\frac{1}{n}}
\\
\\
\qquad\qquad\leq~ \displaystyle c(Q_1,Q_2)~e^{-\alpha t}~\Vert z_1-z_2 \Vert +  c_n(Q_1,Q_2)\,e^{-\alpha t}\,
\left(1+\left\Vert x_1- \scrX_0\right\Vert\right)\,\Vert Q_1-Q_2\Vert
\end{array}
\end{equation}
with some rate $\alpha>0$ and the finite constants $c(Q_1,Q_2),c_n(Q_1,Q_2)>0$.

\section{Kalman-Bucy Diffusion Processes}

For any probability measure $\eta$ on $\bbR^d$ we let $\calP_{\eta}$ denote the $\eta$-covariance
\begin{equation}
\eta\mapsto\calP_{\eta}:=\eta\left([\iota-\eta(\iota)][\iota-\eta(\iota)]'\right)
\end{equation}
with the identity function $\iota(x):=x$ and the column vector $\eta(f):=\int f\, d\eta$ for some measurable function $f:\bbR^d\rightarrow\bbR^d$.

We now consider three different cases of a conditional nonlinear McKean-Vlasov-type diffusion process,
\begin{align}\label{Kalman-Bucy-filter-nonlinear-ref}
(\texttt{F1})\qquad d\calX_t~=&~A\,\calX_t~dt~+~R^{1/2}\,d\calV_t+\calP_{\overline{\eta}_t}~H^\prime\,R_1^{-1}~\left[d\scrY_t-\left(H\calX_tdt+R_1^{1/2}~
d\calW_{t}\right)\right]\nonumber\\
(\texttt{F2})\qquad  d\calX_t~=&~A\,\calX_t~dt~+~R^{1/2}\,d\calV_t+\calP_{\overline{\eta}_t}\,H^\prime\,R_1^{-1}\left[d\scrY_t-H\left(\frac{\calX_t+
\overline{\eta}_t(\iota)}{2}\right)dt\right] \\
(\texttt{F3})\qquad  d\calX_t~=&~ A\,\calX_t~dt~+~R\,\calP_{\overline{\eta}_t}^{-1}\left(\calX_t-\overline{\eta}_t(\iota)\right)dt+\calP_{\overline{\eta}_t}\,H^\prime\,R_1^{-1}\left[d\scrY_t-H\left(\frac{\calX_t+
\overline{\eta}_t(\iota)}{2}\right)dt\right]\nonumber
\end{align}
where
\begin{equation}\label{def-nl-cov}
	\overline{\eta}_t := \mbox{\rm Law}(\calX_t~|~\calY_t)
\end{equation}
and thus the diffusions in (\ref{Kalman-Bucy-filter-nonlinear-ref}) depend in some nonlinear fashion on the conditional law of the diffusion process itself. In all three cases  $(\calV_t,\calW_t,\calX_0)$ are independent copies of $(\scrV_t,\scrW_t,\scrX_0)$. These diffusions are time-varying Ornstein-Uhlenbeck processes \cite{DelMoral/Tugaut:2016} and consequently $\overline{\eta}_t$ is Gaussian; see also \cite{Bishop/DelMoral:2016}. These Gaussian distributions have the same conditional mean $\overline{\eta}_t(\iota)$ and conditional covariance $\calP_{\overline{\eta}_t}$. 

\begin{proposition}[\cite{DelMoral/Tugaut:2016,Bishop/DelMoral:2016}]
We have
\begin{equation}
	\overline{\eta}_t := \mbox{\rm Law}(\calX_t~|~\calY_t) \,= \, \mbox{\rm Law}(\scrX_t~|~\calY_t) =: {\eta}_t
\end{equation}
and ${X}_t:=\overline{\eta}_t(\iota)={\eta}_t(\iota)$ and $P_t=\calP_{\overline{\eta}_t}=\calP_{\eta_t}$ where $X_t$ and $P_t$ correspond to the Kalman-Bucy filter update and Riccati equations in (\ref{nonlinear-KB-mean}) and (\ref{nonlinear-KB-Riccati}).
\end{proposition}

We may refer to this specific class (\ref{Kalman-Bucy-filter-nonlinear-ref}) of McKean-Vlasov-type diffusion as a Kalman-Bucy diffusion process \cite{Bishop/DelMoral:2016}. The case (\texttt{F1}) corresponds to the limiting object that is sampled in the continuous-time version of the `vanilla' {\tt EnKF} \cite{evensen2009book}; while (\texttt{F2}) is the continuous-time limiting object that is sampled in the `deterministic' {\tt EnKF} of \cite{sakov2008deterministic}, see also \cite{Reich2013}; and (\texttt{F3}) is a fully deterministic transport-inspired equation \cite{Reich2013,Taghvaei2016ACC}. Note that in this case (\texttt{F3}) the existence of the inverse of $\calP_{\overline{\eta}_t}$ is given by the positive-definiteness properties of the solution of the Riccati equation in (\ref{nonlinear-KB-Riccati}). In the next section we detail the Monte-Carlo ensemble filters derived from these Kalman-Bucy diffusion processes. 

Note we may define a generalised version of case (\texttt{F3}) by,
\begin{align}\label{Kalman-Bucy-filter-nonlinear-ref-specialf3}
(\texttt{F3}')\qquad  d\calX_t~=&~ A\,\calX_t~dt~+~R\,\calP_{\overline{\eta}_t}^{-1}\left(\calX_t-\overline{\eta}_t(\iota)\right)dt \nonumber\\
&\qquad+\calP_{\overline{\eta}_t}\,H^\prime\,R_1^{-1}\left[d\scrY_t-H\left(\frac{\calX_t+
\overline{\eta}_t(\iota)}{2}\right)dt\right] +G_t\,\calP_{\overline{\eta}_t}^{-1}\left(\calX_t-\overline{\eta}_t(\iota)\right)dt
\end{align}
for any skew symmetric matrix $G^\prime_t=-G_t$ that may also depend $\overline{\eta}_t$. This added tuning parameter may be related to an optimality metric, when deriving this transport equation from an optimal transport beginning. We may also write similar generalised versions (\texttt{F1}$'$) and  (\texttt{F2}$'$) by adding $G_t\,\calP_{\overline{\eta}_t}^{-1}\left(\calX_t-\overline{\eta}_t(\iota)\right)$ to (\texttt{F1}) and  (\texttt{F2}); though practically it likely makes little sense.

\section{Ensemble Kalman-Bucy Filtering}

Ensemble Kalman-Bucy filters ({\tt EnKF}) coincide with the mean-field particle interpretation of the nonlinear diffusion processes defined in \eqref{Kalman-Bucy-filter-nonlinear-ref}. 

Let $(\calV^i_t,\calW^i_t,\calX_0^i)$ with ${1\leq i\leq \sfN+1}$ be $(\sfN+1)$ independent copies of $(\calV_t,\calW_t,\calX_0)$. Again, we consider three different cases of Kalman-Bucy-type interacting diffusion process,
\begin{align}\label{EnKF-sampled-diffusions}
(\texttt{F1})\qquad d\calX_t^i~=&~A\,\calX_t^i\,dt~+~R^{1/2}\,d\calV^i_t+\widehat{P}_{t}\,H^\prime\,R_1^{-1}\left[d\scrY_t-\left(H\calX_t^i\,dt+R_1^{1/2}\,
d\calW_{t}^i\right)\right]\nonumber\\
(\texttt{F2})\qquad  d\calX_t^i~=&~A\,\calX_t^i\,dt~+~R^{1/2}\,d\calV^i_t+\widehat{P}_{t}\,H^\prime\,R_1^{-1}\left[d\scrY_t-H\left(\frac{\calX^i_t+
 \widehat{X}_{t}}{2}\right)dt\right] \\
(\texttt{F3})\qquad  d\calX_t^i~=&~ A\,\calX_t^i\,dt~+~R\,{\widehat{P}_{t}}^{-1}\left(\calX_t^i- \widehat{X}_{t}\right)dt+\widehat{P}_{t}\,H^\prime\,R_1^{-1}\left[d\scrY_t-H\left(\frac{\calX^i_t+
 \widehat{X}_{t}}{2}\right)dt\right]\nonumber
\end{align}
with $1\leq i\leq \sfN+1$ and the rescaled (particle) sample mean and covariance
\begin{equation}\label{fv1-3-2}
\begin{array}{l}
\displaystyle\widehat{\eta}_{t}:={\eta}^\sfN_{t}=\frac{1}{\sfN+1}\sum_{i=1}^{\sfN+1}\delta_{\calX_t^i} \\
	~\qquad\qquad\qquad\qquad\Longrightarrow\quad~\displaystyle
\widehat{X}_{t} :=X^\sfN_t=\frac{1}{\sfN+1}\sum_{i=1}^{\sfN+1}\calX_t^i \quad\mathrm{and}\quad \displaystyle \widehat{P}_{t} := P^\sfN_t=\frac{\sfN+1}{\sfN}\,\calP_{\widehat{\eta}_{t}} 
\end{array}\end{equation}
In cases $(\texttt{F1})$ and $(\texttt{F2})$ we have $\sfN\geq1$ and in case $(\texttt{F3})$ we require $\sfN\geq d$ for the almost sure invertibility of $\widehat{P}_{t}$ (although in case $(\texttt{F3})$ one may substitute a pseudo-inverse of $\widehat{P}_{t}$ without changing the mathematical analysis). The scaling factor on the sample covariance ensures unbiasedness. A sampled version of case $(\texttt{F3}')$ may also be derived in the same way. 

The filters of (\ref{EnKF-sampled-diffusions}) are mean-field approximations of those in \eqref{Kalman-Bucy-filter-nonlinear-ref}. In (\ref{EnKF-sampled-diffusions}) we see the utility of the Kalman-Bucy filter formulation in \eqref{Kalman-Bucy-filter-nonlinear-ref}. In particular, in (\ref{EnKF-sampled-diffusions}) we have eliminated the classical Riccati matrix differential equation completely, and replaced it with an ensemble of (interacting) particle flows and the computation of a sample covariance matrix from this ensemble. The sample mean and covariance of \eqref{fv1-3-2} can also be used for inference or decision making, etc.

\subsection{Vanilla Ensemble Kalman-Bucy Filter}

 The vanilla {\tt EnKF}, denoted by {\tt VEnKF}, is associated with the first case $(\texttt{F1})$ of nonlinear process $\calX_t$ in (\ref{Kalman-Bucy-filter-nonlinear-ref}) and is defined by the Kalman-Bucy-type interacting diffusion process $(\texttt{F1})$ in (\ref{EnKF-sampled-diffusions}). We then have the following key result.

\begin{proposition}[\cite{DelMoral/Tugaut:2016}]
	Let $\sfN\geq1$. The stochastic flow of the sample mean satisfies,
 \begin{align}
 d\widehat{X}_{t} &~\stackrel{ law}{=}~ \left(A-\widehat{P}_{t}\,S\right)\widehat{X}_{t}\,dt+ \widehat{P}_{t}\,H^{\prime}\,R_1^{-1} d\scrY_t+\frac{1}{\sqrt{\sfN+1}}\left(R+\widehat{P}_{t}\,S\,\widehat{P}_{t}\right)^{1/2} d\calB_t \label{flow-sm}
\end{align}
where $\calB_t$ is an independent $d$-dimensional Brownian motion. 

The sample covariance evolves according to a so-called matrix-valued Riccati diffusion process of the form,
\begin{equation}\label{f21}
	d\widehat{P}_{t} ~\stackrel{ law}{=}~ \mathrm{Ricc}(\widehat{P}_{t})\,dt+\frac{2}{\sqrt{\sfN}}\left[{\widehat{P}_{t}}^{1/2}\,d\calM_t\,\left(R+\widehat{P}_{t}\,S\,\widehat{P}_{t}\right)^{1/2}\right]_{\mathrm{sym}}
\end{equation}
where $\calM_t$ is a $(d\times d)$-matrix with independent Brownian entries (also independent of $\calB_t$). 

\end{proposition}

We see that for the vanilla {\tt EnKF}, the convergence of $\widehat{X}_{t}\rightarrow X_{t}$ and $\widehat{P}_{t}\rightarrow P_{t}$ as $\sfN\rightarrow\infty$ follows immediately. This result follows via the martingale representation theorem, e.g. Theorem 4.2 in \cite{karatzas}, see also~\cite{doob}.

\subsection{`Deterministic' Ensemble Kalman-Bucy Filter}

The `deterministic' {\tt EnKF}, denoted {\tt DEnKF}, is associated with the second case $(\texttt{F2})$ of nonlinear process $\calX_t$ in (\ref{Kalman-Bucy-filter-nonlinear-ref}), and is defined by the Kalman-Bucy-type interacting diffusion process $(\texttt{F2})$ in (\ref{EnKF-sampled-diffusions}). This `deterministic' epithet in the {\tt DEnKF} follows because the update `part' of the particle flow is deterministic and does not rely on the stochastic perturbations by $\calW_t^i$ appearing in the {\tt VEnKF}. This name and idea was taken from \cite{sakov2008deterministic}; see also \cite{Bergemann2012,Reich2013} and \cite{tippetsqrt2003,Lange2021b}. We have the following key result.

\begin{proposition}[\cite{2017arXiv171110065B,Bishop/DelMoral/multiDimRicc}]
	Let $\sfN\geq1$. The stochastic flow of the sample mean satisfies,
 \begin{align}
 d\widehat{X}_{t} &~\stackrel{ law}{=}~ \left(A-\widehat{P}_{t}\,S\right)\widehat{X}_{t}\,dt+ \widehat{P}_{t}\,H^{\prime}\,R_1^{-1} d\scrY_t+\frac{1}{\sqrt{\sfN+1}}\,R^{1/2} d\calB_t \label{denkf-sm}
\end{align}
where $\calB_t$ is an independent $d$-dimensional Brownian motion. 

The sample covariance evolves according to a so-called matrix-valued Riccati diffusion process of the form,
\begin{equation}\label{denkf-scov}
	d\widehat{P}_{t} ~\stackrel{ law}{=}~ \mathrm{Ricc}(\widehat{P}_{t})\,dt+\frac{2}{\sqrt{\sfN}}\left[{\widehat{P}_{t}}^{1/2}\,d\calM_t\,R^{1/2}\right]_{\mathrm{sym}}
\end{equation}
where $\calM_t$ is a $(d\times d)$-matrix with independent Brownian entries (also independent of $\calB_t$). 

\end{proposition}

Again, for the {\tt DEnKF}, the convergence of $\widehat{X}_{t}\rightarrow X_{t}$ and $\widehat{P}_{t}\rightarrow P_{t}$ as $\sfN\rightarrow\infty$ follows immediately. Note the simplified diffusion weighting(s) in the case of the {\tt DEnKF}, as compared to the {\tt VEnKF}.

\subsection{Transport-Inspired Ensemble Transport Filter} 

The fully deterministic ensemble transport filter {\tt DEnTF} is associated with the third case $(\texttt{F3})$; defined by the Kalman-Bucy-type interacting diffusion process $(\texttt{F3})$ in (\ref{EnKF-sampled-diffusions}). In this case, we have the special result.

\begin{proposition}[\cite{Reich2013,Taghvaei2016ACC}]
	Let $\sfN\geq1$. The flow of sample mean is given by,
 \begin{align}
 d\widehat{X}_{t} &~{=}~ \left(A-\widehat{P}_{t}\,S\right)\widehat{X}_{t}\,dt+ \widehat{P}_{t}\,H^{\prime}\,R_1^{-1} d\scrY_t, \quad \widehat{X}_{0} :=\frac{1}{\sfN+1}\sum_{i=1}^{\sfN+1}\calX_0^i \label{dentf-sm}
\end{align}

The sample covariance evolves according to the deterministic Riccati equation,
\begin{equation}\label{dentf-scov}
	d\widehat{P}_{t} ~{=}~ \mathrm{Ricc}(\widehat{P}_{t})\,dt, \qquad \widehat{P}_{0} := \frac{\sfN+1}{\sfN}\,\calP_{\widehat{\eta}_{0}} 
\end{equation}
\end{proposition}

Note that the particle mean $\widehat{X}_{t}$ and the particle covariance $\widehat{P}_{t}$ associated with the particle interpretation $(\texttt{F3})$ discussed in (\ref{EnKF-sampled-diffusions}) satisfy exactly the equations of the Kalman-Bucy filter with the associated deterministic Riccati equation. 

The ``randomness'' in this case only comes from the initial conditions. The stability analysis of this class of {\tt DEnTF} model resumes to the one of the Kalman-Bucy filter and the associated Riccati equation. Thus, the results, e.g., in (\ref{ref-phi-1-stability}), (\ref{KB-bias-conv}), (\ref{event-control-KB-bis}) and (\ref{control-KF-error}) hold immediately; see also \cite{Bishop/DelMoral:2016} in the linear-Gaussian setting. In \cite{deWiljes2018,deWiljes2019} this filter is analysed in the case of a nonlinear signal, but fully observed (linear observation) model. The fluctuation analysis in this case can also be developed easily by combining certain stability results w.r.t. the initial state (see \cite{Bishop/DelMoral:2016}) with conventional sample estimates based on independent copies of the initial states (see e.g.~\cite{wishart-bdn-17} for estimates associated with classical sample covariance estimates). Consequently, we do not consider this class of model going forward, but recommend \cite{Bishop/DelMoral:2016,deWiljes2018,deWiljes2019}. 

When $\sfN$ is small compared to $d$, the inverse of the sample covariance defining the {\tt DEnTF} is ill-posed and this is likely a limiting factor in the applicability of this method in high-dimensional applications with stochastic state evolutions. With non-Gaussian signal noise, one may also prefer the stochastic perturbation method in the {\tt DEnKF}.

\subsection{Nonlinear Ensemble Filtering in Practice}

In practice, the ensemble Kalman filtering methodology is applied in high-dimensional, nonlinear state-space models, e.g. see \cite{evensen03,evensen2009book} and the application references listed in the introduction. 

It is rather straightforward to extend the algorithmic particle methods in (\ref{EnKF-sampled-diffusions}) to nonlinear systems as we now outline. Consider a time-invariant nonlinear diffusion model of the form,
\begin{equation}\label{nonlinear-diffusion-filtering}
\begin{split}
d\scrX_t~&=~a(\scrX_t)\,dt \,+\, R^{1/2}\,d\scrV_t \\
d\scrY_t~&=~h(\scrX_t)\,dt \,+\, R_1^{1/2}\,d\scrW_t
\end{split}
\end{equation}
where $a:\bbR^d\rightarrow\bbR^d$ and $h:\bbR^d\rightarrow\bbR^{d_y}$ are the nonlinear signal and sensor model functions of some sufficient regularity.

Let $(\calV^i_t,\calW^i_t,\calX_0^i)$ with ${1\leq i\leq \sfN+1}$ be $(\sfN+1)$ independent copies of $(\scrV_t,\scrW_t,\scrX_0)$. We consider the three {\tt EnKF} variants as before and define the flow of particles by,
\begin{align}\label{EnKF-sampled-diffusions-nonlinear}
(\texttt{NF1})\qquad d\calX_t^i~=&~a(\calX_t^i)\,dt~+~R^{1/2}\,d\calV^i_t+\widehat{P}^{\,h}_{t}\,R_1^{-1}\left[d\scrY_t-\left(h(\calX_t^i)\,dt+R_1^{1/2}\,
d\calW_{t}^i\right)\right]\nonumber\\
(\texttt{NF2})\qquad  d\calX_t^i~=&~a(\calX_t^i)\,dt~+~R^{1/2}\,d\calV^i_t+\widehat{P}^{\,h}_{t}\,R_1^{-1}\left[d\scrY_t-\left(\frac{h(\calX_t^i)+
 \widehat{h}_{t}}{2}\right)dt\right] \\
(\texttt{NF3})\qquad  d\calX_t^i~=&~ a(\calX_t^i)\,dt~+~R\,{\widehat{P}_{t}}^{-1}\left(\calX_t^i- \widehat{X}_{t}\right)dt+\widehat{P}^{\,h}_{t}\,R_1^{-1}\left[d\scrY_t-\left(\frac{h(\calX_t^i)+
 \widehat{h}_{t})}{2}\right)dt\right]\nonumber
\end{align}
with $1\leq i\leq \sfN+1$ and the (particle) sample mean $\widehat{X}_{t}$ and sample covariance $\widehat{P}_{t}$ defined as usual, e.g. see (\ref{fv1-3-2}), and with the observation function sample mean and sample cross-covariance defined as,
\begin{equation}\label{fv1-3-2-nonlinear}
\widehat{h}_{t} :=\frac{1}{\sfN+1}\sum_{i=1}^{\sfN+1}h(\calX_t^i) ~~\quad\mathrm{and}\quad~~ \displaystyle \widehat{P}^{\,h} := \frac{1}{\sfN}\,\sum_{i=1}^{\sfN+1}\left[\calX_t^i-\widehat{X}_{t}\right]\left[h(\calX_t^i)-\widehat{h}_{t}\right]^{\prime}
\end{equation}

The mean-field limit of these interacting nonlinear conditional particle diffusion systems (\ref{EnKF-sampled-diffusions-nonlinear}) is studied in \cite{deWiljes2018,Lange2020}. The (conditional) law of these mean field McKean-Vlasov diffusions may even be given in terms of a Kushner/Fokker-Planck-type partial differential equation, e.g. see \cite{deWiljes2018,Lange2020}. However, if the mean-field limit in this nonlinear setting is denoted by, say, $\calX_t$, then it is certainly true that,
\begin{equation}
	 \mbox{\rm Law}(\calX_t~|~\calY_t) \,\neq \, \mbox{\rm Law}(\scrX_t~|~\calY_t) =: {\eta}_t
\end{equation}
in the nonlinear model setting. Said differently, even with infinite computational power, the {\tt EnKF} methods as applied in this nonlinear model setting do not converge to the optimal nonlinear Bayes filter. As noted earlier, and again later, the {\tt EnKF} in this nonlinear model setting is probably best viewed in practice as a type of (random) sample-based (point-valued) state estimator or a stochastic observer. In general it should not be seen as an approximation of the optimal Bayesian filter.

We discuss connections and extensions of our results to the nonlinear model setting, including different instances of the {\tt EnKF} in these settings, in a later section (at the end of this article).

\section{Theory in the Linear-Gaussian Setting}

Going forward, we consider only the {\tt VEnKF} (case (\texttt{F1})) and {\tt DEnKF} (case (\texttt{F2})) since as noted the theory of the {\tt DEnTF} in the linear-Gaussian setting reverts to that of the standard Kalman-Bucy filter as detailed in \cite{Bishop/DelMoral:2016}. The parameter $\kappa\in\{0,1\}$ will distinguish the two cases ($\kappa=1$ in case (\texttt{F1}), and $\kappa=0$ in case (\texttt{F2})) throughout. 

We may unify the analysis via the following representation,
\begin{align} 
d\widehat{X}_t~=&~(A-\widehat{P}_tS)~\widehat{X}_t~dt+\widehat{P}_t~H^{\prime} R_1^{-1}~d\mathscr{Y}_t+ \frac{1}{\sqrt{\sfN+1}}\,\Sigma^{1/2}_{\kappa}(\widehat{P}_t)~d\calB_t \label{EnKF-3}
 \quad\\
d\widehat{P}_t~=&~\mathrm{Ricc}(\widehat{P}_t)~dt+\frac{2}{\sqrt{\sfN}}\left[
\widehat{P}_t^{1/2}~d\calM_t~\Sigma^{1/2}_{\kappa}(\widehat{P}_t)\right]_{\mathrm{sym}} \label{EnKF-1}
\end{align}
with the mapping,
\begin{equation}
 \Sigma_{\kappa}(Q)\,:=\,R+\kappa\,QSQ\qquad\mbox{\rm with}\quad \kappa=\left\{
\begin{array}{rl}
1&\mbox{\rm in case (\texttt{F1})}\\
0&\mbox{\rm in case (\texttt{F2})}
 \end{array}\right. \label{sigma-mapping}
\end{equation}

Let $\widehat{Z}_t:=(\widehat{X}_t-\mathscr{X}_t)$ and observe that
\begin{eqnarray}
d\widehat{Z}_t&=&
(A-\widehat{P}_tS)\,\widehat{Z}_t\,dt+
\widehat{P}_t~H^{\prime} R_1^{-1/2}~d\mathscr{W}_t-R^{1/2}\,d\mathscr{V}_t+ \frac{1}{\sqrt{\sfN+1}}\,\Sigma^{1/2}_{\kappa}(\widehat{P}_t)~d\calB_t \nonumber\\
&\stackrel{ law}{=}&(A-\widehat{P}_tS)\,\widehat{Z}_t\,dt+\Omega^{1/2}_{\kappa} (\widehat{P}_t)\,d\widehat{\scrB}_t  \label{EnKF-2}
\end{eqnarray}
for some independent $d$-dimensional Wiener process $\widehat{\scrB}_t$ and with,
\begin{equation}
	\Omega_{\kappa} \,:=\, \Sigma_{1} + \frac{1}{{\sfN+1}}\, \Sigma_{\kappa}
\end{equation}
Note we often refer to the flows $\widehat{Z}_t$ or ${Z}_t$ as error flows.

We also underline that
\begin{equation}
\widehat{Z}_t - {Z}_t \,=\, (\widehat{X}_t-\mathscr{X}_t)-({X}_t-\mathscr{X}_t)\,=\,\widehat{X}_t-{X}_t
\end{equation}
so that the difference between the noisy error flow $\widehat{Z}_t$ and the classical Kalman-Bucy error flow $Z_t$ is equal to the difference between the {\tt EnKF} (sample mean) state estimate and the classical Kalman-Bucy state estimate.

Let $\widehat{\phi}_t(Q):=\widehat{P}_t$ denote the flow of the Riccati diffusion equation in (\ref{EnKF-1}) with $\widehat{P}_0=Q\in\bbS^0_d$. Let $\widehat{\psi}_t(z,Q):=\widehat{Z}_t$ denote the flow of the stochastic error (\ref{EnKF-2}) with $\widehat{Z}_0=z=(x-\scrX_0)\in\bbR^d$ and $\widehat{P}_t= \widehat{\phi}_t(Q)$. Finally, we denote the flow of the sample mean in (\ref{EnKF-3}) with $\widehat{X}_0=x\in\bbR^d$ by $\widehat{\chi}_t(x,Q):=\widehat{X}_t$. 

We underline further that the difference between two error flows satisfies,
\begin{equation}
	\widehat{\psi}_{t}(z_1,Q_1) - \widehat{\psi}_{t}(z_2,Q_2) = \widehat{\chi}_t(x_1,Q_1) - \widehat{\chi}_t(x_2,Q_2) \label{diff-error-flows-sample-means}
\end{equation}
and is thus equal to the difference between the two corresponding sample means (with compatible starting points). Studying the difference between two error flows $(\widehat{\psi}_{t}(z_1,Q_1) - {\psi}_{t}(z_2,Q_2))$ subsumes the study of something like $(\widehat{\chi}_t(x_1,Q_1) - {\chi}_t(x_2,Q_2))$ which is the difference between the {\tt EnKF} (sample mean) state estimate and the classical Kalman-Bucy state estimate (with different initial conditions).

For any $s\leq t$ and $Q\in\bbS_{d}^0$ we define the stochastic state-transition matrix,
\begin{equation}\label{KB-semigroup-def-main-stoch}
\widehat{\calE}_{s,t}(Q):=\exp{\left(\oint_s^t\left(A-\widehat{\phi}_u(Q)\,S\right)du\right)} ~~~\Longleftrightarrow~~~ \partial_t \widehat{\calE}_{s,t}(Q)=\left(A-\widehat{\phi}_u(Q)S\right)\widehat{\calE}_{s,t}(Q)
\end{equation}

As with the classical Kalman-Bucy filter, e.g. see (\ref{closed-form-OU}) and (\ref{closed-form-ricc-diff}), the convergence and stability properties of the ensemble Kalman-Bucy filter and the associated Riccati diffusion equation are directly related to the contraction properties of the stochastic state-transition matrix $\widehat{\calE}_{s,t}(Q)$. For example, the flow of the stochastic error equation (\ref{EnKF-2}) is given by,
\begin{equation}\label{closed-form-OU-stoch}
\widehat{\psi}_t(z,Q) \,=\, \widehat{\calE}_{s,t}(Q)\,\widehat{\psi}_s(z,Q)+\int_s^t~\widehat{\calE}_{u,t}(Q)\,\Omega^{1/2}_{\kappa} (\widehat{\phi}_{u}(Q))\,d\widehat{\scrB}_u
\end{equation}
and the stochastic flow of the matrix Riccati diffusion (\ref{EnKF-1}) is given implicitly by
\begin{align}
	\widehat{\phi}_t(Q) ~=&~ \widehat{\calE}_{s,t}(Q)\,\widehat{\phi}_{s}(Q)\,\widehat{\calE}_{s,t}(Q)' + \int_{s}^t\,\widehat{\calE}_{u,t}(Q)\,\Sigma_{1}(\widehat{\phi}_{u}(Q))\,\widehat{\calE}_{u,t}(Q)'\,du \nonumber\\
	&~\qquad\qquad +\frac{2}{\sqrt{\sfN}}\,\int_s^t~\widehat{\calE}_{u,t}(Q)\, \left[
\widehat{\phi}^{\,1/2}_{u}(Q)~d\calM_u~\Sigma^{1/2}_{\kappa}(\widehat{\phi}_{u}(Q))\right]_{\mathrm{sym}}   \,\widehat{\calE}_{u,t}(Q)^{\prime} \label{closed-form-ricc-diff-stoch}
\end{align}
for any $s\leq t$. We denote by $\widehat{\Pi}_t$ the Markov semigroup of $\widehat{\phi}_t(Q)$ defined for any bounded measurable function $F$ on $\bbS_d$ and any $Q\in \bbS_d^0$ with the property that,
\begin{equation}
	\widehat{\Pi}_t(F)(Q) \,:=\, \bbE\left[F(\widehat{\phi}_t(Q))\right] \,=\, \int\,
\widehat{\Pi}_t(Q,dP)\,F(P)
\end{equation}
When $Q$ is random with distribution $\Gamma(dQ)$ on $\bbS_d^+$, by Fubini's theorem we have,
\begin{equation}
	(\Gamma\widehat{\Pi}_t)(F) \,:=\, \int\,\Gamma(dQ)\,\widehat{\Pi}_t(F)(Q) \,:=\, \bbE\left[F(\widehat{\phi}_t(Q))\right] \,=\, \int\left(\int\,\Gamma(dQ)~
\widehat{\Pi}_t(Q,dP)\right)\,F(P)
\end{equation}
This yields the formula
\begin{equation}
(\Gamma\widehat{\Pi}_t)(dP)=\int~\Gamma(dQ)~
\widehat{\Pi}_t(Q,dP)
\end{equation}
for the distribution of $\widehat{\phi}_t(Q)$ on $\bbS_d^+$.

We then have the first result concerning the quadratic, matrix-valued, Riccati diffusion process (\ref{closed-form-ricc-diff-stoch}).

\begin{theorem} \label{theo-existence-s-ric}
For any $\sfN\geq 1$ the Riccati diffusion (\ref{closed-form-ricc-diff-stoch}) has a unique weak solution on $\bbS^0_d$. For
 $\sfN\geq d+1$ there exists a unique strong solution on $\bbS^+_d$. Moreover, $\widehat{\Pi}_t(Q,dP)$ is a strongly Feller and irreducible semigroup with a unique invariant probability measure $\widehat{\Gamma}_{\infty}$ on $\bbS^+_d$. This measure admits a positive density with respect to the natural Lebesgue measure on $\bbS_d$.
\end{theorem}

Given the existence of a solution to the Riccati diffusion (\ref{EnKF-1}), it follows a solution for $\widehat{X}_t$ in (\ref{EnKF-3}) or a solution $\widehat{Z}_t$ in (\ref{EnKF-2}) exists and is unique. This result is proven in \cite[Theorem 2.1]{Bishop/DelMoral/multiDimRicc}.

Once the problem of existence and uniqueness is tackled, one major problem in this equation is the behavior at infinity: existence of a stationary measure and speed of convergence towards this stationary measure or even distance between two solutions starting at different points.

We will make wide use of the following two assumptions in the remainder of this article. 
\setcounter{assumption}{14}
\begin{assumption}\label{mainAssumpObs}
	The matrix $S:=H^{\prime}R_1^{-1}H$ is strictly positive-definite, i.e. $S\in\bbS_d^+$. This is a strong form of observability, and it implies classical observability as defined in (\ref{def-contr-obs}). 
\end{assumption}

\setcounter{assumption}{2}
\begin{assumption}\label{mainAssumpCon}
	The pair $(A,R^{1/2})$ is controllable, as defined in (\ref{def-contr-obs}).
\end{assumption}

Under both Assumptions \ref{mainAssumpObs} and \ref{mainAssumpCon} it follows that $P_\infty\in\bbS_d^+$ and $\mathrm{Absc}(A-P_\infty S)<0$, see the earlier discussion on this topic. We may relax the controllability Assumption \ref{mainAssumpCon} to just stabilisability. We discuss Assumption \ref{mainAssumpObs} more later as it (re-)appears throughout our presentation and is more restrictive than the classical observability/detectability assumptions in classical Kalman filtering (noting again it implies observability/detectability).

We emphasise the following:

\begin{quote}
	Suppose Assumptions \ref{mainAssumpObs} and \ref{mainAssumpCon} hold. Then there exists some logarithmic norm, which we denote by $\overline{\mu}(\cdot):\mathbb{M}_d\rightarrow\mathbb{R}$, with the property that $\overline{\mu}(A-P_\infty S)<0$.
\end{quote}

Proof of this statement follows from the fact that $\mathrm{Absc}(A-P_\infty S)<0$ under just detectability and stabilisability model conditions, and then an application of \cite[Theorem 5]{strom1975logarithmic}. The logarithmic norm $\overline{\mu}(\cdot)$ is not necessarily unique, but any particular chosen logarithmic norm $\overline{\mu}(\cdot)$ is indexed to the model parameters $(A,H,R,R_1)$. We use the notation $\overline{\mu}(\cdot)$ to distinguish the log-norms for which $\overline{\mu}(A-P_\infty S)<0$ whenever $\mathrm{Absc}(A-P_\infty S)<0$ holds, or more specifically throughout this work whenever Assumptions \ref{mainAssumpObs} and \ref{mainAssumpCon} hold. 

In prior work \cite{DelMoral/Tugaut:2016,Bishop/DelMoral/Niclas:2017,2017arXiv171110065B,Bishop/DelMoral/STV2018,Bishop/DelMoral/multiDimRicc} and even the first draft of this article, we state certain results in terms of ${\mu}(A-P_\infty S)$, and under the assumption ${\mu}(A-P_\infty S)<0$; for some, but we don't care which, logarithmic norm $\mu(\cdot)$. We knew of course that certain observability and controllability model conditions ensured $\mathrm{Absc}(A-P_\infty S)<0$. However, it was unclear that negativity of the spectral abscissa translated in general to ${\mu}(A-P_\infty S)<0$ for some version of the logarithmic norm. Thus in many results we start with the assumption ${\mu}(A-P_\infty S)<0$ in prior work \cite{DelMoral/Tugaut:2016,Bishop/DelMoral/Niclas:2017,2017arXiv171110065B,Bishop/DelMoral/STV2018,Bishop/DelMoral/multiDimRicc}, and claimed somewhat informally that this amounts to asking for a strong form of observability and controllability (given its similarity to $\mathrm{Absc}(A-P_\infty S)<0$, but without actually giving testable model conditions). Owing to \cite[Theorem 5]{strom1975logarithmic}, we can begin results simply with some form of observability and controllability assumption (typically we need the stronger observability Assumption \ref{mainAssumpObs}, for different reasons) and state results in terms of the special class of logarithmic norms $\overline{\mu}(A-P_\infty S)<0$; which we know is negative because $\mathrm{Absc}(A-P_\infty S)<0$. This is a significant relaxation of the conditions precedent in many of the subsequent results; and places these results back in the testable and relatable context of classical controllability and observability assumptions.

In Table \ref{tableflows} we denote the relevant flows and notation of interest going forward. This notation allows us to relate (for example) the flow of the approximation relative to the true object with respect to their initial conditions, e.g. fluctuation-type results: $\widehat{\chi}_t(x,Q) - {\chi}_t(x,Q)$; or (for example) the flow of two approximated objects with respect to different initial positions, e.g. stability/contraction-type results: $\widehat{\psi}_{t}(z_1,Q_1) - \widehat{\psi}_{t}(z_2,Q_2)$. 

\begin{table}[!h]
\begin{center}
    \begin{tabular}{| c | c | p{8.5cm} |}
    \hline
    True Flow & Approximated Flow & Description  \\ \hline
    ${\chi}_t$ & $\widehat{\chi}_t$ & Flow of the Kalman state estimate (the true conditional mean) (\ref{nonlinear-KB-mean}); and the sample mean (\ref{EnKF-3})  \\ \hline
    ${\phi}_t$ & $\widehat{\phi}_t$ & Flow of the Kalman state estimate error covariance, i.e. the Riccati differential equation (\ref{nonlinear-KB-Riccati}); and the sample covariance (\ref{EnKF-1}) \\ \hline
    ${\psi}_t$ & $\widehat{\psi}_t$ & Flow of the Kalman state error (\ref{kf-error1}); and the sample error (\ref{EnKF-2})  \\ \hline
     ${\calE}_{s,t}$ & $\widehat{\calE}_{s,t}$ & The state transition matrix (\ref{KB-semigroup-def-main}); and the stochastic, approximated, state transition matrix (\ref{KB-semigroup-def-main-stoch})  \\ \hline
     N.A. & $\widehat{\Pi}_{t}$ & The Markov transition kernel for $\widehat{\phi}_t$ on $\bbS_{d}^0$. \\ \hline
    \end{tabular} 
\end{center}
\caption{Table of flow notations.}\label{tableflows}
\end{table}

In Figure \ref{fig:flowchart} we plot the flow of some of the subsequent sections and the main results. The presentation ordering is given mostly in terms of the dependencies and natural progression of the derivations. We discuss briefly the dependencies and reasoning as we progress.

\begin{figure}[!ht]
	\centering
    	\includegraphics[width=0.9\textwidth]{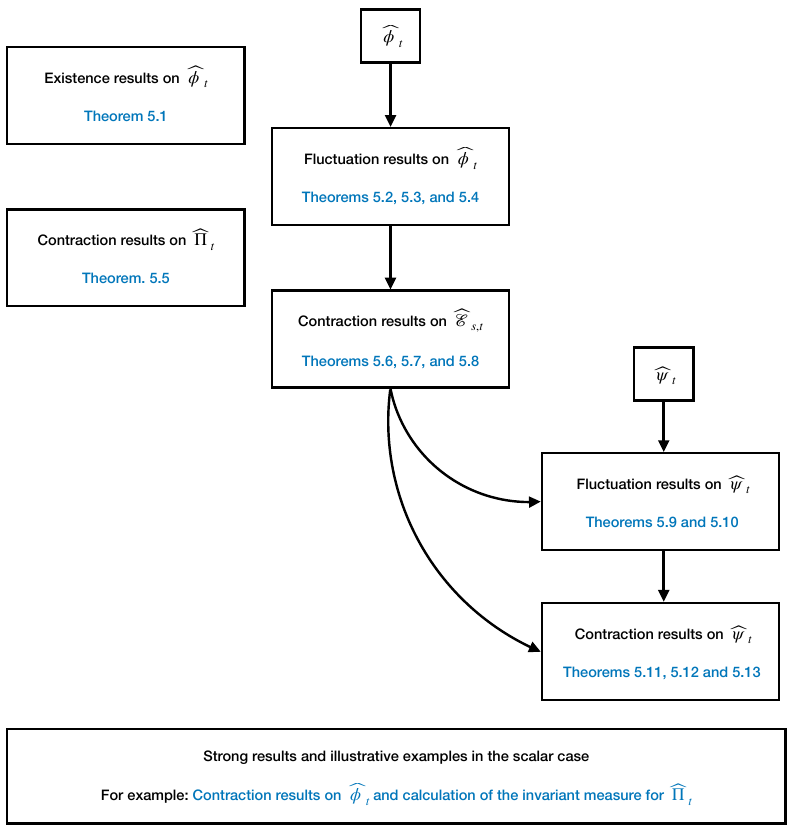}
	\caption{Flowchart of the general result and topics in this article. Although detailed proof in most cases is not given, the arrows and top-down direction in the flowchart depict both the presentation of the results in this article, and the dependency in terms of the proof and derivation of the results }
	\label{fig:flowchart}
\end{figure}

 \subsection{Fluctuation and Contraction Results for the Riccati Diffusion}

 \subsubsection{Fluctuation Properties of the Riccati Diffusion}

In this section we consider the fluctuation of $\widehat{\phi}_t(Q)$ about ${\phi}_t(Q)$ and of $\widehat{\psi}_t(z,Q)$ about ${\psi}_t(z,Q)$.

The fluctuation properties and moment boundedness properties of $\widehat{\phi}_t(Q)$ and $\widehat{\psi}_t(z,Q)$ depend naturally on the size on the fluctuation as determined by $\sfN$. 

Typically, we will write either of the the following expressions in stating our results,
\begin{equation}
''\sfN~\mbox{\rm is sufficiently large}''~\qquad\mathrm{or}\qquad~''\sfN\geq1'' \label{fluc-param}
\end{equation}
In case (\texttt{F1}) with $\kappa=1$ there is often a minimum threshold on $\sfN$ needed to prove the results. In case (\texttt{F1}) this lower threshold on $\sfN$ may be large. In case (\texttt{F2}) with $\kappa=0$, these same results typically hold; but moreover, we can often refine the relevant results and at the same time relax the conditions on $\sfN$, often needing just $\sfN\geq1$. This is a significant analytical advantage of the {\tt DEnKF} over the {\tt VEnKF}. In some cases, this advantage is practically realised and provable (and not just a by-product of analysis methods). For example, we will show later that some moments of the {\tt VEnKF} sample covariance in one-dimension provably do not exist in the steady-state without a sufficient number of particles; whereas in the {\tt DEnKF} these moments always exist with $\sfN\geq1$. In some cases, the results stated in this work are only known for the {\tt DEnKF}. If we do not specify a particular case, or a value for $\kappa\in\{0,1\}$, then the stated results may be assumed to hold for both the {\tt VEnKF} and the {\tt DEnKF}.

We start with the following under-bias estimate on the sample covariance which holds for both the {\tt VEnKF} and the {\tt DEnKF}.

\begin{theorem}\label{underbiastheorem}
For any $t\geq 0$, any $Q\in \bbS_d^0$, and any $\sfN\geq1$, we have the uniform under-bias estimate,
\begin{equation}\label{control-intro-bias}
\bbE\left[\widehat{\phi}_t(Q)\right]\,\leq~ \phi_{t}\left(Q\right)\,\leq~ c\,(1+\Vert Q\Vert)\,I
\end{equation}
for a finite constant $c>0$ that doesn't depend on the time horizon. 
\end{theorem}

We may refine this under-bias result as is done in \cite{Bishop/DelMoral/multiDimRicc}. For example, if we assume further that $S\in\bbS_d^+$, i.e. under Assumption \ref{mainAssumpObs}, then for any $t\geq 0$ we also have the refined bias estimates,
\begin{equation}\label{ref-vp-max}
0~\leq\, \phi_{t}\left(Q\right)- \bbE\left[\widehat{\phi}_t(Q)\right]\,\leq\, c(Q)\,\frac{1}{N}\, I 
\end{equation}
when $\sfN$ is sufficiently large in case (\texttt{F1}), $\kappa=1$; or for any $\sfN\geq1$ in case (\texttt{F2}), $\kappa=0$. The proof of this refinement, and details on the constant $c(Q)$, is in \cite[Theorem 2.3]{Bishop/DelMoral/multiDimRicc} and in \cite{Bishop/DelMoral/Niclas:2017}.  

We will see subsequently that Assumption \ref{mainAssumpObs}, i.e. the condition $S\in\bbS_d^+$, ensures that for any $n\geq 1$, the $n$-th moments of the trace of the sample covariance are uniformly bounded w.r.t. the time horizon (with a sufficient number of particles) even when {\em the matrix $A$ is unstable}. 

The next theorem concerns these time-uniform moment estimates on the stochastic Riccati flow in (\ref{f21}), i.e. on the flow of the sample covariance matrix. 

\begin{theorem}\label{theo-existence-s-ric-proof-bis}
Suppose Assumption \ref{mainAssumpObs} holds. For any $n\geq 1$, $t\geq 0$, any $Q\in \bbS_d^0$, and any $\sfN$ sufficiently large, we have the uniform estimate,
\begin{equation}\label{trace-Phi-inverse}
\bbE\left[\Vert \widehat{\phi}_t(Q) \Vert^n\right]^{1/n}
\,\leq\, c_{n}\,(1+\Vert Q\Vert)
\end{equation}
Furthermore, for any time horizon $t\geq \tau>0$ we also have the uniform estimates
\begin{equation}\label{trace-Phi-inverse-bis}
\bbE\left[\Vert \widehat{\phi}_t(Q) \Vert^n\right]^{1/n}\,\leq\, c_{n,\tau}~
\end{equation}
In addition, in case (\texttt{F2}), for any $\sfN\geq1$, any $n\geq 1$, $Q\in \bbS_d^0$, $t\geq 0$ and any $s\geq \tau>0$  we have the refined estimates,
\begin{equation}\label{ref-trace-unif-intro}
\bbE\left[\Vert \widehat{\phi}_t(Q) \Vert^n\right]^{1/n} \,\leq\, c\,(1+ \Vert Q\Vert)\,(1+\sqrt{\frac{n}{\sfN}})
\quad\mbox{and}\quad
\bbE\left[\Vert \widehat{\phi}_s(Q) \Vert^n\right]^{1/n} \,\leq\, c_{\tau}\,(1+\sqrt{\frac{n}{\sfN}})
\end{equation}

\end{theorem}

The proof of this result is provided in \cite[Theorem 2.2]{Bishop/DelMoral/multiDimRicc} where a precise description of the (finite) parameters $c_{n},c_{n,\tau},c,c_{\tau}>0$ is also provided. The first estimate in (\ref{trace-Phi-inverse}) also holds without Assumption \ref{mainAssumpObs}, and even if $S=0$, when $\mathrm{Absc}(A)<0$. The proof of this Theorem is based on a reduction of (\ref{f21}) to a scalar Riccati diffusion, a novel representation of its $n$-th powers, and a comparison of its moments to a judiciously designed deterministic scalar Riccati equation. We discuss this proof later, but this scalar reduction necessitates the condition $S\in\bbS_d^+$, i.e. Assumption \ref{mainAssumpObs}. The proof is conservative by nature (due to the scalar reduction and comparison).

Now we turn to quantifying the fluctuations of the matrix Riccati diffusions around their limiting (deterministic) values as found when $\sfN$ tends to $\infty$. That is, we quantify the fluctuation of the {\tt EnKF} sample covariance about the limiting covariance of the classical Kalman-Bucy filter.

\begin{theorem}\label{theo-3-intro}
Suppose Assumption \ref{mainAssumpObs} holds. For any $n\geq 1$, $t\geq 0$, any $Q\in \bbS_d^0$, and any $\sfN$ sufficiently large we have the uniform estimates,
\begin{equation}
\bbE\left[\Vert \widehat{\phi}_t(Q)-\phi_{t}(Q) \Vert^n\right]^{1/n} \,\leq\, c_n\,\frac{1}{\sqrt{N}}\,(1+ \Vert Q\Vert^7)
\label{ref-phi-1-est}
\end{equation}
In case (\texttt{F2}), for any $\sfN\geq1$, any $n\geq 1$, $t\geq 0$, and any $Q\in \bbS_d^0$, we have
 \begin{equation}\label{ref-V-0-1}
\bbE\left[\Vert \widehat{\phi}_t(Q)-\phi_{t}(Q) \Vert^n\right]^{1/n}
\,\leq\, c\, \frac{1}{\sqrt{\sfN}}\,(1+ \Vert Q\Vert^5)\,\left(1+\sqrt{\frac{n}{\sfN}}\right)^5
\end{equation}
\end{theorem}

The estimates in Theorem \ref{theo-3-intro} do not depend on $Q\in \bbS_d^0$ when $t\geq \tau$ for any $\tau>0$ and with $c_n,c$ replaced with $c_{n,\tau},c_\tau$; e.g. similarly to (\ref{trace-Phi-inverse-bis}) in Theorem \ref{theo-existence-s-ric-proof-bis}.  

The proof of the preceding Theorem is provided in \cite[Theorem 2.3]{Bishop/DelMoral/multiDimRicc} and in \cite{Bishop/DelMoral/Niclas:2017}. The proof follows from a second-order expansion of the stochastic flow $\widehat{\phi}_t$ about the deterministic flow $\phi_t$ and then an appropriate bounding of the first and second order stochastic terms. More generally, in \cite{Bishop/DelMoral/Niclas:2017} we consider a Taylor-type perturbation expansion of the form,
\begin{equation}\label{intro-main-Taylor}
	\widehat{\phi}_t=\phi_t+\sum_{1\leq k<n}~\frac{{{\sfN}}^{-k/2}}{k!}\,{\varphi}^{(k)}_t+\frac{1}{{\sfN}^{n/2}}\widehat{{\varphi}}^{\,(n)}_t
\end{equation}
for any $n\geq 1$, and a stochastic flow ${\varphi}^{(k)}_t$ {\em whose values don't depend on the ensemble size $\sfN$}, and a stochastic remainder term $\widehat{{\varphi}}^{\,(n)}_t$. Odd order stochastic terms ${\varphi}^{(k)}_t$, with $k$ odd, are zero mean (i.e. centred). This representation allows us in \cite{Bishop/DelMoral/Niclas:2017} to present sharp and non-asymptotic expansions of the matrix moments of the matrix Riccati diffusion with respect to $\sfN$.

In \cite{Bishop/DelMoral/Niclas:2017} we provide uniform estimates of the stochastic flow ${\varphi}^{(k)}_t$ w.r.t. the time horizon {\em even when the matrix $A$ is unstable}. These estimates are stronger than the conventional functional central limit theorems for stochastic processes. For example, these results imply the almost sure central limit theorem on the sample covariance,
\begin{equation}
	\sqrt{\sfN}\left[\widehat{\phi}_t-\phi_t\right]~\longrightarrow_{\sfN\rightarrow \infty}~{\varphi}_t
\end{equation}
Bias and variance estimates based on the expansion (\ref{intro-main-Taylor}) are also given in \cite{Bishop/DelMoral/Niclas:2017}. See also in particular \cite[Section 1.3]{Bishop/DelMoral/Niclas:2017} for detailed exposition of this functional central limit theorem and the bias and variance estimates. In the scalar case, we explore this expansion (\ref{intro-main-Taylor}) up to second-order in detail in a later section to illustrate this form.

The under bias result (\ref{control-intro-bias}) holds with any $\sfN\geq1$ in both the {\tt VEnKF} of case (\texttt{F1}), and in the {\tt DEnKF} of case (\texttt{F2}). This under-bias is a motivation for so-called sample covariance regularisation in practice; e.g. so-called sample covariance inflation or localisation methods \cite{andersonanderson1999,houtekamer2001seq,hamill01,Mitchell2002,evensen03}. Later we discuss the effects of inflation in particular. 

As with the deterministic Riccati equation, we may bound the moments of the inverse of the stochastic Riccati flow $\widehat{\phi}_t(Q)$ under stronger conditions on the number of particles $\sfN$ required; e.g. see \cite{Bishop/DelMoral/multiDimRicc}. It follows that with $Q\in \bbS_d^+$ and with additional conditions on $\sfN$, that for $t\geq \tau>0$ there exists a uniform positive definite lower bound on $\bbE[\widehat{\phi}_t(Q)]$.

A number of basic corollaries follow the proofs in \cite{Bishop/DelMoral/Niclas:2017,Bishop/DelMoral/multiDimRicc}, for instance, we have the monotone property,
\begin{equation}
	\bbS_d^0\,\ni\, Q_1 \,\leq\, Q_2 \quad\Longrightarrow\quad \bbE\left[\widehat{\phi}_t(Q_1)\right] \,\leq\, \bbE\left[\widehat{\phi}_t(Q_2)\right]
\end{equation}
and, for any $Q\in \bbS_d^0$, the fixed upper bound,
\begin{equation}
	\bbE\left[\widehat{\phi}_t(Q)\right]\,\leq\, {P}_\infty+\mathcal{E}_{t}({P}_\infty)\,(Q-{P}_\infty)\,\mathcal{E}_{t}({P}_\infty)^{\prime}
\end{equation}
These estimates hold for any $\sfN\geq1$ without any additional assumptions, as in Theorem \ref{underbiastheorem}. 

Several spectral estimates can be deduced from the estimates (\ref{ref-vp-max}), (\ref{ref-phi-1-est}) and (\ref{ref-V-0-1}). For example, in case (\texttt{F2}), with $\kappa=0$ and $\sfN\geq 1$ then combining (\ref{ref-V-0-1}) with the $n$-version of the Hoffman-Wielandt inequality we have the uniform estimate,
\begin{equation}
	\sup_{1\leq i\leq r}\, \bbE\left[ \left \|\lambda_i\left( \widehat{\phi}_t(Q)\right)-\lambda_i\left(\phi_t(Q)\right)\right\|^n\right]^{1/n} \,\leq\, c_n(Q)\, \frac{1}{\sqrt{N}}
\end{equation}

Finally, it is worth noting briefly that all moment boundedness and fluctuation results stated in this section hold with any $\sfN\geq 1$ and without further assumptions, if one replaces the constants $c,c_{n}, c_{\tau},c_{n}(Q),\ldots$ with functions that now depend on (and grow with) the time horizon $t\geq0$. However, if these bounds depend exponentially on time (as is quite typical in analysis), an exponent of the form $(\alpha\,t)>200$ induces an exceedingly pessimistic estimate larger than the estimated number of elementary particles of matter in the visible universe. In this sense, non-time-uniform bounds of this form are clearly impractical from a numerical user-case perspective.

 \subsubsection{Contraction and Long Time Properties of the Riccati Diffusion}

With $Q\in \bbS_d^+$ we set $\Lambda(Q):=\Vert Q\Vert_2+\Vert Q^{-1}\Vert_2$ and we consider the collection of 
$\Lambda$-norms on the set of probability measures $\Gamma_1,\Gamma_2$ on $\bbS_d^+$, indexed by $\hbar>0$, and defined by,
\begin{equation}
	 \Vert \Gamma_1-\Gamma_2\Vert_{\hbar,\Lambda} \,:=\, \sup{\vert \Gamma_1(F)-\Gamma_2(F)\vert}
\end{equation}
In the above display, the supremum is taken over all measurable function $F$ on $\bbS_d$ such that
\begin{equation}
\Vert F\Vert_{\Lambda}:=\sup_{Q\in \bbS_d^+}\frac{\vert F(Q)\vert}{1+\hbar\, \Lambda(Q)}\leq 1
\end{equation}
  
 It is known that the deterministic Riccati equation that describes the flow of the covariance matrix in classical Kalman-Bucy filtering tends to a fixed point $P_\infty$ for any initial point $Q\in\bbS_d^0$ when the (time-invariant) model (\ref{lin-Gaussian-diffusion-filtering}) is detectable and stabilisable; e.g. see (\ref{ref-phi-1-stability}) and \cite{Bishop/DelMoral:2016}. The next result is the analogue of this idea in the {\tt EnKF} setting and describes the stability of the flow of the sample covariance.
 
  \begin{theorem}\label{theo-stab-intro}
  Assume the fluctuation parameter $\sfN$ is sufficiently large such that $\bbE[\|\widehat{\phi}_t(Q)\|]$ and $\bbE[\|\widehat{\phi}^{-1}_t(Q)\|]$ are uniformly bounded (e.g. as in Theorem \ref{theo-existence-s-ric-proof-bis} for bounds on $\bbE[\|\widehat{\phi}_t(Q)\|]$). Then, there exists some finite constants $c, \alpha,\hbar>0$ such that for any $t\geq 0$ and probability measures $\Gamma_1,\Gamma_2$ on $\bbS_d^+$, we have the $\Lambda$-norm contraction inequality
 \begin{equation}\label{expo-decays}
 \Vert \Gamma_1\,  \widehat{\Pi}_{t}-\Gamma_2\,  \widehat{\Pi}_{t}\Vert_{\hbar,\Lambda} ~\leq~ c\,e^{-\alpha\, t}\,\,\Vert \Gamma_1-\Gamma_2\Vert_{ \hbar,\Lambda}
 \end{equation}
 \end{theorem}

Of course, setting $\Gamma_2=\widehat{\Gamma}_{\infty}$ where $\widehat{\Gamma}_{\infty}$ is the unique invariant probability measure described in Theorem \ref{theo-existence-s-ric} implies that for any initial probability measure $Q\sim\Gamma$ on $\bbS^+_d$ we have that $\widehat{\phi}_t(Q)$ tends to be distributed according to $\widehat{\Gamma}_{\infty}$. The proof of the above theorem is provided in \cite[Theorem 2.4]{Bishop/DelMoral/multiDimRicc} and is based on matrix-valued Lyapunov and minorisation conditions (choosing the Lyapunov candidate, $\Lambda(\cdot)$). 
 
 For one-dimensional models, the article~\cite{2017arXiv171110065B} provides explicit analytical expressions for the reversible measure of $\widehat{P}_t$ in terms of the model parameters. As expected, heavy tailed reversible measures arise when $\kappa=1$, and weighted Gaussian distributions when $\kappa=0$. The article \cite{2017arXiv171110065B} also provides sharp exponential decay rates to equilibrium, in the sense that the decay rates tend to those of the limiting deterministic Riccati equation when $\sfN$ tends to $\infty$.

In a later section, we explore the one-dimensional case in more detail and explicitly examine the invariant measures in each model $\kappa\in\{0,1\}$. The contrast between the steady-state invariant measures in each case $\kappa\in\{0,1\}$ provides some insight into various phenomenon seen in practice we believe, e.g. so-called catastrophic divergence, and fluctuations of the sample covariance, etc. We also state the strong $\bbL_n$-type contraction of $\widehat{\phi}_t(Q)$ in both cases (\texttt{F1}) and (\texttt{F2}).

\subsection{Contraction Properties of Exponential Semigroups}\label{sec-stability}

Recall that the stability properties of the deterministic ($\sfN=\infty$) semigroups $\mathcal{E}_{s,t}(Q)$ associated with the classical Kalman-Bucy filter are rather well understood; e.g. see (\ref{ref-E-1}), (\ref{ref-krause-inf}), and (\ref{expo-Estable}) and also \cite{Bishop/DelMoral:2016,bd-CARE}. We emphasise that in the deterministic case, stability of the matrix-valued Riccati differential equation, e.g. as in (\ref{ref-phi-1-stability}), follows from the contraction properties of $\mathcal{E}_{s,t}(Q)$ in (\ref{ref-E-1}); see \cite{Bishop/DelMoral:2016,bd-CARE} for the derivation. Some intuition for this follows from the implicit form for the solution in (\ref{closed-form-ricc-diff}). Similarly, in classical Kalman-Bucy filter, the stability properties of the error flow (\ref{kf-error1}) are related to the contraction properties of the state-transition matrix $\calE_{s,t}(Q)$. Again, the intuition follows from the solution form in (\ref{closed-form-OU}). The stability properties of the classical Kalman-Bucy error flow are given in, e.g., (\ref{KB-bias-conv}) and (\ref{control-KF-error}); see \cite{Bishop/DelMoral:2016}.

We come now to the contractive properties of $\widehat{\mathcal{E}}_{s,t}(Q)$ defined in (\ref{KB-semigroup-def-main-stoch}). The stability of $\widehat{\mathcal{E}}_{s,t}(Q)$ will naturally play a role in the derivation of contraction results on, e.g., the sample error flow $\widehat{\psi}_t(z,Q)$, see (\ref{closed-form-OU-stoch}). Indeed, we also require stability of $\widehat{\mathcal{E}}_{s,t}(Q)$ to derive fluctuation results on the sample error flow $\widehat{\psi}_t(z,Q)$. Note we did not need stability of the exponential semigroup to derive fluctuation results on the sample covariance $\widehat{\phi}_t(z,Q)$ earlier.

Firstly, we remark that if $S\in\bbS_d^+$, then up to a change of basis we can always assume that $S=I$. Then, for any $s,t\in[0,\infty[$ we immediately have the rather crude almost sure estimate
\begin{equation}\label{change-basis-intro-formula}
\mu\left(A\right)<0 \qquad\Longrightarrow\qquad
\left\Vert \widehat{\mathcal{E}}_{s,s+t}(Q)\right\Vert_2 ~\leq~ e^{\,t\,\mu(A)}~\longrightarrow_{t\rightarrow \infty}~0
\end{equation}
for any logarithmic norm. Note again that if $\mathrm{Absc}(A)<0$, then $\mu\left(A\right)<0$ for some log-norm. In any case, in general, asking for $A$ to be stable is a very strong and restrictive condition. We typically seek contraction results on $\widehat{\mathcal{E}}_{s,t}(Q)$ that accomodate arbitrary $A\in\bbM_d$ matrices; in particular, we seek to accommodate unstable signal matrices $A$, i.e. matrices with (some) non-negative eigenvalues. To this end, fix $Q\in \bbS^0_d$ and consider the process $\widehat{\calA}$ defined by
\begin{equation}\label{def-Aa}
 \widehat{\calA}\,:\,t\in[0,\infty[\,~~\mapsto\,~~ \widehat{\calA}_t\,:=\,A-\widehat{\phi}_t(Q)S
\end{equation}
We write $\calA$ for the analogous process driven by $\phi_t(Q)$, i.e. with $\sfN=\infty$; which we know under just detectability conditions is a time-varying stabilising matrix process \cite{VanHandel2009}.

We seek to characterise, in a useful manner, the fluctuation of the stochastic process $ \widehat{\calA}$ about $\calA$; with the hope that the contractive properties of $\widehat{\mathcal{E}}_{s,t}$ can then be in some sense related to the established contractive properties of ${\mathcal{E}}_{s,t}$.

For example, given Assumption \ref{mainAssumpObs} and $\kappa=0$, combining (\ref{ref-trace-unif-intro}) (\ref{ref-V-0-1}) and (\ref{ref-krause-inf}) with Krause's inequality \cite{Krause1994} for any $nd\geq 1$ we have the uniform fluctuation estimate,
\begin{equation}\label{krause-ref-c}
 \bbE\left[\left\Vert\,\mathrm{SpecDist}\left(\calA_t,\widehat{\calA}_t\right)\right\Vert^{nd}\right]^{1/{(nd)}}~\leq~ c_{n}(Q)\frac{1}{\sqrt{\sfN}}
\end{equation}
where we define the optimal matching distance between the spectrum of matrices $A,B\in\bbM_d$ by
\begin{equation}\label{optimal-match-d}
\mathrm{SpecDist}\left(A,B\right)=\min_{\mathrm{perm(\cdot)}}\,{\max_{1\leq i\leq d}\vert \lambda_i(A)-\lambda_{\mathrm{perm}(i)}(B)\vert}
 \end{equation}
where the minimum is taken over the set of $d!$ permutations of $\{1,\ldots,d\}$. This spectral estimate is of interest on its own, but is not immediately usable for controlling the contraction properties of the exponential semigroups.

By Theorem~\ref{theo-existence-s-ric-proof-bis} and Theorem~\ref{theo-3-intro}, under Assumption \ref{mainAssumpObs}, the collection of processes $(\calA,\widehat{\calA})$ satisfy the following regularity properties:
\begin{itemize}
\item {\em Case $\kappa\in\{1,0\}$}: 
For any $n\geq 1$, $t\geq0$, $Q\in \bbS_d^0$, and any $\sfN$ sufficiently large we have the uniform estimates
\begin{equation}\label{fluc-resultsAall}
\sqrt{\sfN}\, \bbE\left[\left\Vert\calA_t-\widehat{\calA}_t \right\Vert^n\right]^{\frac{1}{n}} \,\leq\, c_n\,(1+\|Q\|^7)~\quad\mbox{\rm and}\quad~
\bbE\left[\left\Vert\widehat{\calA}_t \right\Vert^n\right]^{\frac{1}{n}} \,\leq\, c_n\,(1+ \|Q\|)
\end{equation}
\item {\em Case $\kappa=0$}:  

For any $n\geq 1$, $t\geq0$, $Q\in \bbS_d^0$, and any $\sfN\geq1$ we have the uniform estimates
\begin{equation}\label{fluc-resultsA}
\begin{split}
\sqrt{\sfN}\,\bbE\left[\left\Vert \calA_t-\widehat{\calA}_t \right\Vert^n\right]^{\frac{1}{n}} \,&\leq\, c\,(1+ \|Q\|^5)\,(1+\frac{\sqrt{n}}{\sqrt{\sfN}})^5,\\~&\mathrm{and}~\\
\sqrt{\sfN} \,\bbE\left[\left\Vert \widehat{\calA}_t \right\Vert^n\right]^{\frac{1}{n}} \,&\leq\, c\,(1+ \|Q\|)\,(1+\sqrt{n})
\end{split}
\end{equation}
\end{itemize}
The stability properties of stochastic semigroups associated with a general collection of stochastic flows $(\calA,\widehat{\calA})$ satisfying fluctuation and moment boundedness properties in a general form accommodating both (\ref{fluc-resultsAall}) and (\ref{fluc-resultsA}) have been developed in our prior work~\cite{Bishop/DelMoral/STV2018}. Several local-type contraction estimates can now be derived. 

\begin{theorem}\label{cor-H1-semigroups}
 Let $\kappa\in\{1,0\}$ and suppose Assumptions \ref{mainAssumpObs} and \ref{mainAssumpCon} hold. Then for any increasing sequence $0\leq s \leq t_k\uparrow_{k\rightarrow\infty}\infty$, and for any $Q\in \bbS_d^0$ , the probability of the following event
 \begin{equation}\label{EA-epsilon-H1}
\limsup_{k\rightarrow\infty}\frac{1}{t_k}\log{\Vert \widehat{\mathcal{E}}_{s,t_k}(Q)\Vert} \,\leq \, \frac{1}{2}\,\overline{\mu}(A-{P}_{\infty}S)\,<0\,\quad \mbox{is greater than $1-\nu$}
 \end{equation}
for any $\nu\in]0,1[$, as soon as $\sfN$ is sufficiently large (as a function of $\nu\in]0,1[$).
\end{theorem}

This log-Lyapunov estimate (\ref{EA-epsilon-H1}) immediately implies the semigroup $\widehat{\mathcal{E}}_{s,t_k}(Q)$ is exponentially contracting with a high probability (in both cases $\kappa\in\{1,0\}$); given a sufficient number of particles, and the observability and controllability Assumptions \ref{mainAssumpObs} and \ref{mainAssumpCon}.

A number of reformulations of this result that offer insight individually are worth stating:
\begin{itemize}
\item Let $\kappa\in\{1,0\}$. For any $0\leq s \leq t_{k_1}\uparrow_{{k_1}\rightarrow\infty}\infty$, there exists a sequence $\sfN:=\sfN_{k_2}\uparrow_{{k_2}\rightarrow\infty} \infty$ such that we have the almost sure Lyapunov estimate
 \begin{equation}\label{EA-epsilon-H1-dd}
\limsup_{{k_2}\rightarrow\infty}\limsup_{{k_1}\rightarrow\infty}\frac{1}{t_{k_1}}\,\log{\Vert \widehat{\mathcal{E}}_{s,s+t_{k_1}}(Q))\Vert}\,\leq\, \frac{1}{2}\,\overline{\mu}(A-{P}_{\infty}S)
 \end{equation} 
\item Let $\kappa\in\{1,0\}$. Then, for any increasing sequence of times $0\leq s \leq t_k\uparrow_{k\rightarrow\infty}\infty$, the probability of the following event,   
 \begin{align}
 \left\{\begin{array}{l}
 \forall 0<\nu_2\leq 1~~~ \exists l\geq 1 ~~~\mbox{such~that}~~~ \forall k\geq l~~~\mbox{it~holds~that~} \\~
 \\
 \qquad\qquad\qquad\qquad\qquad\qquad\qquad\displaystyle\frac{1}{t_k}\log{\Vert \widehat{\mathcal{E}}_{s,t_k}(Q)\Vert} \,\leq\,   \frac{1}{2}\,(1-\nu_2)\,\overline{\mu}(A-{P}_{\infty}S)
 \end{array}\right\} \label{EA-epsilon-H1-cor}
\end{align}
is greater than $1-\nu_1$, for any $\nu_1\in]0,1[$, as soon as $\sfN$ is sufficiently large (as a function of $n\geq1$ and $\nu_1\in]0,1[$).
 
\item Let $\kappa\in\{1,0\}$. Consider any $s\geq 0$, any increasing sequence of time horizons $t_k\uparrow_{{k_1}\rightarrow\infty}\infty$, and any sequence $\sfN:=\sfN_{k_2,n}\uparrow_{{k_2}\rightarrow\infty} \infty$ such that $\sum_{{k_2}\geq 1}1/\sqrt{\sfN_{k_2,n}}<\infty$ for some $n\geq 1$. Then, we have the almost sure Lyapunov estimate, 
\begin{align}
 \left\{\begin{array}{l}
 \forall 0<\nu\leq 1~~~ \exists l_1,l_2\geq 1 ~~~\mbox{such~that}~~~ \forall k_1\geq l_1,~\forall k_2\geq l_2~~~\mbox{it~holds~that~} \\~
 \\
 \qquad\qquad\qquad\qquad\qquad\qquad\quad\displaystyle \frac{1}{t_{k_1}}\log{\Vert \widehat{\mathcal{E}}_{s,s+t_{k_1}}(Q)\Vert} \,\leq\,   \frac{1}{2}\,(1-\nu)\,\overline{\mu}(A-{P}_{\infty}S)
\end{array}\right\} \label{EA-epsilon-H1-cor-v2}
\end{align}
\end{itemize}

The first dot-point result captured by (\ref{EA-epsilon-H1-dd}) is derived from (\ref{EA-epsilon-H1}) in Theorem \ref{cor-H1-semigroups} via the Borel-Cantelli lemma. The next two dot-point results provide some reformulation of the supremum limit estimates (\ref{EA-epsilon-H1}) and (\ref{EA-epsilon-H1-dd}) in terms of random relaxation time horizons and random relaxation-type fluctuation parameters. The last reformulation in (\ref{EA-epsilon-H1-cor-v2}) underlines the fact that after some random time (i.e. determined by $l_1$), and given some randomly sufficiently large number of particles (determined by $l_2$) the semigroup $\widehat{\calE}_{s,t}(Q)$ is exponentially contractive. We have no direct control over the parameters $l_1$ and $l_2$ in (\ref{EA-epsilon-H1-cor-v2}) which depend on the randomness in any realisation.

Stronger results hold if we restrict $\kappa=0$, i.e. in case (\texttt{F2}). We have the following immediate corollary of our prior work in \cite{Bishop/DelMoral/STV2018} and the earlier fluctuation analysis leading to (\ref{fluc-resultsA}):

\begin{theorem}\label{cor-H2-semigroups}
 Let $\kappa=0$ and suppose Assumptions \ref{mainAssumpObs} and \ref{mainAssumpCon} hold. Then, for any $n\geq 1$, $s\geq0$, $Q\in \bbS_d^0$, there is some time horizons $\mathfrak{t}_n<\widehat{\mathfrak{t}}_n\rightarrow_{\sfN\rightarrow \infty}\infty$ such that for any $\mathfrak{t}_n\leq t\leq \widehat{\mathfrak{t}}_n$ we have
 \begin{equation}\label{EA-epsilon-H2}
 \frac{1}{t}\log{\bbE\left[\Vert \widehat{\mathcal{E}}_{s,s+t}(Q)
\Vert^{n}\right]} \,\leq\, 
 \frac{n}{4}\,\overline{\mu}(A-{P}_{\infty}S)\,<\,0
  \end{equation}
whenever $\sfN$ is sufficiently large such that $\widehat{\mathfrak{t}}_n>\mathfrak{t}_n$; see \cite{Bishop/DelMoral/STV2018} for details on these time parameters.
\end{theorem}

Importantly, in this last result we have $\widehat{\mathfrak{t}}_n\longrightarrow_{\sfN\rightarrow \infty}\infty$ and thus we can control (via $\sfN$) the horizon on which the semigroup $\widehat{\mathcal{E}}_{s,t}(Q)$ is asymptotically $\bbL_n$-stable for any $n\geq 1$ when $\kappa=0$. In other words, the estimate (\ref{EA-epsilon-H2}) ensures that the stochastic semigroup $\widehat{\mathcal{E}}_{s,t}(Q)$ is stable {\em on arbitrary long finite time horizons}, as soon as $\kappa=0$, and when the ensemble size is sufficiently large. We have the following fact immediate from Theorem \ref{cor-H2-semigroups}:

\begin{itemize}
\item Assume $\kappa=0$. For any $n\geq 1$, $s\geq0$, we have
$$
\limsup_{\sfN\rightarrow\infty}\,
\frac{1}{\widehat{\mathfrak{t}}_n}\,\log{\bbE\left[\Vert \widehat{\calE}_{s,s+\widehat{\mathfrak{t}}_n}(Q)\Vert^{n}\right]} \,\leq \,
 \frac{n}{4}\,\overline{\mu}(A-{P}_{\infty}S) 
$$
\end{itemize}

Combining Theorem \ref{cor-H1-semigroups} and Theorem \ref{cor-H2-semigroups} we may draw the basic (qualitative) conclusion that, after some initial time period, and given enough particles, the (noisy) exponential semigroups $\widehat{\mathcal{E}}_{s,t}(Q)$ are exponentially contractive (in some sense, e.g. almost-sure or $\bbL_n$-type) at a rate related to a logarithmic norm $\overline{\mu}(A-{P}_{\infty}S)$. 

We remind the reader again that weak detectability and stabilisability assumptions ensure $\mathrm{Absc}(A-{P}_{\infty}S)<0$ and consequently, via the earlier discussion and \cite[Theorem 5]{strom1975logarithmic}, there exists some logarithmic norm such that $\overline{\mu}(A-{P}_{\infty}S) <0$. Assumptions \ref{mainAssumpObs} and \ref{mainAssumpCon} imply weak detectability and stabilisability.

Finally, we also have the following new result which extends the exponential decay results for one-dimensional models presented in~\cite{2017arXiv171110065B} to the determinant of the matrix-valued Riccati diffusions considered herein. This is a type of stochastic Liouville formula.

\begin{theorem}\label{det-E-theo}
Suppose Assumptions \ref{mainAssumpObs} and \ref{mainAssumpCon} hold. Then,
for any $n\geq1$, $t\geq0$, any $Q\in \bbS_d^+$, and $\sfN$ sufficiently large we have the exponential decay estimate
\begin{align}\label{ars-0}
\bbE\left[\mbox{\rm det}(\widehat{\calE}_t(Q))^{n}\right]^{1/n} \,&=\,\bbE\left[\exp{\left(n\int_0^t\mathrm{Tr}(A-\widehat{\phi}_s(Q)S)\,ds\right)}\right]^{1/n}\nonumber \\
\,& \leq\, c_{n}(Q)\,
\exp{\left(
-t\,\sqrt{\mathrm{Tr}\left(\widehat{R}_{n}\widehat{S}_{n}\right)}
\right)}
\end{align}
with
 \begin{equation}\label{def-RS-repsilon}
\widehat{R}_{n} \,:=\, R \left( 1- \frac{1}{\sfN}(2n+ d+1)\right) \,>\, 0
\quad \mbox{and}\quad
\widehat{S}_{n} \,:=\, S \left(1- \frac{1}{\sfN}(2n+ d+1)\kappa\right)\, \,>\, 0
\end{equation} 
In addition, there exists some function $\widehat{\nu}_{n}$ with $\lim_{\sfN\rightarrow \infty}\widehat{\nu}_{n}=0$ such that
\begin{equation}\label{ars}
\bbE\left[\mbox{\rm det}(\widehat{\calE}_t(Q))^{n}\right]^{1/n}
\displaystyle\,\leq\,  c_{n}(Q)\,
\exp{\left(
-t\,(1-\widehat{\nu}_{n})\,\sqrt{\mathrm{Tr}(A)^2+\mathrm{Tr}(RS)}
\right)}
\end{equation}
\end{theorem}
The proof of this theorem is in \cite[Theorem 2.7]{Bishop/DelMoral/multiDimRicc}. In the one-dimensional case, $d=1$, this result collapses to capture the strong exponential contraction results presented in \cite{2017arXiv171110065B}. Indeed in one dimension, Theorem \ref{det-E-theo} can be seen as a significant improvement over both Theorem \ref{cor-H1-semigroups} and Theorem \ref{cor-H2-semigroups} in both theoretical development and practical usability. 

In the scalar case, strong stability results on the stochastic Riccati flow $\widehat{\phi}_t$ analogous to the deterministic setting, e.g. (\ref{ref-phi-1-stability}), also follow from Theorem \ref{det-E-theo}; see also \cite{2017arXiv171110065B} and the results and illustrative examples in a later section in this article.

 \subsection{Fluctuation and Stability of the Ensemble Kalman-Bucy Filter}

In this section we consider the fluctuation of the sample mean $(\widehat{\chi}_t(x,Q):=\widehat{X}_t$ with $x\in\mathbb{R}^d$ and $\widehat{{P}}_{0} = Q\in\bbS_d^0$; or more typically the sample mean error $\widehat{\psi}_t(z,Q) := \widehat{Z}_t =(\widehat{X}_t-\mathscr{X}_t)$ with $\widehat{Z}_0 =(x-\mathscr{X}_0)=z\in\mathbb{R}^d$. We also consider the the contraction properties of the error flow of $\widehat{\psi}_t(z,Q)$. This flow may be related to the Ornstein-Uhlenbeck process (\ref{EnKF-2}) and whose solution can be written more generally as in (\ref{closed-form-OU-stoch}). 

The first result is a fluctuation result of the ensemble sample mean about the Kalman-Bucy filter estimate, i.e. the true conditional mean; and also a conditional bias, or fluctuation, result on the conditional expectation of the ensemble sample mean given the observation sequence, with respect to the true conditional mean given by the Kalman-Bucy filter. 

The first result is given under the strong assumption that the latent signal is stable, i.e. $\mathrm{Absc}(A)<0$, and this result holds for both the {\tt VEnKF} and the {\tt DEnKF}.

\begin{theorem}\label{biastheoremmean}
Let $\kappa\in\{1,0\}$ and suppose Assumption \ref{mainAssumpObs} holds and $\mathrm{Absc}(A)<0$. For any $n\geq 1$, any $x\in\mathbb{R}^d$, any $Q\in \bbS_d^0$, and for $\sfN\geq1$ sufficiently large, we have the fluctuation estimate,
\begin{equation}\label{meanfluceqn}
	\bbE\left[ \Big\Vert\, \widehat{\chi}_t(x,Q)  \,-\, \chi_t(x,Q)\, \Big\Vert^n\right]^{1/n} \,\leq\, c_n(x,Q)\,\frac{1}{\sqrt{\sfN}}
\end{equation}
We also have the conditional bias estimate,
\begin{equation}\label{condbiasmeaneqn}
	\bbE\left[ \Big\Vert\, \bbE\left[\widehat{\chi}_t(x,Q) \,|\, \calY_t\right] \,-\, \chi_t(x,Q)\, \Big\Vert^n\right]^{1/n} \,\leq\, c_n(x,Q)\,\frac{1}{\sfN}
\end{equation}
\end{theorem}

Proof of the fluctuation estimate (\ref{meanfluceqn}) is given in \cite{DelMoral/Tugaut:2016}. Proof of the conditional bias estimate (\ref{condbiasmeaneqn}) is given in \cite[Theorem 2.4]{CrisanDM2022}. The latter result (\ref{condbiasmeaneqn}) is used in \cite{CrisanDM2022} to study the estimation of the log-normalization constant associated with a class of continuous-time filtering models.

The next theorem concerns time-uniform moment estimates on the sample mean error; and the fluctuation of the sample mean error around its limiting value (found when $\sfN$ tends to $\infty$). The next result relaxes the assumption that the latent signal be stable.

\begin{theorem}\label{theo-fluc-sample-mean}
Consider only case (\texttt{F2}) and suppose Assumptions \ref{mainAssumpObs} and \ref{mainAssumpCon} hold. Then for any $n\geq 1$, $z\in\mathbb{R}^d$, $Q\in \bbS_d^0$, there exists a time $\widehat{\mathfrak{t}}_n\rightarrow_{\sfN\rightarrow \infty}\infty$ such that for any $0\leq t\leq \widehat{\mathfrak{t}}_n$ we have,
\begin{equation}\label{trace-Psi}
\bbE\left[ \Vert \widehat{\psi}_t(z,Q) \Vert^n\right]^{1/n}
\,\leq\, c_{n}(z,Q)
\end{equation}
and
\begin{equation}
\bbE\left[ \Vert \widehat{\psi}_t(z,Q)-{\psi}_t(z,Q)\Vert^n\right]^{1/n} \,\leq\, c_n(z,Q)\,\frac{1}{\sqrt{\sfN}}
\label{ref-psi-1-est}
\end{equation}
See \cite{Bishop/DelMoral/STV2018} for details on the time parameter $\widehat{\mathfrak{t}}_n\longrightarrow_{\sfN\rightarrow \infty}\infty$.
\end{theorem}

Note again the difference $(\widehat{\psi}_{t}(z_1,Q_1) - {\psi}_{t}(z_2,Q_2))$ resumes to that of $(\widehat{\chi}_t(x_1,Q_1) - {\chi}_t(x_2,Q_2))$. Thus, e.g., (\ref{ref-psi-1-est}) is comparable to (\ref{meanfluceqn}), under different antecedent conditions.

Unlike Theorem \ref{theo-existence-s-ric-proof-bis} and Theorem \ref{theo-3-intro}, the proof of both Theorem \ref{biastheoremmean} and Theorem \ref{theo-fluc-sample-mean} requires contraction properties to be established a priori for the stochastic transition matrix $\widehat{\mathcal{E}}_{s,t}(Q)$ defined in (\ref{KB-semigroup-def-main-stoch}). Hence, in Theorem \ref{biastheoremmean} we rely on $\mathrm{Absc}(A)<0$ which ensures the contractive property holds for $\widehat{\mathcal{E}}_{s,t}(Q)$, see (\ref{change-basis-intro-formula}). In Theorem \ref{theo-fluc-sample-mean}, we rely on Theorem \ref{cor-H2-semigroups} which establishes the $\bbL_n$-contractivity of $\widehat{\mathcal{E}}_{s,t}(Q)$ without asking for $A$ to be stable, but only in the case of the {\tt DEnKF} with $\kappa=0$, at least in the multi-dimensional setting.

The proof of Theorem \ref{theo-fluc-sample-mean} is provided in \cite{2017arXiv171110065B} in the one-dimensional setting where a detailed description of the (finite) parameters $c_{n}(z,Q)>0$ are provided. The multi-dimensional result follows using similar proof methods to those used in \cite{2017arXiv171110065B} in combination with the contraction properties of the transition matrix $\widehat{\mathcal{E}}_{s,t}(Q)$ established in Theorem \ref{cor-H2-semigroups}. In the one-dimensional setting studied in \cite{2017arXiv171110065B}, contraction of $\widehat{\mathcal{E}}_{s,t}(Q)$ is given under very general model conditions which also accommodate both the  {\tt VEnKF} and the {\tt DEnKF}. Consequently, in one dimension Theorem \ref{theo-fluc-sample-mean} holds on an infinite time horizon for any $t\geq0$ and with any $\kappa\in\{0,1\}$.

One may consider a perturbation expansion of the sample mean flow as
\begin{equation}\label{intro-main-Taylor-sample-mean}
	\widehat{\psi}_t=\psi_t+\sum_{1\leq k<n}~\frac{{{\sfN}}^{-k/2}}{k!}\,{\vartheta}^{(k)}_t+\frac{1}{{\sfN}^{n/2}}\widehat{{\vartheta}}^{\,(n)}_t
\end{equation}
for any $n\geq 1$, and some stochastic flow ${\vartheta}^{(k)}_t$ {\em that does not depend on the ensemble size $\sfN$}, and some stochastic remainder term $\widehat{{\vartheta}}^{\,(n)}_t$. This implies the almost sure central limit theorem on the sample mean,
\begin{equation}
	\sqrt{\sfN}\left[\widehat{\psi}_t-\psi_t\right]~\longrightarrow_{\sfN\rightarrow \infty}~{\vartheta}_t
\end{equation}
See in particular \cite[Section 1.3]{Bishop/DelMoral/Niclas:2017} for detailed exposition of this functional central limit theorem.

Uniform propagation of chaos follows from the proceeding central limit theorems and the development in this subsection. In particular we have,
\begin{equation}
	\mathrm{Law}(\calX_t^i)~~\longrightarrow_{\sfN\rightarrow \infty}~~\mathrm{Law}(\calX_t),~~\forall t\geq0,~\forall i\in\{1,\ldots,N+1\}
\end{equation}
in some suitable metric (e.g. Wasserstein).

Now we turn to the stability of the error flow $\widehat{\psi}_t(z,Q)$ and its contraction properties. The subsequent study on the stability of $\widehat{\psi}_t(z,Q)$ relies again on the contraction of $\widehat{\calE}_{s,t}$ studied previously.

The following uniform error contraction estimate follows from (\ref{closed-form-OU-stoch}) and Theorem \ref{cor-H2-semigroups},
\begin{equation}\label{EnKF-bias-conv}
\sup_{Q\in\bbS^0_{d}}\, \left\Vert\,  \bbE\left[ \widehat{\psi}_{t}(z,Q) \,\vert\, \scrX_0\right]\, \right\Vert
~\leq~ c\,e^{t\,\alpha\,\overline{\mu}(A-{P}_{\infty}S) } \,\Vert\,x-\scrX_0 \Vert
\end{equation}
and holds for the {\tt DEnKF}, with $\kappa=0$, for some $\alpha,c>0$, and under conditions compatible with the conditions in Theorem \ref{cor-H2-semigroups}. This contraction result is analogous to (\ref{KB-bias-conv}) for the classical Kalman-Bucy filter; but under stronger conditions dictated by the available results on the contraction properties of $\widehat{\calE}_{s,t}$ stated in Theorem \ref{cor-H2-semigroups}. In particular, our methods prove this contraction (\ref{EnKF-bias-conv}) only in the case of the {\tt DEnKF}, with $\kappa=0$, with $\sfN$ sufficiently large, and on time horizons compatible with those detailed in Theorem \ref{cor-H2-semigroups}. 

If $\mathrm{Absc}(A)<0$ and Assumption \ref{mainAssumpObs} holds, then (\ref{EnKF-bias-conv}) holds on any infinite time horizon for both the {\tt VEnKF} and {\tt DEnKF}; because in this case $\widehat{\mathcal{E}}_{s,t}(Q)$ is contractive from (\ref{change-basis-intro-formula}). This is analogous to the setting of Theorem \ref{biastheoremmean}, as compared to that of Theorem \ref{theo-fluc-sample-mean}; in line with the earlier discussion on the conditions leading to stability of $\widehat{\mathcal{E}}_{s,t}(Q)$.

The next results on the stability of $\widehat{\psi}_t(z,Q)$ similarly follow immediately from those stability results in the preceding section, but are stated at the level of the process $\widehat{\psi}_t(z,Q)$ itself, rather than the stochastic exponential semigroup $\widehat{\calE}_{s,t}$.

\begin{theorem}\label{cor-H1-errorflows}
 Let $\kappa\in\{1,0\}$ and suppose Assumptions \ref{mainAssumpObs} and \ref{mainAssumpCon} hold. Then for any increasing sequence of times $t_k\uparrow_{k\rightarrow\infty}\infty$, any $z_1\not=z_2$ and any $Q\in\bbS_d^0$, the probability of the following event
 \begin{equation}\label{error-epsilon-H1-cor-v3}
\limsup_{k\rightarrow\infty}\frac{1}{t_k}\log{\Vert \widehat{\psi}_{t_k}(z_1,Q) - \widehat{\psi}_{t_k}(z_2,Q) \Vert} \,<\, \alpha\,\overline{\mu}(A-{P}_{\infty}\,S)\quad \mbox{is greater than $1-\nu$}
 \end{equation}
for any $\nu\in]0,1[$ and some $\alpha>0$, as soon as $\sfN$ is sufficiently large (as a function of $\nu$).
 \end{theorem}

Two reformulations of this result may shed insight individually and are worth highlighting:
\begin{itemize}
\item Let $\kappa\in\{1,0\}$. For $0 \leq t_{k_1}\uparrow_{{k_1}\rightarrow\infty}\infty$, there exists a sequence $\sfN:=\sfN_{k_2}\uparrow_{{k_2}\rightarrow\infty} \infty$ such that we have the almost sure Lyapunov estimate
 \begin{equation}\label{error-epsilon-H1-cor-v3-dd}
\limsup_{{k_2}\rightarrow\infty}\limsup_{{k_1}\rightarrow\infty}\frac{1}{t_{k_1}}\,\log{\Vert \widehat{\psi}_{t_{k_1}}(z_1,Q) - \widehat{\psi}_{t_{k_1}}(z_2,Q) \Vert}\,<\, \alpha\,\overline{\mu}(A-{P}_{\infty}S)
 \end{equation}

\item Let $\kappa\in\{1,0\}$. Consider any increasing sequence of time horizons $t_k\uparrow_{{k_1}\rightarrow\infty}\infty$, and any sequence $\sfN:=\sfN_{k_2,n}\uparrow_{{k_2}\rightarrow\infty} \infty$ such that $\sum_{{k_2}\geq 1}1/\sqrt{\sfN_{k_2,n}}<\infty$ for some $n\geq 1$. Then, we have the almost sure Lyapunov estimate, 
\begin{align}
 \left\{\begin{array}{l}
 \forall 0<\nu\leq 1~~~ \exists l_1,l_2\geq 1 ~~~\mbox{such~that}~~~ \forall k_1\geq l_1,~\forall k_2\geq l_2~~~\mbox{it~holds~that~} \\~
 \\
 \qquad\qquad\qquad \displaystyle \frac{1}{t_{k_1}}\log{\Vert \widehat{\psi}_{t_{k_1}}(z_1,Q) - \widehat{\psi}_{t_{k_1}}(z_2,Q)\Vert} \,\leq\,   \alpha\,(1-\nu)\,\overline{\mu}(A-{P}_{\infty}S)
\end{array}\right\} \label{error-epsilon-H1-cor-v3-ddd}
\end{align}
\end{itemize}

Again we emphasise that the reformulation in (\ref{error-epsilon-H1-cor-v3-ddd}) highlights that after some random time (i.e. determined by $l_1$), and given a random sufficiently large number of particles (determined by $l_2$) the difference of error flows (or sample means; see (\ref{diff-error-flows-sample-means})) is exponentially stable.

We have stronger $\bbL_n$-type stability results in those settings analogous to Theorems \ref{biastheoremmean} and \ref{theo-fluc-sample-mean} and in line with the discussion after Theorem \ref{theo-fluc-sample-mean} on the conditions for of $\widehat{\mathcal{E}}_{s,t}(Q)$.

\begin{theorem}\label{theo-psi-error-con1}
Let $\kappa\in\{1,0\}$ and suppose Assumption \ref{mainAssumpObs} holds and $\mathrm{Absc}(A)<0$. Then for any $n\geq 1$, any $z_1\not=z_2$, and any $Q\in \bbS_d^0$ we have the stability estimate,
\begin{equation}\label{uniformly-Lipschitz-OU}
{\bbE\left[ \Vert \widehat{\psi}_{t}(z_1,Q) - \widehat{\psi}_{t}(z_2,Q)\,\Vert^n \right]}^{1/n} \,\leq\,  c_{n}(z_1,z_2,Q)\Vert z_1-z_2\Vert\, e^{\,t\,\overline{\mu}(A-{P}_{\infty}\,S)} 
\end{equation}
whenever $\sfN$ is sufficiently large.
\end{theorem}

In the case (\texttt{F2}), i.e. for the {\tt DEnKF} only, when $\kappa=0$, we can relax the strong assumption that the latent signal be stable.

\begin{theorem}\label{cor-H2-errorflows}
 Let $\kappa=0$ and suppose Assumptions \ref{mainAssumpObs} and \ref{mainAssumpCon} hold. Then for any $n\geq 1$, any $z_1\not=z_2$, and any $Q\in \bbS_d^0$, there exists some time horizons $\mathfrak{t}_n<\widehat{\mathfrak{t}}_n\longrightarrow_{\sfN\rightarrow \infty}\infty$ such that for any $\mathfrak{t}_n\leq t\leq \widehat{\mathfrak{t}}_n$ we have the stability estimate,
\begin{equation}\label{uniformly-Lipschitz-moments-OU-Lip}
{\bbE\left[ \Vert \widehat{\psi}_{t}(z_1,Q) - \widehat{\psi}_{t}(z_2,Q)\Vert^n \right]}^{1/n} \,\leq\,  c_{n}(z_1,z_2,Q)\Vert z_1-z_2\Vert\, e^{\,t\,\overline{\mu}(A-{P}_{\infty}\,S)} 
\end{equation}
whenever $\sfN$ is sufficiently large such that $\widehat{\mathfrak{t}}_n>\mathfrak{t}_n$; see \cite{Bishop/DelMoral/STV2018} for details on these time parameters.
\end{theorem}

We emphasise again that $\widehat{\mathfrak{t}}_n\longrightarrow_{\sfN\rightarrow \infty}\infty$. With regards to qualitative reasoning, we may combine Theorem \ref{cor-H1-errorflows} and Theorem \ref{cor-H2-errorflows} and draw the basic (qualitative) conclusion that, after some initial time period, and given enough particles, the difference in (noisy) error flows $(\widehat{\psi}_{t}(z_1,Q) - \widehat{\psi}_{t}(z_2,Q))$, or the difference in sample means $(\widehat{\chi}_t(x_1,Q) - \widehat{\chi}_t(x_2,Q))$, is exponentially stable (in some sense) with a rate related to $\overline{\mu}(A-{P}_{\infty}S)$. 

 In the scalar case $d=1$, stronger stability results on the error flow $\widehat{\psi}_{t}(z,Q)$ follow from the contraction properties in Theorem \ref{det-E-theo} under weaker model and ensemble size assumptions. The strong $\bbL_n$-type stability results in the scalar $d=1$ case are quantitative and hold over infinite horizons for both the {\tt VEnKF} and the {\tt DEnKF}, i.e. with $\kappa\in\{0,1\}$, with unstable latent signals, with differing initial variance states, and with exponential rates that collapse to the optimal deterministic rates (explicitly computable when $d=1$) as $\sfN\rightarrow\infty$. See \cite{2017arXiv171110065B}; and the results, and illustrative examples in the next section.

\section{Strong Results in One-Dimensional Illustrative Examples}

Throughout this section we let $d=1$ and $R\wedge S>0$. The latter condition $R\wedge S>0$ is both necessary and sufficient for observability and controllability to hold in one dimension; and besides, in some cases, conditions on $\sfN\geq1$ no other conditions are needed in this section (and we emphasise that the latent signal may be unstable). The purpose of this section is to illustrate in more detailed quantitative terms some of the more abstract or qualitative results given in the general multi-variate setting. In some cases, the derivation of a multi-variate counterpart of a result in this section remains an open problem. In the scalar setting, the analysis of the {\tt EnKF} is rather complete in the linear-Gaussian case.

When $P_0\in[0,\infty[$, the deterministic Riccati equation defined on $[0,\infty[$, in the classical Kalman-Bucy filter, satisfies the quadratic differential equation (\ref{nonlinear-KB-Riccati}) which may be written also as,
\begin{equation}
\partial_tP_t\,=\,\mathrm{Ricc}(P_t) \,=\, -S\,(P_t-\varrho_+)\,(P_t-\varrho_-),
\label{nonlinear-KB-Riccati-modes}
\end{equation}
with the equilibrium states $(\varrho_-,\varrho_+)$ defined by
   \begin{equation}\label{def-varpi}
S\,\varrho_-:=A-\sqrt{A^2+RS}~<~0~<~S\,\varrho_+:=A+\sqrt{A^2+RS} 
 \end{equation}
 With $P_0\in[0,\infty[$, we have $P_t\rightarrow_{t\rightarrow\infty}P_\infty=\varrho_+$. It follows that,
 \begin{equation}
 	A - P_\infty\,S \,=\, -\sqrt{A^2+RS} 
 \end{equation}
 and thus simplifying, e.g. (\ref{ref-E-1}), we have the equality,
 \begin{equation}\label{ref-E-1-scalar}
 \mathcal{E}_{t}(Q) \,=\, c_t(Q)\, \mathcal{E}_{t}({P}_{\infty}) \,\leq\, c(Q)\, \mathcal{E}_{t}({P}_{\infty}) \qquad \mbox{\rm and}\qquad  \mathcal{E}_{t}({P}_{\infty})\,= \,  e^{-t\,\sqrt{A^2+RS} } 
\end{equation}
 where $-\sqrt{A^2+RS}$ may be viewed explicitly as the optimal semigroup contraction rate in the scalar case. The explicit form of the constants $c_t(Q)$, $c(Q)$ is also available in the scalar case, see \cite{2017arXiv171110065B} and also the general Floquet-type multivariate result in \cite{bd-CARE}.

 The Riccati drift function $\mathrm{Ricc}(\cdot)$ is also the derivative of the double-well potential function
\begin{equation}
	F(Q)=-\frac{S}{3}~Q~(Q-\zeta_{-})~(Q- \zeta_+) 
 \end{equation}
with the roots
\begin{equation}
\zeta_{-}:=\frac{3A}{2S}-\left[
\left(\frac{3A}{2S}\right)^2+\frac{3R}{S}\right]^{1/2}< 0< \zeta_+:=\frac{3A}{2S}+\left[
\left(\frac{3A}{2S}\right)^2+
\frac{3R}{S}\right]^{1/2}
 \end{equation}
In this situation, the general Riccati diffusion (\ref{EnKF-1}) describing the flow of the sample covariance in both case (\texttt{F1}) and case (\texttt{F2}) resumes to the Langevin-Riccati drift-type diffusion process,
\begin{equation}\label{definition-Langevin-Riccati}
	d\widehat{P}_t \,=\, \partial F(\widehat{P}_t)\,dt+
		\frac{1}{\sqrt{\sfN}}\,\widehat{P}_t^{1/2}\,\Sigma^{1/2}_{\kappa}(\widehat{P}_t) \,d\calM_t
\end{equation}
with the mapping $\Sigma_\kappa$ defined in (\ref{sigma-mapping}). Recall that case (\texttt{F1}) corresponds to the vanilla {\tt EnKF}, denoted by {\tt VEnKF}, and case (\texttt{F2}) corresponds to the `deterministic' {\tt EnKF}, denoted by {\tt DEnKF}. Also observe that $\partial F>0$ on the open interval $]0,\zeta_+[$ and $\partial F(0)=R>0=\sigma(0)$ so that the origin is repellent and instantaneously reflecting. 

At any time $t\geq0$ we may comment on the boundedness of certain moments of the sample variance and the fluctuation of the sample variance and sample mean about their limiting (classical Kalman-Bucy variance and mean) values. 

For example, we have the under-bias result $\bbE[\widehat{\phi}_t(Q)]\leq \phi_{t}$ with any $\sfN\geq1$ in both the {\tt VEnKF} of case (\texttt{F1}), and in the {\tt DEnKF} of case (\texttt{F2}). This under-bias motivates so-called variance/covariance regularisation methods in practice; e.g. so-called sample covariance inflation or localisation methods. Later we discuss the effects of inflation in particular. However, more generally, in the scalar case we have the result of Theorem \ref{theo-existence-s-ric-proof-bis} with explicit expressions on the ensemble size, i.e. we have for any $n\geq 1$, $t\geq 0$, $Q\in [0,\infty[$, and any $\sfN\geq 1\vee2\,\kappa\,(n-1)$, the uniform estimates,
\begin{equation}\label{trace-Phi-inverse-scalar}
\bbE\left[\widehat{\phi}_t(Q)^{n}\right]^{1/n}
\,\leq\, c_{n}(Q)
\end{equation}
We also have bounds on the inverse Riccati flow (leading to lower bounds on the sample covariance) under stronger conditions on $\sfN$; see \cite{2017arXiv171110065B}. We remark here, and again later when we explicitly examine the invariant measure for $\widehat{\phi}_t$, that these conditions on $\sfN$ while explicit, may still conservative (in the case of the {\tt VEnKF}). From Theorem \ref{theo-3-intro} and the scalar exposition in \cite{2017arXiv171110065B} we have the uniform fluctuation estimate $\bbE[\widehat{\phi}_t(Q)-\phi_{t}(Q)^{n}]^{1/n} \leq c_n(Q)/\sqrt{N}$ with the explicit $\sfN\geq 1\vee2\,\kappa\,(n-1)$. The constant $c_n(Q)$ is also studied in \cite{2017arXiv171110065B} with $d=1$ in explicit detail. 

Note that we may expand the stochastic flow of the sample variance as in (\ref{intro-main-Taylor}). Exploring this idea further in the scalar case for illustrative purposes, we may write the first and second-order fluctuations as,
\begin{eqnarray}
\widehat{\varphi}_t&:=&\sqrt{\sfN}\,[\widehat{\phi}_{t}-\phi_t],
\\
\widehat{{\varphi}}^{(2)}_t&:=& \sqrt{\sfN}[\widehat{\varphi}_t-{\varphi}_t]
\end{eqnarray}
where in the second line we emphasise the superscript $\cdot^{(2)}$ is an order index (not a power) and where,
\begin{eqnarray}
{\varphi}_t(Q)&:=&\int_0^t\,\left(\partial \phi_{t-s}\right)(\phi_{s}(Q))\,\Sigma^{1/2}_{\kappa}\left(\phi_{s}(Q)\right)\,d\calM_s
\end{eqnarray}
and the derivatives $\partial^k \phi_{t}$ of any order are explicitly given in \cite{2017arXiv171110065B}. In this case, $\partial \phi_{t}(Q) = \calE_t^2(Q)$ (where the superscript here is now a power). We then have,
\begin{equation}\label{scalar-Taylor}
	\widehat{\phi}_t=\phi_t+ \frac{1}{\sqrt{\sfN}}\,{\varphi}_t+\frac{1}{{\sfN}}\widehat{{\varphi}}_t^{(2)}
\end{equation}
The natural central limit theorem follows, i.e. $\sqrt{\sfN}[\widehat{\phi}_t-\phi_t]\,\longrightarrow_{\sfN\rightarrow \infty}\,{\varphi}_t$. The (non-)asymptotic variance is estimated in \cite{2017arXiv171110065B,Bishop/DelMoral/Niclas:2017}.

The expansion (\ref{scalar-Taylor}) allows ones to better understand the bias properties of the sample covariance $\widehat{\phi}_t$. Writing the third-order fluctuation as,
\begin{eqnarray}
	\widehat{\varphi}^{(3)}_t&:=&\sqrt{\sfN}\,[\widehat{\varphi}^{(2)}_t- {\varphi}^{(2)}_t/2]
\end{eqnarray}
and expanding and taking expectations,
\begin{eqnarray}\label{scalar-Taylor-expectations}
	\bbE\left[\widehat{\phi}_t\right] &=& \phi_t+ \frac{1}{2{\sfN}}\bbE\left[ {{\varphi}}_t^{(2)} \right]  + \frac{1}{{\sfN^{3/2}}}\bbE\left[ \widehat{{\varphi}}_t^{(3)} \right]
\end{eqnarray}
and limits we have the dominating ($N$-order-asymptotic) bias is given by
\begin{equation}
	\sfN \left(\bbE\left[\widehat{\phi}_t\right] - {\phi}_t\right) ~~\longrightarrow_{\sfN\rightarrow\infty}~~~~\,\frac{1}{2} \,\int_0^t\,\left(\partial ^2\phi_{t-s}\right)(\phi_{s}(Q))~\Sigma_{\kappa}\left(\phi_{s}(Q)\right)~ds~<0
\end{equation}
which is always negative (agreeing with the under-bias result $\bbE[\widehat{\phi}_t(Q)]\leq \phi_{t}$). See \cite{2017arXiv171110065B,Bishop/DelMoral/Niclas:2017} for further exploration of these general expansions. A detailed study of these expansions may aid in the development and tuning of (adaptive) sample covariance regularisation methods. 

Significantly generalising Theorem \ref{theo-fluc-sample-mean} in the scalar case \cite{2017arXiv171110065B}, we have for any $t\geq 0$, and any $\sfN>2(4n+1)(1+4\kappa)$, the uniform bound $\bbE[\widehat{\psi}_t(z,Q)^{n}]^{1/n} \leq c_{n}(z,Q)$. We also have the generalisation that for any $t\geq 0$, and any $\sfN>2(6n+1)(1+4\kappa)$, the uniform fluctuation estimate,
\begin{equation}
\bbE\left[\widehat{\psi}_t(Q)-\psi_{t}(Q)^{n}\right]^{1/n} \,\leq\, \frac{c_n(z,Q)}{\sqrt{N}}
\label{ref-psi-1-est-scalar}
\end{equation}
holds. 

The expansion (\ref{intro-main-Taylor-sample-mean}) of the sample mean (error) may be explored similarly to the above expansion of the sample covariance. The first order terms $\vartheta_t$ in (\ref{intro-main-Taylor-sample-mean}) related to the central limit theorem are studied in \cite[Section 1.3]{Bishop/DelMoral/Niclas:2017}.

The infinitesimal generator of the diffusion (\ref{definition-Langevin-Riccati}) on $]0,\infty[$ is given 
in Sturm-Liouville form by the equation
\begin{equation}
L(f)\,=\,\frac{2}{\sfN}\,\iota\,\Sigma_\kappa \,e^{V}\,\partial\left(e^{-V}\,\partial f\right)\quad\mbox{\rm with}\quad
V(\cdot)=-\frac{\sfN}{2}\int_{\delta}^\cdot~\partial F(x)\,\iota^{-1}(x)\,\Sigma^{-1}_\kappa(x)\,dx
\end{equation}
for any $\delta>0$ and where we recall the identity function $\iota(x):=x$. This implies that a reversible measure of the Riccati diffusion (\ref{EnKF-1}) in the scalar $d=1$ case is given by the formula
\begin{equation}
	\widehat{\Gamma}_{\infty}(dx)\,\propto\, 1_{]0,\infty[}(x)\,\frac{\sfN}{4\,\iota(x)\,\Sigma_{\kappa}(x)}\,\exp{\left(-V(x)\right)}\,dx
\end{equation}

In case (\texttt{F1}) corresponding to the {\tt VEnKF}, we have that $L$ is reversible w.r.t. the probability measure $\widehat{\Gamma}_\infty$ on $]0,\infty[$ defined by,
\begin{equation}\label{invariant-measure-venkf}
	\widehat{\Gamma}_\infty(dx)~\propto~1_{]0,\infty[}(x)~\exp{\left(\sfN \frac{A}{\sqrt{RS}}\,\tan^{-1}\left(x~\sqrt{\frac{S}{R}}\right)\right)} \left(\frac{x}{R+Sx^2}\right)^{\frac{\sfN}{2}}\frac{1}{x(R+Sx^2)}\,dx
\end{equation}
See also \cite{2017arXiv171110065B} for alternate derivations/forms of this heavy tailed invariant measure. The heavy tailed nature of the stationary measure implies that for the $n$-th moment to exist one requires $\sfN>0\vee 2(n-2)$. As expected this condition on $\sfN$ is generally weaker than that required for $n$-th moment boundedness at any time $t\geq0$ in (\ref{trace-Phi-inverse-scalar}) in terms of the {\tt VEnKF}. In Figure \ref{fig:momentExistEnKF} we plot the line defined by $(2n-4)/\sfN$ for various $\sfN$ values. With any $\sfN\geq1$, we have existence of the first two moments. 

\begin{figure}[!ht]
	\centering\resizebox*{0.65\textwidth}{0.25\textheight}{\includegraphics{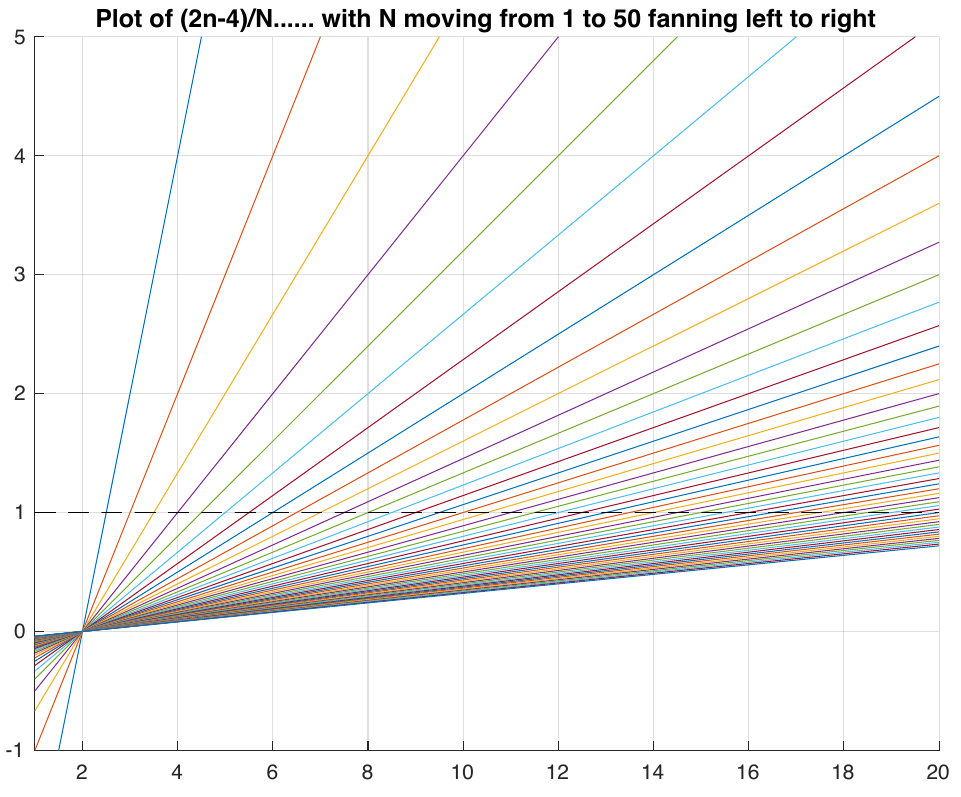}}
	\caption{Existence of moments for the {\tt VEnKF}. Each line corresponds to some number $\sfN$ with $\sfN$ moving from $1$ to $50$ fanning left to right. The `x-axis' corresponds to moment orders $n$ and a moment $n$ exists whenever the line $(2n-4)/\sfN$ is strictly less than one. }
	\label{fig:momentExistEnKF}
\end{figure}

Higher-order moments even in one dimension are still troublesome (for the {\tt VEnKF}, $\kappa=1$). In fact, the diffusion $\widehat{P}_t$ for the sample variance in case (\texttt{F1}) does not have any exponential moments in the stationary regime for any finite $\sfN\geq1$. That is, for any $t\geq 0$ and any finite $\alpha>0$ we have
\begin{equation}
 \mbox{\rm Law}(Q)=\widehat{\Gamma}_\infty~~~~\Longrightarrow~~~~\bbE\left[\exp{\left(\alpha\,\Vert \widehat{\phi}_{t}(Q)\Vert\right)}\right]=\infty
\end{equation}
 for any $\sfN\geq1$.

We also remark that the heavy tailed nature of this stationary distribution, in the case of the {\tt VEnKF}, implies that numerical stability in practice may be worrisome. In the stationary regime, it is realistic to expect samples from the tails in this case, and these may be large enough and/or frequent enough to cause numerical divergence. This property may lead to so-called catastrophic divergence as studied in, e.g., \cite{kelly2015}. In \cite{Harlim2010,Gottwald2013,kelly2015} mechanisms for catastrophic divergence are studied in complex nonlinear systems. Here we argue that even in linear systems, the heavy-tailed nature of the invariant measure of the sample covariance may lead to samples numerically large enough to cause numerical catastrophe in any practical computing system.

In case (\texttt{F2}) corresponding to the {\tt DEnKF}, we have that $L$ is reversible w.r.t. the probability measure $\widehat{\Gamma}_\infty$ on $]0,\infty[$ defined by,
\begin{equation}\label{invariant-measure-denkf}
		\widehat{\Gamma}_\infty(dx)~\propto~1_{]0,\infty[}(x)\,x^{\frac{\sfN}{2}-1}\,\exp{\left(-\frac{S\,\sfN}{4\,R}\,\left(x-2~\frac{A}{S}\right)^2\right)}\,dx
\end{equation}
Note this measure has Gaussian tails, and we contrast this with the heavy tailed nature of (\ref{invariant-measure-venkf}). This is significant, since it implies that the sample variance (and mean) of this {\tt DEnKF} will exhibit smaller fluctuations than the {\tt VEnKF}, and that all moments (including exponential moments) exist in this case for any choice of $\sfN\geq1$. This latter result is consistent with Theorem \ref{theo-existence-s-ric-proof-bis} at any time $t\geq0$ (and in the general multivariate setting). We can also expect better numerical stability (e.g. less outliers); including better time-discretisation properties \cite{hutzenthaler} in case (\texttt{F2}). These better fluctuation properties are already apparent in the preceding results (e.g. see Theorem \ref{theo-existence-s-ric-proof-bis}, \ref{theo-3-intro} and \ref{cor-H2-errorflows}) in the full multi-dimensional setting.

 As an illustrative example, take $A=20$ (i.e. the underlying signal model is highly unstable), $R=S=1$ and $\sfN=6$. In Figure \ref{fig:invariantmeasureEnKFvsDEnKF} we compare the invariant measure for the flow of the sample variance in each case.

\begin{figure}[!ht]
	\centering\resizebox*{0.75\textwidth}{0.25\textheight}{\includegraphics{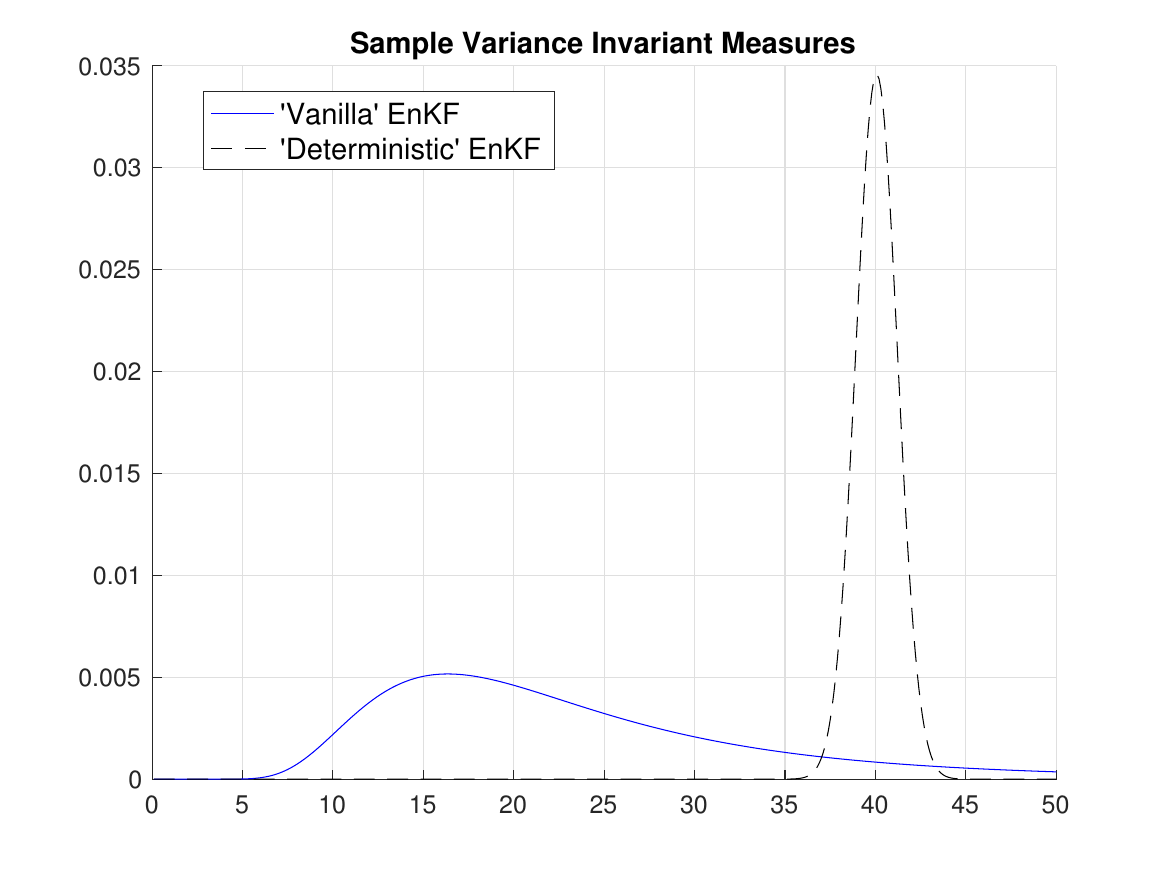}}
	\caption{The invariant measure of the sample variance of the `vanilla' {\tt EnKF} in case (\texttt{F1}), versus that of the `deterministic' {\tt EnKF} in case (\texttt{F2}).}
	\label{fig:invariantmeasureEnKFvsDEnKF}
\end{figure}

We see in Figure \ref{fig:invariantmeasureEnKFvsDEnKF} the heavy tails of the invariant measure (\ref{invariant-measure-venkf}) for the vanilla {\tt EnKF} sample variance, and conversely the Gaussian-type tails in the case (\ref{invariant-measure-denkf}) of the `deterministic' {\tt EnKF}. Note also the positioning of the mode/mean in each case. In case (\texttt{F1}) of the {\tt VEnKF}, $n$-th order moments exist only when $(2n-4)/\sfN$ is strictly less than one (in this case for $n<5$); while all moments exist in case (\texttt{F2}) for the {\tt DEnKF}.

The benefit and real interest in the scalar case is the ability to explicate the convergence rates, e.g. as in (\ref{ref-E-1-scalar}). We finally tun to the convergence/stability properties of the {\tt EnKF} sample variance and sample mean. In the case of the sample variance, we know from Theorem \ref{theo-stab-intro} that convergence of $\widehat{\phi}_t$ to its invariant measure $\widehat{\Gamma}_\infty$ (e.g. as depicted in Figure \ref{fig:invariantmeasureEnKFvsDEnKF} and described by (\ref{invariant-measure-venkf}) or (\ref{invariant-measure-denkf})) holds if $\sfN>4+(\kappa\,S)/(2\,R)$. Proof of this condition on $\sfN$ follows from Theorem \ref{theo-stab-intro}, the original multivariate statement of the same result in \cite[Theorem 2.4]{Bishop/DelMoral/multiDimRicc} and bounds on the mean of the sample variance flow and its inverse  \cite{2017arXiv171110065B,Bishop/DelMoral/multiDimRicc}. In \cite{2017arXiv171110065B} we also consider contraction and stability properties of the distribution of the sample covariance with respect to a particular Wasserstein metric; as opposed to the $\Lambda$-norm contraction used in Theorem \ref{theo-stab-intro}. An interesting result from \cite{2017arXiv171110065B} is that when $\kappa=0$, and for stable signal models (i.e. $A<0$), the Riccati diffusion (\ref{definition-Langevin-Riccati}) (describing the flow of the sample covariance) may converge faster to its invariant measure in (\ref{invariant-measure-denkf}), than the deterministic Riccati (\ref{nonlinear-KB-Riccati-modes}) does to its fixed point in (\ref{def-varpi}).

In one-dimensional ($d=1$) settings we may say more on the (stochastic) stability of the {\tt EnKF} sample covariance $\phi_t$ and sample mean $\psi_t$ based on the contraction properties of the stochastic transition matrix $\widehat{\calE}_t(Q)$ defined in (\ref{KB-semigroup-def-main-stoch}). It follows from Theorem \ref{det-E-theo} that we have the exponential decay estimate with $\sfN> 4 \vee(4n-2)\kappa$ which comes from \cite[Theorem 2.7]{Bishop/DelMoral/multiDimRicc},
\begin{equation}\label{ars-0-scalar}
\bbE\left[\widehat{\calE}_t(Q)^{n}\right]^{1/n}=\,\bbE\left[\exp{\left(n\int_0^t (A-\widehat{\phi}_s(Q)S)\,ds\right)}\right]^{1/n}
\displaystyle \leq\, c_{n}(Q)\,
\exp{\left(
-t\,\sqrt{\widehat{R}_{n}\widehat{S}_{n}}
\right)}
\end{equation} 
where $\widehat{R}_{n}$ and $\widehat{S}_{n}$ follow from (\ref{def-RS-repsilon}). In addition, there exists some function (of $\sfN$) $\lim_{\sfN\rightarrow \infty}\widehat{\nu}_{n}=0$ such that
\begin{equation}\label{ars-scalar}
\bbE\left[\widehat{\calE}_t(Q)^{n}\right]^{1/n}
\displaystyle\,=\,  c_{n}(Q)\,
\exp{\left(
-t\,(1-\widehat{\nu}_{n})\,\sqrt{A^2+RS} \right)}
\end{equation}
which we may relate (or contrast) with the exact contraction rate of the exponential semigroup associated with the deterministic Riccati equation in (\ref{ref-E-1-scalar}) describing the true filter variance in the classical Kalman-Bucy filter.  The rate parameter $\widehat{\nu}_{n}$ is different between the {\tt VEnKF} and {\tt DEnKF}. Details on the parameter $\widehat{\nu}_{n}$ are given in \cite{2017arXiv171110065B} but importantly for both $\kappa\in\{0,1\}$ we recover naturally the convergence rate of the deterministic Riccati flow in (\ref{ref-E-1-scalar}).

The exponential decay of the exponential semigroup $\widehat{\calE}_t(Q)$ plays a central role in the stability of the pair of processes $(\widehat{\phi}_t,\widehat{\psi}_t)$. For large time horizons the Lyapunov exponent can be estimated by the formula
\begin{equation}\label{Lyap-expo-def}
\frac{1}{t}\log{\widehat{\calE}_t(Q)}=
\frac{1}{t}~\int_0^t (A-\widehat{\phi}_s(Q)S)~ds~~\longrightarrow_{t\rightarrow\infty}~~ A-\widehat{\Gamma}_\infty(\iota)S,
\end{equation}
where $\widehat{\Gamma}_\infty$ denotes the reversible measure (\ref{invariant-measure-venkf}) or (\ref{invariant-measure-denkf}). We also have the following estimates of the Lyapunov exponent (\ref{Lyap-expo-def}) from \cite{2017arXiv171110065B}, and that relate also to the under-bias $\bbE[\widehat{\phi}_t(Q)]\leq \phi_{t}$. Let $\kappa=0$ and let $\mbox{\rm Law}(Q)=\widehat{\Gamma}_\infty$ be the reversible probability measure defined in (\ref{invariant-measure-denkf}). Then, for any $t\geq0$, we have
\begin{equation}
\sfN> 4~~\Longrightarrow~~ -\sqrt{A^2+RS}\leq~ A-\bbE[\widehat{\phi}_t(Q)]\,S~\leq~ -\sqrt{A^2+RS\left(1-\frac{4}{N}\right)}~<~0
\end{equation}
Similarly assuming $\kappa=1$ with $\mbox{\rm Law}(Q)=\widehat{\Gamma}_\infty$ and $\widehat{\Gamma}_\infty$ defined in (\ref{invariant-measure-venkf}) we have for any $t\geq0$,
\begin{equation}
\sfN> 4~~\Longrightarrow~~ -\sqrt{A^2+RS}\leq~ A-\bbE[\widehat{\phi}_t(Q)]\,S~\leq~ -\tfrac{\sqrt{A^2+RS\left(1-\left({4}/{N}\right)^2\right)}-{4A}/{N}}{1+4/N}~<~0
\end{equation}
As noted, the left hand inequalities in the preceding two equations follows immediately from the under-bias result $\bbE[\widehat{\phi}_t(Q)]\leq \phi_{t}$ 

From the contraction properties on $\bbE[\widehat{\calE}_t(Q)^{n}]$ we may deduce, in the scalar setting, strong stability results on the stochastic Riccati flow $\widehat{\phi}_t$ analogous to the deterministic setting, e.g. (\ref{ref-phi-1-stability}).  Similarly, strong stability results on the error flow $\widehat{\psi}_{t}$ follow from the contraction properties of $\bbE[\widehat{\calE}_t(Q)^{n}]$. Importantly, in the scalar $d=1$ case of $\widehat{\psi}_{t}(z,Q)$ we may relax the multivariate results like Theorem \ref{cor-H1-errorflows} and Theorem \ref{cor-H2-errorflows} which require more restrictive model (e.g. the strong observability/stability $\mu(A-{P}_{\infty}S)<0$ condition) and ensemble (particle) size assumptions.

From \cite[Theorem 5.10]{2017arXiv171110065B} we have that for any $\sfN>4\vee 4\kappa(n-1)$,
\begin{equation}\label{uniformly-Lipschitz-moments-riccati-scalar}
	{\bbE\left[ \Vert \widehat{\phi}_{t}(Q_1) - \widehat{\phi}_{t}(Q_2)\Vert^n \right]}^{1/n} \,\leq\, c_n\,\Vert Q_1-Q_2\Vert\, \exp{\left(
-t\,(1-\widehat{\nu}_{n})\,\sqrt{A^2+RS} \right)}
\end{equation}
for some function (of $\sfN$) $\lim_{\sfN\rightarrow \infty}\widehat{\nu}_{n}=0$. Note we have found no analogue of this result in the multivariate setting.

From \cite[Theorem 6.1]{2017arXiv171110065B} we have that for any $\sfN>2(4n+1)(1+4\kappa)$,
\begin{align}\label{uniformly-Lipschitz-moments-OU-scalar}
	&{\bbE\Big[ \Vert \widehat{\psi}_{t}(z_1,Q_1) -  \widehat{\psi}_{t}(z_2,Q_2)\Vert^n \Big]}^{1/n} \,\leq\, \nonumber\\
	&\qquad\qquad\qquad c_n(z_1,z_2,Q_1,Q_2)\,(\Vert z_1-z_2\Vert+\Vert Q_1-Q_2\Vert)\, \exp{\left(
-t\,(1-\widehat{\nu}_{n})\,\sqrt{A^2+RS} \right)}
\end{align}
for some function (of $\sfN$) $\lim_{\sfN\rightarrow \infty}\widehat{\nu}_{n}=0$. We may contrast this result with the more restrictive Theorem \ref{cor-H2-errorflows} in the multivariate setting. Note in the scalar setting we accommodate both the {\tt VEnKF} and {\tt DEnKF}, different initial variance states, and we recover, over fully infinite horizons, a continuous relationship with the optimal stability rates of (\ref{ref-E-1-scalar}).

The constants in (\ref{ars-scalar}), (\ref{uniformly-Lipschitz-moments-riccati-scalar}), and (\ref{uniformly-Lipschitz-moments-OU-scalar}) are given explicitly in terms of the model parameters in \cite{2017arXiv171110065B}. We remark that across these three stability results, the details of $\widehat{\nu}_{n}$ vary \cite{2017arXiv171110065B}, but importantly we recover the optimal (classical Kalman-Bucy) rates $\lim_{\sfN\rightarrow \infty}\widehat{\nu}_{n}=0$.

We consider an illustration of the fluctuation and stability properties of the sample variance in the different {\tt EnKF} variants. Consider again the model leading to Figure \ref{fig:invariantmeasureEnKFvsDEnKF}, and let $\widehat{{P}}_0=0$. The deterministic Riccati flow ($\sfN=\infty$, in (\ref{nonlinear-KB-Riccati-modes})) of the classical Kalman-Bucy filter and with the chosen model parameters ($A=20$, $R=S=1$) is given in Figure \ref{fig:detRiccatiflow}, along with $100$ sample paths of the sample variances for both the {\tt VEnKF} and the {\tt DEnKF} (with $\sfN=6$).

\begin{figure}[!ht]
	\centering
	\resizebox*{0.32\textwidth}{0.2\textheight}{\includegraphics{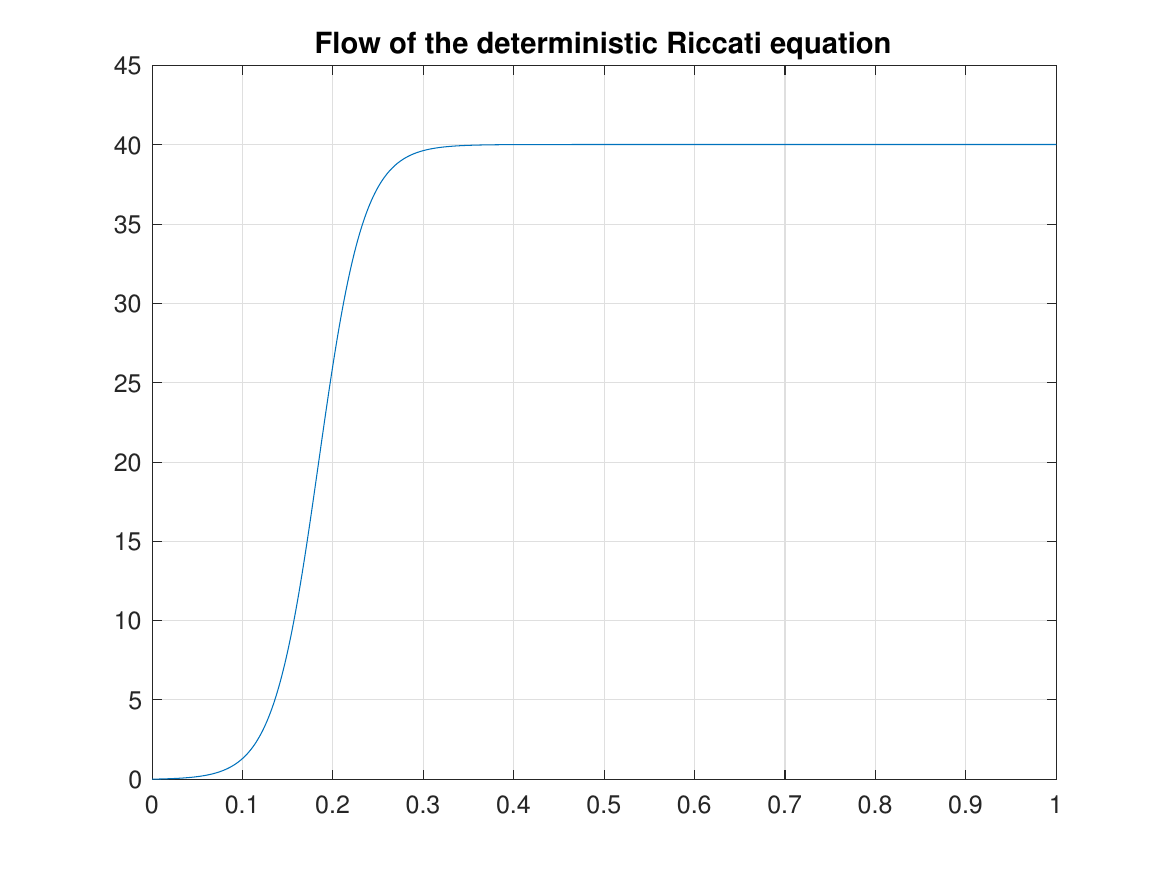}}
	\resizebox*{0.32\textwidth}{0.2\textheight}{\includegraphics{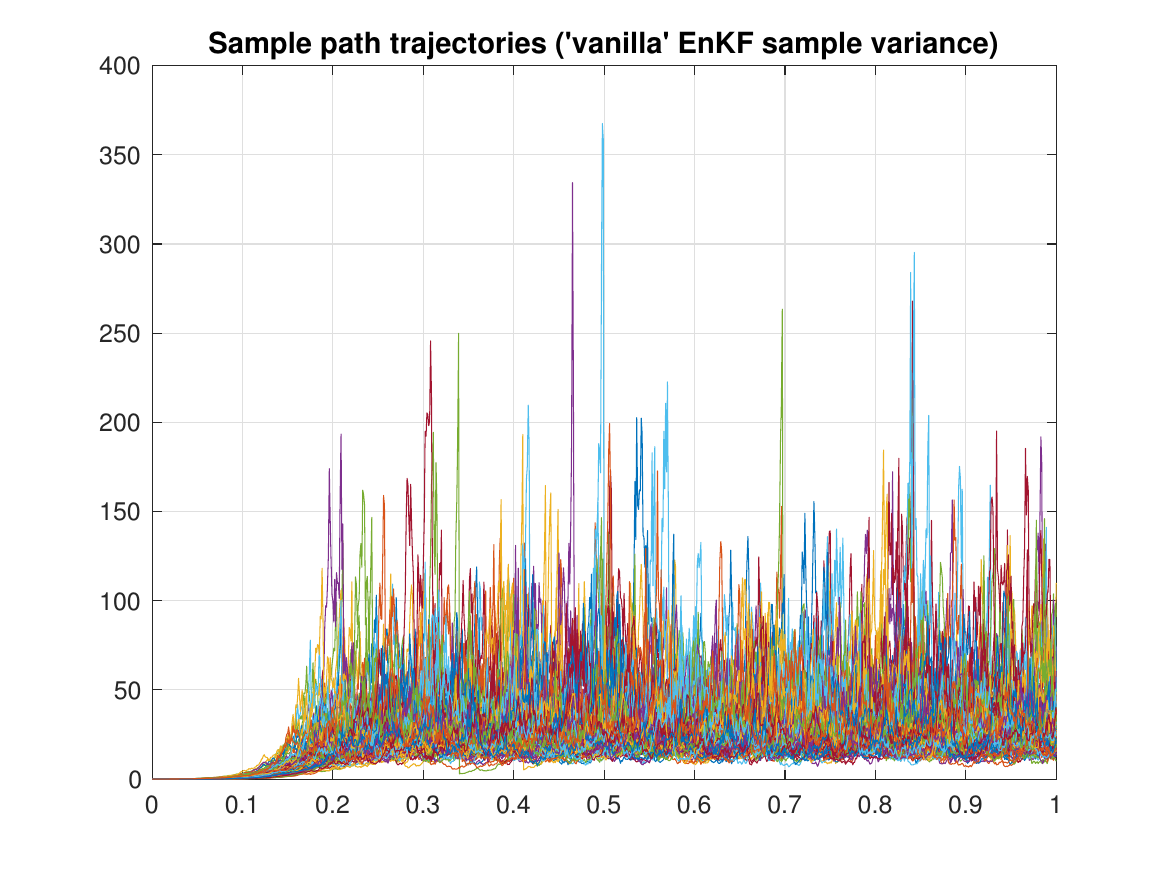}}
	\resizebox*{0.32\textwidth}{0.2\textheight}{\includegraphics{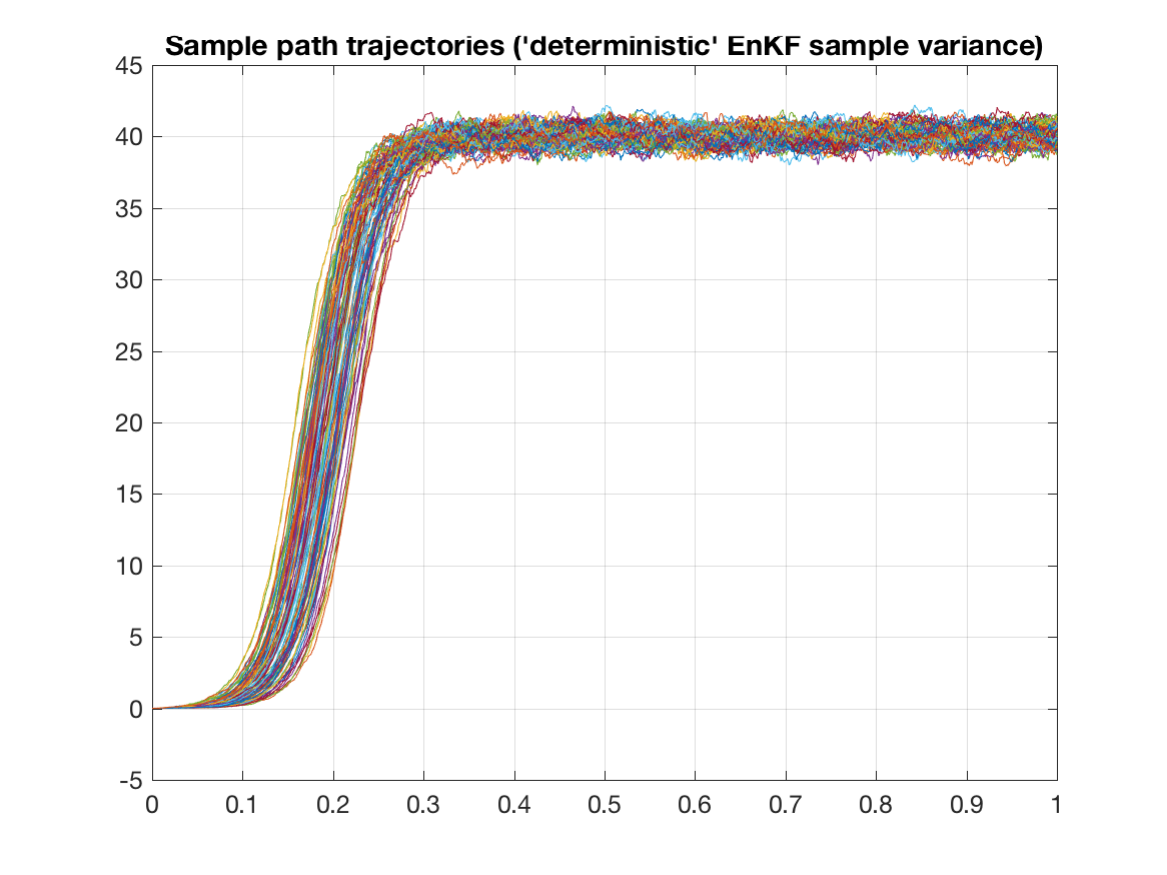}}
	\caption{Flow of the deterministic Riccati equation, and $100$ sample paths of the {\tt VEnKF} sample variance of case (\texttt{F1}), and $100$ sample paths of the {\tt DEnKF} sample variance of case (\texttt{F2}).}
	\label{fig:detRiccatiflow}
\end{figure}

Note in Figure \ref{fig:detRiccatiflow} the drastically reduced fluctuations in the `deterministic' {\tt EnKF} sample variance sample paths. At equilibrium, these fluctuations are related to the invariant measures of the two {\tt EnKF} varieties in (\ref{invariant-measure-venkf}) and (\ref{invariant-measure-denkf}).

In Figure \ref{fig:momentFlow} we plot the flow of the first two central moments and the $3$rd through the $9$th standardised central moments for both the {\tt VEnKF} and {\tt DEnKF} sample variance distribution. Recall that $\sfN=6$ in this example, and we expect moments of the {\tt VEnKF} sample variance in case (\texttt{F1}) to exist up to $n=4$ with $n=5$ the boundary case; while all moments exist for the {\tt DEnKF} of case (\texttt{F2}).

\begin{figure}[!ht]
	\centering
	\resizebox*{0.495\textwidth}{.35\textheight}{\includegraphics{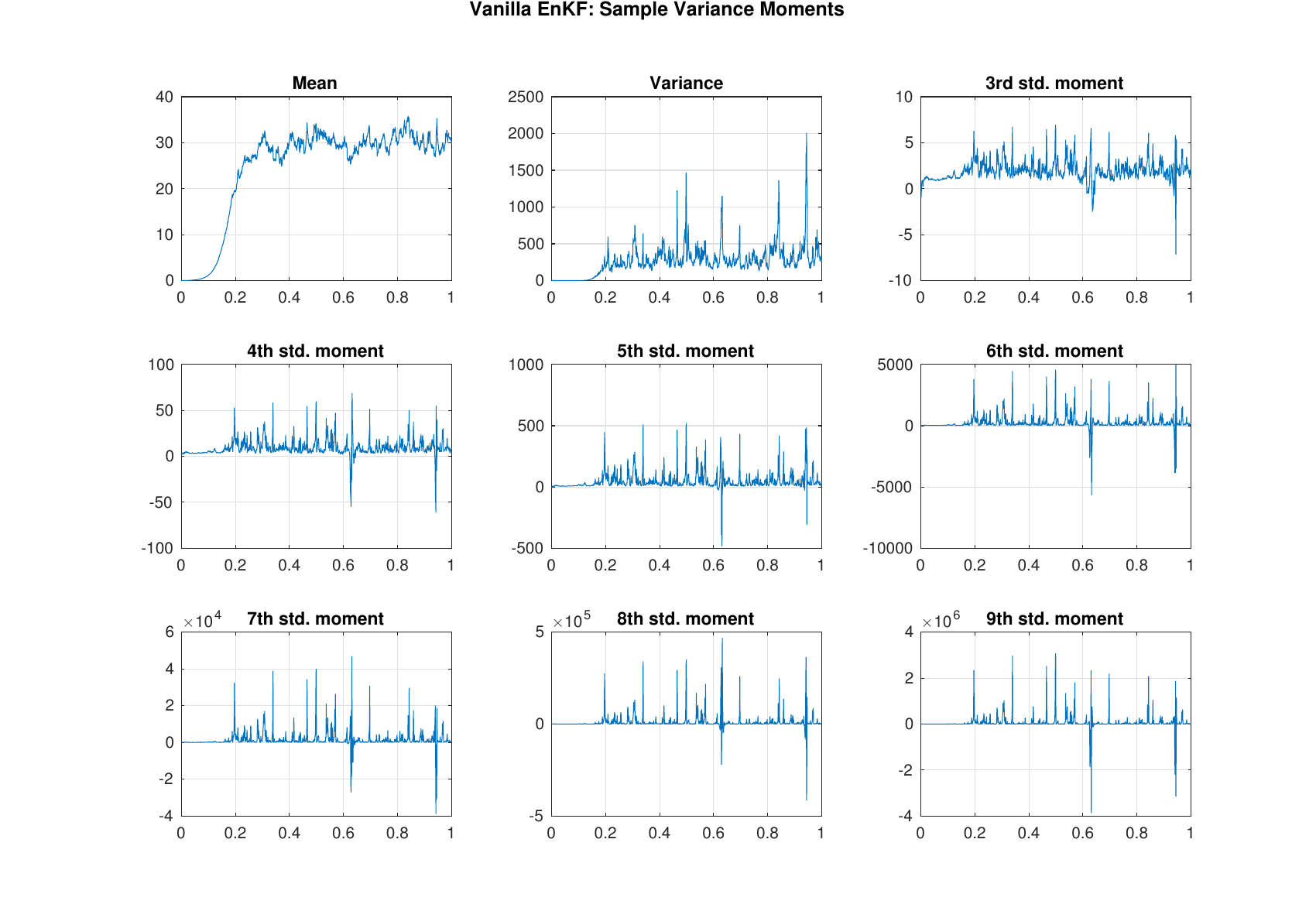}}
	\resizebox*{0.495\textwidth}{.35\textheight}{\includegraphics{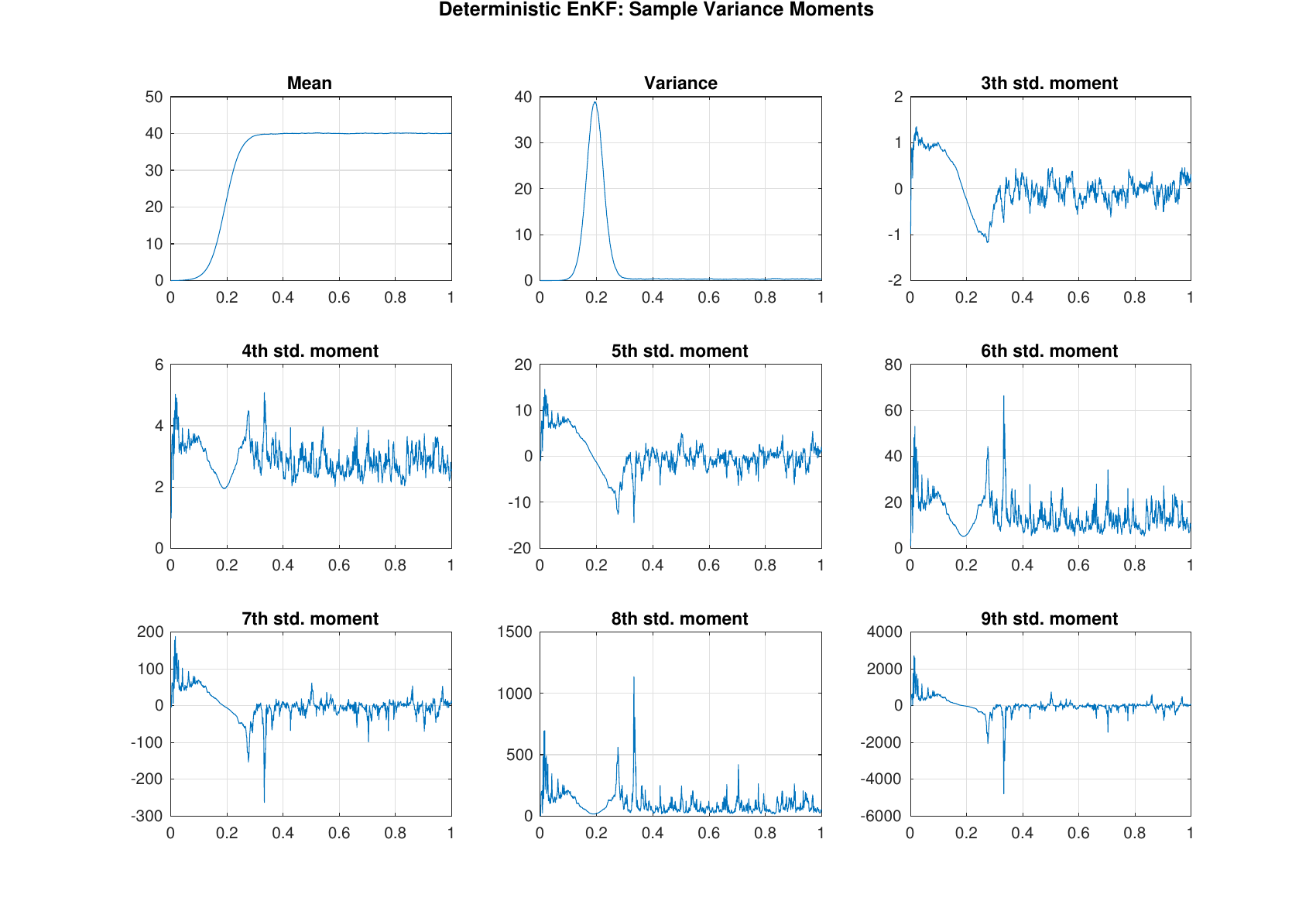}}\vspace{-.5cm}
	\caption{Flow of the sample variance moments for the {\tt VEnKF} and the {\tt DEnKF}.}
	\label{fig:momentFlow}
\end{figure}

We note in Figure \ref{fig:momentFlow} that the sample variance moments for the {\tt VEnKF} in case (\texttt{F1}) begin to destabilize around the $5th/6th$ moments as expected. Importantly, the mean of the sample variance for the {\tt VEnKF} is very negatively biased in this case, while the mean of the {\tt DEnKF} in case (\texttt{F2}) is quite accurate. Note also the very large variance in the sample variance for the {\tt VEnKF}.

Lastly, observe that (\ref{definition-Langevin-Riccati}) has non-globally Lipschitz coefficients in case (\texttt{F1}). The drift is quadratic, while the diffusion has a polynomial growth of order $3/2$. It follows by \cite{hutzenthaler} that the naive Euler-type time-discretization may blow up, regardless of the boundedness properties of the limiting (continuous-time) diffusion.

\section{Regularisations and Comparisons}

\subsection{Covariance Regularisation: Inflation}\label{sec-Inflation}

Let $(\calV^i_t,\calW^i_t,\calX_0^i)$ with ${1\leq i\leq \sfN+1}$ be $(\sfN+1)$ independent copies of $(\calV_t,\calW_t,\calX_0)$.  Now consider a modification of the individual particle update equations in the two cases of interest,
\begin{align}
(\texttt{F1})\qquad d\calX_t^{i,\varepsilon}~=&~A\,\calX_t^{i,\varepsilon}\,dt\,+\,R^{1/2}\,d\calV^i_t+\left(\widehat{P}^\varepsilon_{t}+ \varepsilon\,T\right)\,H^\prime\,R_1^{-1}\left[d\scrY_t-\left(H\calX_t^{i,\varepsilon}\,dt+R_1^{1/2}\,
d\calW_{t}^i\right)\right]\nonumber\\~\label{EnKF-sampled-diffusions-inflation}\\
(\texttt{F2})\qquad  d\calX_t^{i,\varepsilon}~=&~A\,\calX_t^{i,\varepsilon}\,dt\,+\,R^{1/2}\,d\calV^i_t+\left(\widehat{P}^{\varepsilon}_{t}+ \varepsilon\,T\right)\,H^\prime\,R_1^{-1}\left[d\scrY_t-H\left(\frac{\calX^{i,\varepsilon}_t+
 \widehat{X}^{\varepsilon}_{t}}{2}\right)dt\right] \nonumber
\end{align}
where $\varepsilon\in [0,\infty[$, and $T\in\bbS_r^0$ is some given reference matrix. Here, $\widehat{P}^\varepsilon_{t}=\widehat{\phi}^\varepsilon_t(Q)$ denotes the sample covariance-type function of $\calX_t^{i,\varepsilon}$ given by,
\begin{equation}\label{fv1-3-2-inflation}
\begin{array}{l}
\displaystyle\widehat{\eta}^{\,\varepsilon}_{t}:={\eta}^{\sfN,\varepsilon}_{t}=\frac{1}{\sfN+1}\sum_{i=1}^{\sfN+1}\delta_{\calX_t^{i,\varepsilon}} \\
	~\qquad\qquad\qquad\Longrightarrow\quad~\displaystyle
\widehat{X}^{\varepsilon}_{t} :=X^{\sfN,\varepsilon}_t=\frac{1}{\sfN+1}\sum_{i=1}^{\sfN+1}\calX_t^{i,\varepsilon} \quad\mathrm{and}\quad \displaystyle \widehat{P}^{\varepsilon}_{t} := P^{\sfN,\varepsilon}_t=\frac{\sfN+1}{\sfN}\,\calP_{\widehat{\eta}^{\,\varepsilon}_{t}} 
\end{array}\end{equation}

Recall the unified representation (for both the {\tt VEnKF} and the {\tt DEnKF}) for the flow of the sample mean, sample covariance, and the sample error flow in equations (\ref{EnKF-3}) through to (\ref{EnKF-2}).

Now consider the modification of the state estimator (sample mean) update equation resulting from the $(\varepsilon\,T)$-modification to the particle updates,
\begin{equation}
d\widehat{X}^\varepsilon_t\,=\,(A-(\widehat{P}^\varepsilon_t + \varepsilon\,T)S)\,\widehat{X}^\varepsilon_t\,dt+(\widehat{P}^\varepsilon_t + \varepsilon\,T)\,H^{\prime} R_1^{-1}\,d\mathscr{Y}_t+ \frac{1}{\sqrt{\sfN+1}}\,\Sigma^{1/2}_{\kappa,\varepsilon}(\widehat{P}^\varepsilon_t)\,d\calB_t \label{EnKF-3-inflation}
\end{equation}
with the mapping,
\begin{equation}
 \Sigma_{\kappa,\varepsilon}(Q)\,:=\,R+\kappa\,(Q+\varepsilon\,T)\,S\,(Q+\varepsilon\,T)\qquad\mbox{\rm with}\quad \kappa=\left\{
\begin{array}{rl}
1&\mbox{\rm in case (\texttt{F1})}\\
0&\mbox{\rm in case (\texttt{F2})}
 \end{array}\right. \label{sigma-mapping-inflation}
\end{equation}

With $\widehat{Z}^\varepsilon_t:=(\widehat{X}^\varepsilon_t-\mathscr{X}_t)$ we then also have that,
\begin{eqnarray}
d\widehat{Z}^\varepsilon_t&=&(A-(\widehat{P}^\varepsilon_t + \varepsilon\,T)S)\,\widehat{Z}^\varepsilon_t\,dt+(\widehat{P}^\varepsilon_t + \varepsilon\,T)\,H^{\prime} R_1^{-1}\,d\mathscr{W}_t \nonumber \\
&& \qquad\qquad\qquad\qquad\qquad\qquad\qquad\qquad\qquad -R^{1/2}\,d\mathscr{V}_t+ \frac{1}{\sqrt{\sfN+1}}\,\Sigma^{1/2}_{\kappa,\varepsilon}(\widehat{P}^\varepsilon_t)\,d\calB_t \nonumber \\
&\stackrel{ law}{=}& (A-(\widehat{P}^\varepsilon_t + \varepsilon\,T)S)\,\widehat{Z}^\varepsilon_t\,dt  +\Omega^{1/2}_{\kappa,\varepsilon} (\widehat{P}^\varepsilon_t)\,d\widehat{\scrB}_t
\label{EnKF-error-inflation}
\end{eqnarray}
with,
\begin{equation}
	\Omega_{\kappa,\varepsilon} \,:=\, \Sigma_{1,\varepsilon} + \frac{1}{{\sfN+1}}\, \Sigma_{\kappa,\varepsilon}
\end{equation}

In un-regularised ensemble Kalman filtering, we approximate $\phi_t$ by the sample covariance $\widehat{\phi}_t$ since the dimension of $\widehat{X}_t$ may be in the many millions; see \cite{evensen03}. However, when computing the sample covariance in high-dimensions, rank deficient estimation is common due to a lack of enough samples. Covariance inflation, leading to an approximation of the form $(\widehat{\phi}^\varepsilon_t(Q) + \varepsilon\,T)$, in the update equation (e.g. (\ref{EnKF-sampled-diffusions-inflation})), is a common, simple means of addressing this rank deficiency \cite{evensen03}. The under-bias result in Theorem \ref{underbiastheorem} or equation (\ref{ref-vp-max}) can also motivate the use of some form of regularisation such as inflation.

Note that the perturbation in the resulting flow of $\widehat{\phi}^\varepsilon_t(Q)$ comes from a (rather delicate) feedback loop adding $\varepsilon\,T$ to the covariance of the signal $\widehat{X}^\varepsilon_t$ at each instant. The flow of $\widehat{\phi}^\varepsilon_t(Q)$ is given by,
\begin{align}
d\widehat{P}^\varepsilon_t~=&~\left[(A-\frac{(1-\kappa)}{2}\varepsilon\, TS-\widehat{P}^\varepsilon_tS)\widehat{P}^\varepsilon_t+\widehat{P}^\varepsilon_t(A - \frac{(1-\kappa)}{2}\varepsilon \,TS-\widehat{P}^\varepsilon_tS)^{\prime}+R  +\kappa\,\varepsilon^2\,T\,S\,T \right]dt\nonumber\\
&\qquad\qquad+\frac{2}{\sqrt{\sfN}}\left[
{{}\widehat{P}^\varepsilon_t}^{1/2}\,d\calM_t~\Sigma^{1/2}_{\kappa,\varepsilon}(\widehat{P}^\varepsilon_t)\right]_{\mathrm{sym}} \label{EnKF-1-inflation}
\end{align}

In the limit $\sfN\rightarrow\infty$ we recover a perturbed, deterministic, Riccati equation that describes the flow of the limiting covariance. This perturbed Riccati equation is studied in \cite{Bishop/DelMoral/Pathiraja:2017,ap-2018-franklin}. For any size $\|\varepsilon T\|<\infty$, the perturbed Riccati flow qualitatively retains all the stability properties of the nominal Riccati flow (e.g. (\ref{ref-phi-1-stability}), but with a different steady state value), and the size of the error between the two grows in a well-quantified continuous way. 

In the limiting case, we have via \cite[Theorem 2.1]{ap-2018-franklin} that $\phi_{t}(Q)\leq \phi^\varepsilon_{t}(Q)$ in case (\texttt{F1}). In case (\texttt{F2}) we have that $ \phi^\varepsilon_{t}(Q)\leq \phi_{t}(Q)$ in the limit $\sfN\rightarrow\infty$.

For any $s\leq t$ and $Q\in\bbS_{d}^0$ we define the stochastic state-transition matrix,
\begin{equation}\label{KB-semigroup-def-main-inflation}
\widehat{\calE}^{\,\varepsilon}_{s,t}(Q):=\exp{\left[\oint_s^t\left(A- \varepsilon\,TS -\widehat{\phi}^{\,\varepsilon}_u(Q)S\right)du\right]} ~~\Leftrightarrow~~ \partial_t \widehat{\calE}^{\,\varepsilon}_{s,t}(Q) = \left(A- \varepsilon\,TS -\widehat{\phi}^{\,\varepsilon}_u(Q)S\right)\widehat{\calE}^{\,\varepsilon}_{s,t}(Q)
\end{equation}
Note that this semigroup $\widehat{\calE}^{\,\varepsilon}_{s,t}(Q)$ is associated with the evolution of the (inflation) regularised sample mean in (\ref{EnKF-3-inflation}) or the error flow (\ref{EnKF-error-inflation}) in both case (\texttt{F1}) and (\texttt{F2}). Unlike the un-regularised setting, this same semigroup is not directly related to the evolution of the sample covariance, in (\ref{EnKF-1-inflation}); for example, in case (\texttt{F1}) the semigroup associated with the evolution of the sample covariance is just $\widehat{\calE}_{s,t}(Q)$ as given in (\ref{KB-semigroup-def-main-stoch}) and studied throughout the preceding section.

We can comment on the effect of inflation regularisation on the contraction properties of $\widehat{\calE}^{\,\varepsilon}_{s,t}(Q)$, as compared e.g. to $\widehat{\calE}_{s,t}(Q)$. Firstly, it is worth noting, given the contraction estimates in Section \ref{sec-stability}, that,
\begin{equation}
	\mu((A-\varepsilon\,T S) - P S)~\leq~\mu(A - P S)
\end{equation}
for any fixed matrix $P\in\bbS^0_d$ and $S\in\bbS^0_d$. Arguing as in (\ref{change-basis-intro-formula}), when $S\in\bbS_d^+$, then up to a change of basis we can always assume that $S=I$. We then have,
\begin{equation}
	\mu(A)<\|\varepsilon\,T\|~~\Longrightarrow~~~~
\Vert \widehat{\calE}^{\,\varepsilon}_{s,t}(Q) \Vert \,\leq\, \exp{\left(\mu(A-\varepsilon\,T) (t-s)\right)}~\longrightarrow_{(t-s)\rightarrow\infty}~0
\end{equation}
which illustrates the added stabilising effects of $\varepsilon\,T$ in the extreme case in which $\widehat{P}^\varepsilon_t\,S$ has no stabilising effect at all. Contrast this with (\ref{change-basis-intro-formula}). Then one interpretation of the preceding relationship is that $\varepsilon\,T$ extends the set of signal matrices $A\in\bbM_d$ for which one may immediately achieve stabilisation (regardless of the effect of $\widehat{P}^\varepsilon_t\,S$). In practice, $\widehat{P}^\varepsilon_t\,S$ will also act to stabilise the filter, see e.g. (\ref{uniformly-Lipschitz-moments-OU-Lip}). Indeed, in the classical Kalman filtering setting (\ref{nonlinear-KB-mean}), (\ref{nonlinear-KB-Riccati}) with $\varepsilon=0$, the time-varying matrix $(A-{P}_tS)$ is stabilising \cite{Bishop/DelMoral:2016} for any $A\in\bbM_d$, even $A$ unstable. In the {\tt EnKF}, we know that $\widehat{{P}}_{t}$ will fluctuate about ${P}_t$, e.g. see Theorem \ref{theo-3-intro}. Therefore, the stabilisation properties of $(A-\widehat{{P}}_{t}S)$ are unclear; indeed the study of $\widehat{\calE}_{s,t}(Q)$ in the preceding Section \ref{sec-stability} is concerned with precisely this issue. The above implies that the addition of $\varepsilon\, T$ can act to counter the negative effects of this fluctuation (and directly add a stabilising effect on the state estimation error).

Finally, we have $\phi_{t}(Q)\leq \phi^\varepsilon_{t}(Q)$ in case (\texttt{F1}) and $\phi^\varepsilon_{t}(Q)\leq \phi_{t}(Q)$ in case (\texttt{F2}). The semigroup associated with the error flow (\ref{EnKF-error-inflation}) in both cases is the same. The inequality $\phi_{t}(Q)\leq \phi^\varepsilon_{t}(Q)$ in case (\texttt{F1}) suggests that the diffusion fluctuation in (\ref{EnKF-3-inflation}) or (\ref{EnKF-1-inflation}) will increase. However, we conversely expect that with $S\in\bbS^+_d$ we have $\mu((A-\varepsilon\,T S) - \phi^\varepsilon_{t}(Q) S)\leq \mu((A-\varepsilon\,T S) - \phi_{t}(Q) S)$ and thus we gain a type of stabilising effect. Inflation in case (\texttt{F1}) is then a delicate balancing tradeoff between adding noise to the diffusion coefficients (which may kill the existence of sample covariance moments, for example), and adding a stabilising effect on the sample mean error flow. When $\varepsilon>0$ is large enough we can achieve added stabilisation in case (\texttt{F2}), as compared to the non inflated case. This is not automatic as in case (\texttt{F1}) because $\phi^\varepsilon_{t}(Q)\leq \phi_{t}(Q)$. However, the fluctuations are (further) decreased with inflation in case (\texttt{F2}).

\subsection{Discretisation Matters}

The present article is primarily concerned with continuous-time filtering and {\tt EnKF} models. In practice, the stochastic models and analysis used for the continuous-time {\tt EnKF} are not applicable without an additional level of time-discretisation approximation. As alluded to earlier, the flow of the sample covariance for the {\tt VEnKF} has quadratic drift coefficients, while the diffusion term has a polynomial growth of order $3/2$. In this particular case, it follows by \cite{hutzenthaler} that a basic Euler time-discretization may blow up, regardless of the boundedness properties of the diffusion.

In contrast with continuous-time models, discrete-time signal and observation models lead to the so-called discrete-time {\tt EnKF}; e.g. see \cite{burgers1998analysis} and \cite{sakov2008deterministic} respectively for the corresponding {\tt VEnKF} and {\tt DEnKF} methods (also referenced in discrete-time earlier). Convergence of the discrete-time {\tt EnKF} models to their continuous-time counterparts (studied herein) with appropriate model time-step parameterisations is studied in \cite{Lange2021a,Lange2021b,Lange2021c}.

The purely discrete-time {\tt EnKF} is not defined by a single coupled diffusion process, but rather by a coupled two-step prediction-updating process (also known as forecast-analysis steps in the {\tt EnKF} and data assimilation literature). Moreover, the Gaussian-nature of the diffusion models (e.g. the Riccati diffusion) arising in the analysis of continuous-time {\tt EnKF} theory is also lost, and an inherent difficulty in discrete-time is the introduction of more sophisticated non-central chi-squared fluctuations. 

We emphasize that the discrete-time {\tt DEnKF} of Sakov and Oke \cite{sakov2008deterministic} is not consistent; i.e. it does not converge to the optimal filter as the number of particles tends to infinity, even in the linear-Gaussian case. The {\tt VEnKF} is consistent in discrete-time, see \cite{burgers1998analysis}. In discrete-time, another class of {\tt EnKF} methods, termed square-root {\tt EnKF} methods \cite{tippetsqrt2003,Lange2021b}, are consistent \cite{Kwiatkowski2015,Lange2020}. The discrete-time square-root and deterministic {\tt EnKF} methods are closely related (see \cite{sakov2008deterministic}) and in the continuous-time limit they converge to the same object \cite{Lange2021b,Lange2021c}, i.e. the continuous-time {\tt DEnKF} studied here. 

Analysis of the discrete-time {\tt VEnKF} \cite{burgers1998analysis} was studied in the linear-Gaussian setting in \cite{dm-horton-21}. That article presents a rather complete analysis of the fluctuations and the long-time behaviour of {\tt VEnKF} for one-dimensional models, including uniform estimates in the presence of transient and unstable latent signals.

\subsection{Particle Filter Comparisons}

We end this section with some theoretical comparisons between the discrete-time {\tt VEnKF} \cite{burgers1998analysis,dm-horton-21} and the particle filter (abbreviated {\tt PF} in this section) \cite{gordon1993,dm-96,DelMoral/Miclo:2000,doucetSMCpaper2000,kitagawa-93,kitagawa-96}.

We remark that the {\tt PF} and the {\tt EnKF} are, in general, built on different prediction/correction mechanisms.
The {\tt EnKF} uses an empirical gain function to weight the observations with the predicted state estimate in a manner akin to the update/correction stage of the classical Kalman filter. This mechanism which involves ``moving'' the corrected state estimate has the potential to stabilise the state estimate about the, possibly unstable, latent signal. This idea has played a central focus in this article. As shown in Theorem~\ref{underbiastheorem}, the empirical gain depends on an under-biased sampled covariance matrix which may fail to correct the effective unstable dimensions. The {\tt EnKF}  literature abounds with inconsistent but judicious ensemble transformations and regularisation methods like inflation/localisation-type procedures, aimed at addressing this issue. One basic inflation methodology is discussed in Section~\ref{sec-Inflation}, see also~\cite{Bishop/DelMoral/Pathiraja:2017}, in the context of continuous-time {\tt EnKF} methods where its action on the stabilisation properties of the {\tt EnKF} are shown. Conversely, both the continuous and the discrete-time {\tt PF} are based on genetic-type mutation-selection transitions: The basic discrete-time {\tt PF} methods, e.g. see \cite{gordon1993,dm-96,kitagawa-93}, evolve as a Markov chain on a product space. During the mutation transition, particles evolve independently according to the signal transition model. When an observation is delivered by the sensor, particles are selected with a probability proportional to their likelihoods. Importance sampling tricks can also be used to avoid degenerate mutations/predictions. 

Next we provide a detailed comparison of the {\tt PF} and {\tt VEnKF} for one-dimensional linear-Gaussian models with a view toward their tracking capability in the event of an unstable latent signal. The one-dimensional discrete-time version of (\ref{lin-Gaussian-diffusion-filtering}) has the following form,
 \begin{equation}\label{lin-Gaussian-diffusion-filtering-discrete}
\begin{split}
\scrX_{t+1}~&=~A\,\scrX_{t}\, \,+\, R^{1/2}\,\scrV_{t+1} \\
\scrY_t~&=~H\,\scrX_t\, \,+\, R_1^{1/2}\,\scrW_t
\end{split}
\end{equation}
where $t\in\bbN$ denotes the discrete time index and $(\scrW_t,\scrV_{t+1})$ is a sequence of $2$-dimensional 
Gaussian random variables with zero-mean and unit variance. The initial condition of the signal $\scrX_0$ is Gaussian with mean and variance denoted by $(X^-_0,P_0^-)$ (independent of $(\scrW_t,\scrV_{t+1})$), and $(A,H,R,R_1)$ are the model parameters. Any non-zero values for $(A,H,R,R_1)$ ensure that the model is (discrete-time) detectable and stabilisable. 

The discrete-time version of the conditional nonlinear McKean-Vlasov-type diffusion process (\texttt{F1}) discussed in (\ref{Kalman-Bucy-filter-nonlinear-ref}) is given by,
\begin{equation}\label{kalman-X-goth}
\left\{~\begin{aligned}
\calX_{t}\,&=\,\calX_t^-+\calG_{\overline{\eta}^-_t}~(\scrY_t-(H\,\calX^-_t+R^{1/2}_1\,\calW_t))\quad\mbox{\rm with}\quad \calG_{\overline{\eta}^-_t} \,:=\, H\calP_{\overline{\eta}^-_t}/(H^2\calP_{\overline{\eta}^-_t}
+R_1) \\
\calX_{t+1}^-\,&=\, A\,\calX_{t}+R^{1/2}\,\calV_{t+1}
\end{aligned}\right.
\end{equation}
In the above display, $\overline{\eta}_t^-$ denotes the conditional distribution of $\calX_t^-$ given $\calY_{t}^-:=(\scrY_0,\ldots,\scrY_{t-1})$ and  $(\calV_t,\calW_t,\calX_0)$ are independent copies of $(\scrV_t,\scrW_t,\scrX_0)$. Using a simple induction argument, it is straightforward to show that,
\begin{equation}
	\overline{\eta}^-_t \,=\, \mbox{\rm Law}(\calX_t^-~|~\calY^-_{t}) \,=\, \mbox{\rm Law}(\scrX_t~|~\calY^-_{t})
\end{equation}
is Gaussian with mean $X^-_t$ and variance $P^-_t$; and,
\begin{equation}
\overline{\eta}_t \,=\, \mbox{\rm Law}(\calX_t~|~\calY_{t}) \,=\, \mbox{\rm Law}(\scrX_t~|~\calY_{t})
\end{equation}
is also Gaussian with mean $X_t$ and variance $P_t$. Moreover, the conditional means $(X_t^-,X_t)$ and the variances
$(P_t^-,P_t)$ obey the Kalman filtering (update/correction and prediction) equations,
\begin{equation}\label{kalman-rec}
\left\{\begin{array}{rcl}
X_t&=&X_t^-+G_t~\left(\scrY_t-HX_t^-\right) \\
 P_t&=&(1-G_t H)P_t^-
\end{array}\right.\quad\mbox{\rm and} \quad\left\{\begin{array}{rclcrcl}
X_{t+1}^-&=&A X_{t}\\
 P_{t+1}^-&=& A^2P_{t}+R,
\end{array}\right.
\end{equation}
where in this section $G_t$ denotes the so-called Kalman gain parameter,
\begin{equation}
G_t:=HP_t^-/(H^2P^-_t+R_1)~~~\Longrightarrow~~~ 1-G_tH=1/(1+SP_t^-)\quad \mbox{\rm with}\quad S:=H^2/R_1
\end{equation}

The particle approximation of the nonlinear Markov chain discussed above is given by an interacting particle system defined sequentially for any $1\leq i\leq N+1$ by the formulae,
\begin{equation}\label{kalman-EnKF-def}
\left\{~\begin{aligned}
\calX^i_{t} \,&=\, \calX^{i-}_{t}+\widehat{G}_t~(\scrY_t-(H\calX^{i-}_t+R^{1/2}_1\calW^i_t))\quad\mbox{\rm with}\quad \widehat{G}_t:=H\widehat{P}^-_t/(H^2\widehat{P}^-_t
+R_1)\\
\calX^{i-}_{t+1} \,&=\, A\,\calX^{\,i}_{t}+R^{1/2}\,\calV^i_{t+1}
\end{aligned}\right.
\end{equation}
where $\widehat{P}^-_t$ denotes for the normalised sample variance
\begin{equation}
	\widehat{P}^-_t:=\frac{1}{N}\sum_{1\leq i\leq N+1}(\calX^{i-}_t-\widehat{X}^-_t)^2\quad \mbox{\rm with}\quad
\widehat{X}^-_t:=\frac{1}{N+1}\sum_{1\leq i\leq N+1}\calX^{i-}_t
\end{equation}
The above model coincides with the discrete-time version of the {\tt VEnKF} (i.e. in case {\tt F1}) in (\ref{EnKF-sampled-diffusions}), and follows from early results in \cite{burgers1998analysis}. The conditional mean $X_t=\bbE(\scrX_t\,|\,\calY_t)$ is approximated by the sample mean,
\begin{equation}
	\widehat{X}_t^{\scalebox{0.6}{\rm{EnKF}}} \,:=\, \frac{1}{N+1}\sum_{1\leq i\leq N+1}\calX^{i}_t
\end{equation}
The next theorem provides uniform mean-error estimates.
\begin{theorem}\label{timeuniform-discrete-time}
For any non-zero values of the model parameters $(A,H,R,R_1)$, any $n\geq 1$, and for $N\geq1$ sufficiently large, we have the uniform estimates
\begin{equation}
	\sup_{t\geq 0}\, \bbE\left[\vert \widehat{X}_t^{\scalebox{0.6}{\rm{EnKF}}}  -X_{t}\vert^n\right]^{1/n} ~\leq~c_n\frac{1}{\sqrt{N}} 
\end{equation}
We also have the conditional bias estimate,
\begin{equation}\label{condbiasmeaneqn-discrete}
	\bbE\left[ \Big|\, \bbE\left[\widehat{X}_t^{\scalebox{0.6}{\rm{EnKF}}}  \,|\, \calY_t\right] \,-\, X_{t}\, \Big|^n\right]^{1/n} \,\leq\, c_n(P_0^-)\,\frac{1}{\sfN}
\end{equation}
\end{theorem}

The proof of this discrete-time, one-dimensional, version of (\ref{ref-psi-1-est-scalar}) is given in \cite{dm-horton-21}; see also \cite{2017arXiv171110065B} for continuous-time analogues. Contrast this result also with the continuous-time multi-variate result in Theorem \ref{biastheoremmean}.

Particle filters, eg. see \cite{gordon1993,dm-96,dm-g-1999, dm-98,DelMoral/Miclo:2000,doucetSMCpaper2000,doucetSMCbook2001,kitagawa-93,kitagawa-96}, are a popular method for approximate filtering in nonlinear state space models in relatively low dimensions. The genetic-type particle filter (also referred to as the bootstrap filter) is a Markov chain with a mutation and a selection transition,
\begin{equation}
	\left(\xi^{i-}_t\right)_{1\leq i\leq N}\in\bbR^N ~~\stackrel{ selection}{-\!\!\!-\!\!\!-\!\!\!-\!\!\!-\!\!\!\longrightarrow} ~~
\left(\xi^i_t\right)_{1\leq i\leq N}\in\bbR^N ~~
\stackrel{ mutation}{-\!\!\!-\!\!\!-\!\!\!-\!\!\!-\!\!\!\longrightarrow} ~~\left(\xi^{i-}_{t+1}\right)_{1\leq i\leq N}
\end{equation}
The initial configuration $\left(\xi^{i-}_0\right)_{1\leq i\leq N}$ is defined by sampling $N$ independent copies of $\scrX_0$.
In its simplest form, the selection transition involves sampling $N$ independent random variables $\left(\xi^i_t\right)_{1\leq i\leq N}$ with the weighted distribution
\begin{equation}
	\sum_{1\leq i\leq N}~\frac{e^{-(\scrY_t-H\xi^{i-}_t)^2/(2R_1)}}{\sum_{1\leq j\leq N}e^{-(\scrY_t-H\xi^{j-}_t)^2/(2R_1)}}~\delta_{\xi^{i-}_t}
\end{equation}
The corresponding mutation transition coincides with prediction stage of the {\tt VEnKF} above; that is for any $1\leq N$ we set
\begin{equation}\label{def-mutation}
\xi^{i-}_{t+1}= A\,\xi^{\,i}_{t}+R^{1/2}\,\calV^i_{t+1}.
\end{equation}
In this context, the conditional means $X_t=\bbE(\scrX_t~|~\calY_t)$ are approximated by the sample means
\begin{equation}
	\widehat{X}_t^{\scalebox{0.6}{\rm{PF}}} \,:=\, \frac{1}{N}\sum_{1\leq i\leq N}\xi^i_t
\end{equation}

A mature literature on the time-uniform contraction/convergence and fluctuation results for the particle filter exists, and a survey of this topic is beyond the scope of this article. However, typically the time-uniform particle filtering estimates discussed in the literature rely on mixing-type or certain contractive conditions on the mutation transitions, e.g. \cite{dm-g-1999,mf-dm-04,Oudjane2005,VanHandel2009bb,mf-dm-13,Douc2014}. In the case of linear-Gaussian state transition models, none of these conditions hold for general unstable transient signals.

A natural question: Can the particle filter track unstable (latent) signals, like the {\tt VEnKF} can (as per Theorem \ref{timeuniform-discrete-time}), in the simple linear-Gaussian setting of this section? Unfortunately, as we now show (and contrast with the {\tt VEnKF} and Theorem \ref{timeuniform-discrete-time}) the answer here is rather negative.
 
Following ideas of Mathieu Gerber (personal communication), whenever $A>1$ we have
\begin{equation}
\xi^{i-}_{t} \,\geq\,  A \inf_{1\leq i\leq N}\xi^{\,i}_{t-1}+R^{1/2}\,\calV^{\star-}_{t}
\qquad
\mbox{\rm and}\qquad
\xi^{i}_t \,\geq\, \inf_{1\leq i\leq N}\xi^{i-}_t\quad
\end{equation}
with
\begin{equation}
	\calV^{\star-}_{t} \,:=\, \inf_{1\leq i\leq N}\calV^i_{t}
\end{equation}
This implies that,
\begin{equation}
\xi^{i-}_t\wedge \xi^{i}_t ~\geq~ X^{\star-}_t\,:=\,AX^{\star-}_{t-1}+R^{1/2}\calV^{\star-}_t
\qquad
\mbox{\rm with}\qquad X^{\star-}_0=\inf_{1\leq i\leq N}\xi^{i-}_{0}
\end{equation}
Thus, for any given initial conditions $\xi^{-i}_0=x^i_0$ with
\begin{equation}
\inf_{1\leq i\leq N}x^i_0 ~\geq~ \epsilon+\frac{R^{1/2}}{A-1}\sqrt{2\log{(N)}}
\end{equation}
for some $\epsilon>0$, we have
\begin{equation}
\inf_{1\leq i\leq N}\bbE(\xi^{i}_t~|~\calY_t)~\geq~ A^t\,\epsilon+A^t\,\frac{R^{1/2}}{A-1}\sqrt{2\log{(N)}}+R^{1/2}
\sum_{0\leq s<t}\,A^{s}~\bbE\left[\calV^{\star-}_{t-s}\right]
\end{equation}
Recalling that $\bbE[\max_{1\leq i\leq N} U_i]\leq\sqrt{2\log(N)}$ for any sequence of $N$ independent centered Gaussian random variables $U_i$ with unit variance, we conclude that
\begin{align}
\inf_{1\leq i\leq N}\bbE[\xi^{i}_t~|~\calY_t] \,&\geq\, A^t\,\epsilon+A^t\,R^{1/2}\sqrt{2\log{(N)}}\,\left(\frac{1}{A-1}-
\sum_{1\leq s\leq t}~A^{-s}~\right) \nonumber\\
\,&=\, A^t\,\epsilon+\frac{R^{1/2}}{A-1}\sqrt{2\log{(N)}}
\end{align}
This yields the almost sure divergence result
\begin{equation}\label{pre-div}
\lim_{t\rightarrow\infty}\bbE\left[\widehat{X}_t^{\scalebox{0.6}{\rm{PF}}}~|~\calY_t\right] \,=\,+\infty
\end{equation}
This result is not restricted to proportional selection, but rather holds for any unbiased selection transition. Importantly, we emphasise that this result is true even when the unstable latent signal moves to $-\infty$; in which case there is very quickly a drastic divergence between the particle filtering estimate and the latent signal (and the optimal filter). This already indicates that the particle filter is not able to track unstable signals.  

\begin{theorem}\label{pfblowup-discrete-time}
For any non-zero values of the model parameters $(H,R,R_1)$, and any $A>1$ and $P_0^->0$, any $n\geq 1$, and for any $N\geq1$ we find,
\begin{equation}
\lim_{t\rightarrow\infty}\bbE\left[\left|\widehat{X}_t^{\scalebox{0.6}{\rm{PF}}}-X_{t}\right|^n\right]^{1/n} \,=\,\infty \,=\,\lim_{t\rightarrow\infty}\bbE\left[\left\vert \widehat{X}_t^{\scalebox{0.6}{\rm{PF}}}-\scrX_{t}\right\vert^n\right]^{1/n}
\end{equation}
\end{theorem}

Before we proceed to the proof, we contrast the preceding result with the time-uniform bound and error control achievable with the {\tt EnKF} detailed in Theorem \ref{timeuniform-discrete-time} above. Theorem \ref{pfblowup-discrete-time} states that there is no hope in stabilising the particle filtering estimate around an unstable and transient (latent) signal when using the same mutation-prediction (\ref{def-mutation}) as the {\tt EnKF} (and which is common in basic particle filtering implementations, cf. \cite{gordon1993,dm-96,DelMoral/Miclo:2000,doucetSMCpaper2000,kitagawa-93,kitagawa-96}). Increasing the number of particles yields no (long-term) benefits here. More specifically, and in contrast with the {\tt EnKF}, the selection-correction stage of the above particle filter cannot compensate or correct for an unstable prediction-mutation.

\begin{proof}
Following the proof of (\ref{pre-div}), for any $A>0$, note that
\begin{equation}
\xi^{i-}_{t} \,\geq\,  A \inf_{1\leq i\leq N}\xi^{\,i}_{t-1}\,+\,R^{1/2}\,\calV^{\star}_{t}
\qquad
\mbox{\rm and}\qquad
\xi^{i}_t\geq \inf_{1\leq i\leq N}\xi^{i-}_t\qquad
\end{equation}
with
\begin{equation}
	\calV^{\star}_{t} \,:=\, (-\scrV_t) \,\wedge\, \calV^{\star-}_{t}
\end{equation}
Define,
\begin{equation}
	 X^{\star}_t \,:=\, A\,X^{\star}_{t-1}+R^{1/2}\calV^{\star}_t
\qquad
\mbox{\rm with}\qquad X^{\star}_0 \,:=\, \inf_{1\leq i\leq N}\xi^{i-}_{0}
\end{equation}
Then it follows that,
\begin{equation}
	\xi^{i-}_t \,\wedge\, \xi^{i}_t ~\geq~ X^{\star}_t
\end{equation}
We also then have,
\begin{align}
	X^{\star}_t-\scrX_{t} \,&=\, A(X^{\star}_{t-1}-\scrX_{t-1})+R^{1/2}\left(\calV^{\star}_t \,+\, (-\scrV_t)\right) \nonumber\\
					\,&\geq\,  A(X^{\star}_{t-1}-\scrX_{t-1}) \,+\, 2R^{1/2}\calV^{\star}_{t}
\end{align}
which implies that,
\begin{equation}
	\frac{X^{\star}_t-\scrX_{t}}{A^t} ~\geq~ (X^{\star}_{0}-\scrX_{0})\,+\, 2R^{1/2}\sum_{1\leq k\leq t}\frac{\calV^{\star}_{k}}{A^{k}}
\end{equation}
Thus, for any $A>1$ and $\epsilon>0$ on the event
\begin{equation}
\Omega_{\epsilon} \,:=\, \left\{\scrX_{0} \,\leq\, 0\quad \mbox{\rm and}\quad
	X^{\star}_0 \,\geq\, \epsilon+\frac{2R^{1/2}}{A-1}~\sqrt{2\log{(N+1)}}\right\}
\end{equation}
 we readily check that
\begin{equation}
	\bbE\left[ \widehat{X}_t^{\scalebox{0.6}{\rm{PF}}}-\scrX_{t}~|~X^{\star}_{0},\scrX_{0}\right] \,\geq\, 
\bbE\left[X^{\star}_t-\scrX_{t}~|~X^{\star}_{0},\scrX_{0}\right] ~\geq~  \epsilon\,A^t~\longrightarrow_{t\rightarrow\infty}~\infty
\end{equation}
Moreover, we have,
\begin{align}
\bbE\left[\left\vert\widehat{X}_t^{\scalebox{0.6}{\rm{PF}}}-\scrX_{t}\right\vert~|~X^{\star}_{0},\scrX_{0}\right]1_{\Omega_{\epsilon}}\,&\geq\,
\left\vert\,\bbE\left[\widehat{X}_t^{\scalebox{0.6}{\rm{PF}}}-\scrX_{t}~|~X^{\star}_{0},\scrX_{0}\right]\,\right\vert~1_{\Omega_{\epsilon}} \nonumber\\
\,&=\, \bbE\left[\widehat{X}_t^{\scalebox{0.6}{\rm{PF}}}-\scrX_{t}~|~X^{\star}_{0},\scrX_{0}\right]\,1_{\Omega_{\epsilon}}\\
\,&\geq\, \epsilon\,1_{\Omega_{\epsilon}}  \, A^t \nonumber
\end{align}
In discrete-time, the variance $P_t=\bbE[(X_t-\scrX_t)^2]$ is also uniformly bounded with respect to any time horizon, e.g. see \cite{anderson79,Lancaster1995} (and similarly to (\ref{upper-bound-Qhi}) in continuous-time), and thus we have,
\begin{equation}
\epsilon\,\bbP\left[\Omega_{\epsilon}\right]\,A^t \,\leq\, \bbE\left[\left\vert \widehat{X}_t^{\scalebox{0.6}{\rm{PF}}}-\scrX_{t}\right\vert\right] \,\leq\, \bbE\left[\left\vert\widehat{X}_t^{\scalebox{0.6}{\rm{PF}}}-X_{t}\right\vert\right]+c
\end{equation}
for some finite constant $c<\infty$.

Finally we confirm the non-zero probability,
\begin{equation}
P_0^->0 \quad \Longrightarrow\quad 
\bbP\left[\Omega_{\epsilon}\right]\,=\,
\bbP\left[\scrX_{0}\leq 0\right]\,\bbP\left[\scrX_{0} \, \geq\, \epsilon+\frac{2R^{1/2}}{A-1}\sqrt{2\log{(N+1)}}\right]^N~>~0
\end{equation}
This ends the proof of the theorem.
\end{proof}

To conclude this discussion, we note briefly that it is possible to stabilise the prediction-mutation step (about an unstable latent signal) at the expense of also changing the selection-correction stage (i.e. via importance sampling). For example, a mutation-prediction step sampled according to the so-called optimal proposal $\sim\,\mathrm{Law}(x_t \,|\,x_{t-1},\scrY_t)$, see \cite{doucetSMCpaper2000}, or earlier in \cite[Example 3]{dm-98},~\cite[Sections 2.4.3 and 12.6.6]{mf-dm-04}, is a stable option whenever $A/(1+H^2R/R_1)<1$. However, in general nonlinear filtering problems these mutation transitions and the corresponding importance selection weights are intractable \cite{doucetSMCpaper2000}. The terminology ``optimal proposal'', see \cite{doucetSMCpaper2000}, is somehow confusing as this importance sampling strategy and the one discussed in (\ref{def-mutation}) have the same Feynman-Kac-type mathematical structure and sampling according to this proposal doesn't minimize the asymptotic variance. In this context, following~\cite[Section 4.2.2]{dm-98}, we can use an auxiliary local particle approximation to sample $\mathrm{Law}(x_t \,|\,x_{t-1},\scrY_t)$ and compute the corresponding importance weights. Given the topic of this article and this section in particular, we note an interesting approach in \cite{Papadakis2010} employing an ensemble Kalman filter to define a proposal distribution (i.e. in the mutation step) that depends on the observation history \cite{dm-98,Papadakis2010}. We do not explore this topic in further generality here.

\section{Some Topics for Discussion}

\subsection{Comments on the Results Presented}

In places, we switch between rather quantitative estimates to those more qualitative in nature. In part this is to simplify presentation, or when the details are (likely) not tight and thus perhaps of little quantitative interest. In some in places it is because we did not obtain more precise descriptions of the estimates involved. Refining these estimates may be of practical interest in some cases; e.g. when deriving estimates on the required number of particles $\sfN$ for stability of the sample covariance (or convergence to its invariant measure). 

The results presented thus far consisted of constants, e.g. $c$, $c_n$, $c_\tau$, etc, that depend on the model parameters $(A,R,S)$, but importantly not on the ensemble size $(\sfN+1)$ or the time horizon $t\in[0,\infty[$. Due to the dependence on the model (e.g. $(A,R,S)$), these constants depend implicitly (via the matrix norms used) on the underlying signal dimension $d$. It would be of interest to pull this dependence out more explicitly depending on the matrix norm we are using, so as to quantify, at least in some general sense, the tradeoff between $\sfN$ and $d$. For example, in Theorem \ref{theo-3-intro} or Theorem \ref{theo-fluc-sample-mean} detailing the fluctuation of the sample covariance and sample mean about their limiting covariance and (Kalman-Bucy) state estimate values, it would be of interest to know how this fluctuation scales with dimension $d$, say e.g. with fixed $\sfN$. Unfortunately, the proof tools used in the development of this work does not lend itself naturally to this analysis.

The matrix $S:=H^{\prime}R_1^{-1}H$ plays a critical role throughout with regards to obtaining time-uniform fluctuation and then subsequently stability/convergence results. In particular, the assumption that $S\in\bbS_d^+$ is strictly positive-definite, i.e. Assumption \ref{mainAssumpObs}, is needed in numerous places. This assumption amounts to a type of strong observability result; e.g. a requirement on the ``fullness'' of the observations and the size and rank of the observation matrix $H$. It is worth emphasising that this assumption appears in many technical articles discussing the performance properties of the ensemble Kalman filter; e.g. \cite{Kelly2014,DelMoral/Tugaut:2016,tong2016stability,DelMoral/Kurtzmann/Tugaut:2017,deWiljes2018,deWiljes2019}.  Typically, the tools used in the proofs in \cite{DelMoral/Tugaut:2016,Bishop/DelMoral/Niclas:2017,2017arXiv171110065B,Bishop/DelMoral/multiDimRicc} are not sophisticated enough to accommodate zero eigenvalues of $S$. A basic example of this deficiency is in the proof of time-uniform moment boundedness of $\widehat{{P}}_{t}$, stated in Theorem \ref{theo-existence-s-ric-proof-bis}. In that proof, we resort to taking trace or eigenvalue-type reductions of the matrix-valued Riccati diffusion and studying a scalar comparison Riccati equation. This scalar reduction means that we must look at the minimum eigenvalue of $S$ (because it appears with a minus sign in the Riccati equation) and thus we cannot allow this value to be zero (because we would lose this term completely in the scalar comparison). To obtain uniform-in-time bounds, one needs the stabilising effect of this non-zero $S$ in the scalar comparison. See the proof in \cite[Theorem 2.2]{Bishop/DelMoral/multiDimRicc} for this very transparent example. In this example, one may relax the condition on $S$ to $S\in\bbS_d^0$ at the expense of time exponentially growing bounds. Related difficulties in allowing $S\in\bbS_d^0$ instead of $S\in\bbS_d^+$ arise in numerous other places (and as noted in other related works \cite{Kelly2014,DelMoral/Tugaut:2016,tong2016stability,DelMoral/Kurtzmann/Tugaut:2017,deWiljes2018,deWiljes2019}). One difficulty is related to stability of the (time-varying) matrix $(A-\widehat{{P}}_{t}\,S)$ and the positive-definiteness properties of product $\widehat{{P}}_{t}\,S$ as discussed subsequently.

We have focused significant effort on relaxing the assumption that the underlying signal is stable. Note that if $A$ is stable, i.e. $\mathrm{Absc}(A)<0$, then the stability of $\mu(A-\widehat{{P}}_{t}\,S)$, for some log-norm, may be trivially inherited whenever $S\in\bbS_d^+$ via a change of coordinates, see \cite{DelMoral/Tugaut:2016}. We see here again the use of $S\in\bbS_d^+$ as it pertains to the product $\widehat{{P}}_{t}\,S$. If $S\in\bbS_d^0$ is only positive semi-definite, then one can construct counterexamples such that even if $\widehat{\phi}_t = (\phi_t + \widehat{\varphi}_t/ \sqrt{N})\in\bbS_d^+$ is positive definite, there exists flows $\widehat{\varphi}_t$ such that $\mu(A-\widehat{{P}}_{t}\,S)=\mu(A-\phi_t\,S-\widehat{\varphi}_t\,S/ \sqrt{N})>0$. The fluctuation term $\widehat{\varphi}_t/ \sqrt{N}$ \textit{might} not interact well with the only positive semi-definite $S\in\bbS_d^0$. The assumption $\mu(A)<0$ is made in \cite{DelMoral/Tugaut:2016} in the linear-Gaussian setting and follow also in, e.g., \cite{Kelly2014,tong2016stability,DelMoral/Kurtzmann/Tugaut:2017,deWiljes2019} when reducing those studies to the linear-Gaussian setting. 

If $A$ is allowed to be unstable, then the asymptotic (time-varying) stability of $(A-{{P}}_{t}\,S)$ in the classical Kalman filter follows under so-called detectability (or observability) conditions \cite{Kwakernaak72,Park97,VanHandel2009}. Detectability intricately relates the relevant rank deficient directions in $P_t$ and $S$ in terms of the unstable directions in $A$ (i.e. it basically ensures those directions of $A$ that are unstable are observed (as captured by $S$) and non-zero weighted in the update Kalman gain via $P_t$). The rank of the sample covariance $\widehat{{P}}_{t}$ is at most $\sfN\geq1$. If $\sfN<d$, then $\widehat{{P}}_{t}$ is almost surely rank deficient and thus has zero eigenvalues in some directions. In general, we cannot control the directions in which the random, sub-rank, $\widehat{{P}}_{t}$ has zero eigenvalues (e.g. to play nicely with $S$ in the sense of detectability). If $A$ is unstable in those directions, the filter is consequently unstable in those directions. Thus, there is a basic, unavoidable, but also transparent tradeoff in requiring either stability of $\mathrm{Absc}(A)<0$ or sufficiently large ensemble sizes $\sfN\geq d$ in the derivation of uniform-in-time stability results for the {\tt EnKF}. 

In the stability results stated in this work, we emphasised unstable models $A$ but required sufficiently large ensemble sizes $\sfN\geq d$. Nevertheless, most stability results stated in this work with the hypothesis that ``$\sfN$ is sufficiently large'' may be restated with this condition replaced with ``$\sfN\geq 1$ and $\mathrm{Absc}(A)<0$'', and such results hold time-uniformly over infinite time horizons. In \cite{Bishop/DelMoral/Niclas:2017,2017arXiv171110065B,Bishop/DelMoral/multiDimRicc} the details on ``$\sfN$ is sufficiently large'' are given more explicitly. Note some results that do not consider or rely on the long time stability behaviour of the samples, e.g. the fluctuation size of the sample covariance about its true value, hold with $\sfN\geq 1$ and any matrix $A$, e.g. this is true for the {\tt DEnKF} in Theorem \ref{theo-3-intro}.

We remark that the assumption that the true Kalman-Bucy filter is stable in the sense $\mu(A-{{P}}_{\infty}\,S)<0$ is used in a number of the fluctuation (on the sample mean) and long-time behavioural results given in the cited, prior, work \cite{DelMoral/Tugaut:2016,Bishop/DelMoral/Niclas:2017,2017arXiv171110065B,Bishop/DelMoral/multiDimRicc,Bishop/DelMoral/STV2018}. It was originally believed by the authors that this condition was a stronger assumption than the more natural condition $\mathrm{Absc}(A-{{P}}_{\infty}\,S)<0$. The latter follows from the very natural model assumptions of detectability and stabilisability, see (\ref{steady-state-eq}) and the discussion following that equation, e.g. \cite[Theorems 9.12, 9.15]{Lancaster1995}. However, the particular logarithmic norm $\mu(\cdot)$ used throughout the prior work in \cite{DelMoral/Tugaut:2016,Bishop/DelMoral/Niclas:2017,2017arXiv171110065B,Bishop/DelMoral/multiDimRicc,Bishop/DelMoral/STV2018} is unimportant, i.e. the matrix norm defining the log-norm can be chosen arbitrarily. From \cite[Theorem 5]{strom1975logarithmic}, it is now known that if $\mathrm{Absc}(A-{{P}}_{\infty}\,S)<0$, then there exists a particular log-norm $\mu(\cdot)$ such that $\mu(A-{{P}}_{\infty}\,S)<0$. Thus, the conditions antecedent in much of the prior work, in which we ask for $\mu(A-{{P}}_{\infty}\,S)<0$ can be replaced with just asking for classical detectability and stabilisability conditions. (In prior work it was discussed and claimed by the authors that asking for $\mu(A-{{P}}_{\infty}\,S)<0$ may be viewed as asking for a type of strong observability and controllability. This may be true if one specifies first the log-norm of interest. But if one does not care which log-norm is used, we can significantly relax the setting and simply ask for detectability/stabilisability which leads to $\mathrm{Absc}(A-{{P}}_{\infty}\,S)<0$). Much of the analysis, as already discussed, requires $S\in\bbS_d^+$ which automatically implies detectability (in fact a much stronger condition than observability).

It is worth noting again that all moment boundedness and fluctuation results stated in this work hold with any $\sfN\geq 1$ and without further assumptions if one replaces the constants $c,c_{n},c_{n}(Q),c_{n}(z,Q)\ldots$ with functions that depend on (and grow with) the time horizon $t\geq0$.

\subsection{Bridging the Gap to Nonlinear Ensemble Filtering}

The focus of this article is ensemble filtering in the linear-Gaussian (continuous-time) setting. The results surveyed herein portray a rather detailed theory of fluctuation and stability/contraction results in that case. In practice, the ensemble Kalman filtering methodology is applied in high-dimensional, nonlinear state-space models \cite{evensen03,evensen2009book}. The evolution equations for each ensemble member in the case of nonlinear state-space models are given in (\ref{EnKF-sampled-diffusions-nonlinear}).

In \cite{CrisanXiong} and in \cite{Yang2016,Taghvaei2016ACC} the novel idea of a McKean-Vlasov-type diffusion which has conditional distribution equal to the true Bayesian filter is studied, see also \cite{Pathiraja2021}. The mean-field approximation of this McKean-Vlasov-type diffusion in \cite{Yang2016,Taghvaei2016ACC}, termed the \textit{feedback particle filter}, resembles somewhat superficially the ensemble filters in (\ref{EnKF-sampled-diffusions-nonlinear}). However, the analogue of the gain function in (\ref{EnKF-sampled-diffusions-nonlinear}) in the feedback particle filter of \cite{Yang2016,Taghvaei2016ACC} is derived as the solution of a certain Poisson-type partial differential equation. In the linear-Gaussian case, the filter of \cite{Yang2016} coincides with the {\tt DEnKF}.

In the nonlinear model setting, the ensemble filters in (\ref{EnKF-sampled-diffusions-nonlinear}) are not derived as sampled versions of an equation whose (conditional) distribution is equal to the Bayesian filter. That is, these filters are not derived as sampled versions of the McKean-Vlasov-type diffusion in \cite{Yang2016,Taghvaei2016ACC}. Conversely, in the limit ($\sfN\rightarrow\infty$) the ensemble filters in (\ref{EnKF-sampled-diffusions-nonlinear}) do not converge to an object with distribution equal to the optimal Bayes filter. In fact, the object these filters converge to has not been rigorously established in general and its properties, as compared to the true Bayesian filter, remain an open topic. Thus, in the nonlinear setting, the ensemble filters discussed in this work, see (\ref{EnKF-sampled-diffusions-nonlinear}), may be viewed as approximations of the feedback particle filter in \cite{Yang2016,Taghvaei2016ACC} only in some very weak sense (despite any superficial resemblance to the contrary). Indeed, the gain function approximation in (\ref{EnKF-sampled-diffusions-nonlinear}) is likely a very poor approximation of the solution of the Poisson-type partial differential equation in \cite{Yang2016}; except of course in linear-Gaussian models. Rather, we may argue, as we have earlier in this article, that the ensemble filters in (\ref{EnKF-sampled-diffusions-nonlinear}) should be viewed in the context of so-called observer theory, and related not to Bayesian filtering but rather to the more general topic of (dynamic) state estimation \cite{anderson79,Baras88}. The goal of state estimation in this context is to design an observer that tracks in some suitable (typically point-wise) sense the underlying signal and perhaps provides some usable measure of uncertainty on this estimate. The goal is not to develop an approximation (at each time) of the true conditional (Bayesian) distribution of the signal given the observations. The latter contains significantly more information than is perhaps needed in many practical applications. Nevertheless, we also argue that the filtering ideas in \cite{Yang2016,Taghvaei2016ACC}, and suitable approximations thereof, are in need of further investigation. 

In \cite{DelMoral/Kurtzmann/Tugaut:2017} a class of so-called ensemble extended Kalman filters ({\tt En-EKF}) is developed that is based on a type of particle approximation of the linearisation-based extended Kalman filter, see \cite{anderson79}. This ensemble filter is interesting because the sample mean is shown to converge (in $\sfN\rightarrow\infty$) to the extended Kalman filter state estimate. This extended Kalman state estimator has been widely studied in nonlinear filtering and control theory \cite{anderson79,Baras88,Reif2000,BishopJensfelt2009,karvonen2018stability}, and may be viewed more as a type of nonlinear state estimator rather than a Bayesian filter \cite{Baras88,Reif2000}. 

When considering nonlinear signal models, the long time behavioural analysis of various {\tt EnKF} methods in \cite{Kelly2014,tong2016stability,DelMoral/Kurtzmann/Tugaut:2017} assumes a strong type of stability property on the signal (which in the linear case would reduce to assuming that $A$ is Hurwitz stable in our model (\ref{lin-Gaussian-diffusion-filtering})). This stability assumption on the true signal is precisely what we aim to relax in our work; albeit limited in our study to linear models. Filter stability without assumptions on the stability of the true signal will ultimately require some control of the fluctuation properties of the sampled observer, e.g. see the discussion in the preceding section on this topic (in the linear-Gaussian model setting). This fluctuation analysis is lacking somewhat in the nonlinear model setting. It is complicated in that case by the absence of any closed-form evolution equations for the sample mean and sample covariance. 

Viewing, or even designing, an ensemble filter (or its sample mean for example) as a (dynamic) state estimator (or observer) may have some benefits. In particular, stability may be a larger design consideration if starting from this viewpoint rather than seeking Bayesian probabilistic properties. It may be possible to then also exploit the properties of existing nonlinear state estimators which have traditionally been rigorously analysed, e.g. \cite{Baras88,Reif2000,BishopJensfelt2009,karvonen2018stability}.

This is exemplified in the ({\tt En-EKF}) in \cite{DelMoral/Kurtzmann/Tugaut:2017} that converges to the extended Kalman filter in the limit $\sfN\rightarrow\infty$. The stability of the extended Kalman filter as a nonlinear observer has been widely studied, e.g. see \cite{Baras88,Reif2000,BishopJensfelt2009,karvonen2018stability}. Although strong signal stability assumptions are taken in \cite{DelMoral/Kurtzmann/Tugaut:2017}, it would be natural to consider the ({\tt En-EKF}) in \cite{DelMoral/Kurtzmann/Tugaut:2017} without the underlying signal stability assumption and look at developing the fluctuation type analysis considered herein in the linear-Gaussian setting. We may then also exploit the stability analysis that already exists \cite{Baras88,Reif2000,BishopJensfelt2009,karvonen2018stability} for the limiting extended Kalman state estimator. This is analogous in many ways to the stability properties and observability/controllability properties used herein in the linear-Gaussian setting.

Inflation is used in \cite{Kelly2014, tong2016inflation} in the nonlinear model setting to aid in stability. This is similar to the study considered herein on stability under inflation in linear-Gaussian models. It seems natural that added inflation acts to stabilise the various ensemble filters. In the context of the preceding discussion, inflation-based state estimators may also be viewed in the context of stable nonlinear observers, rather than heuristic adaptions of approximate Bayesian filters.

Finally, we remark that the transport-based ensemble filter {\tt DEnTF}, see case ({\tt NF3}) in (\ref{EnKF-sampled-diffusions-nonlinear}), is studied in \cite{deWiljes2018} in a particular nonlinear setting. Non-asymptotic (i.e. finite sample) uniform-in-time accuracy and stability of the {\tt DEnTF} is studied in \cite{deWiljes2018} under the assumptions of small observation noise, and a square observation matrix, or in other words with linear observations and a change of coordinate so that $H=I$. Note this latter assumption is made also in \cite{Kelly2014,tong2016stability,DelMoral/Kurtzmann/Tugaut:2017,deWiljes2019} which otherwise consider certain classes of nonlinear signals and different {\tt EnKF} variants. Thus, this strong (and linear) observability assumption seems key to analysis in the ensemble filtering literature even when moving away from the linear signal model.

\subsection{Other Related Literature}

The focus of this article is a detailed fluctuation and contraction analysis of the relevant ensemble Kalman filtering terms (e.g. the exponential semigroups, sample mean, and sample covariance) in the linear-Gaussian setting. There is considerable work on the periphery of this rather specific topic and analysis. A broad overview of the literature on filtering is not possible, but we note below some topics and literature for further study. 

For example, the introduction of this article was focused primarily on (specific) related literature in ensemble Kalman approximation methods, and so-called data assimilation. The topic of filter stability in the case of the true nonlinear filter (e.g. as given by the Kallianpur-Striebel formula \cite{kallianpur}) has been studied widely; see a broad but necessarily incomplete snapshot in \cite{Kunita1971,oconepardoux,Atar1998,Bhatt2000,Budhiraja2003,Baxendale2004,VanHandel2009,VanHandel2009ccc}. Known results in general suggest that sufficiently informative observation processes and/or sufficiently contractive/ergodic latent signals translate to contractive stability of the nonlinear filter. Emphasis on the continuous-time, linear-Gaussian model setting is studied in \cite{anderson71,oconepardoux,VanHandel2009,Bishop/DelMoral:2016,bd-CARE} where deeply understood observability and controllability model conditions explicitly formalise the properties leading to filter stability. See also \cite{Bonnabel2013} for a contraction analysis of a low-rank Kalman-Bucy filter particularly relevant in the application domain of this article.

The fluctuation and stability of other approximations schemes in nonlinear filtering have also been studied; e.g. see \cite{LaScala1995,Budhiraja1999,mf-dm-04,Chopin2004,Oudjane2005,Heine2008,VanHandel2009bb,Whiteley2013,Douc2014,BishopBonilla2022} for a snapshot of some of these methods. A detailed discussion of these approximation methods and their fluctuation and stability properties is beyond the scope of the ensemble Kalman-type methods studied here. However, the monograph \cite{mf-dm-04} provides a detailed study of the fluctuation and contractive properties of numerous particle filtering methods. We note in passing that in general when studying the stability of filtering approximation schemes it is rather common to assume the latent signal is stable/ergodic and/or the observation sequence is stationary, e.g. see \cite{Budhiraja1999,mf-dm-04,VanHandel2009bb,Douc2014,BishopBonilla2022}. In this article, and in prior work \cite{DelMoral/Tugaut:2016,Bishop/DelMoral/Niclas:2017,2017arXiv171110065B,Bishop/DelMoral/multiDimRicc}, we relax those assumptions and primarily rely on a strong form of observability (i.e. a strong form of observational informativeness). The results presented here \textit{do not} rely on \textit{any} form of latent signal stability in general. This latter fact distinguishes this work, but also the ensemble Kalman approximation method in its ability to handle totally unstable latent signals. The latter point was illustrated above via comparisons with the bootstrap particle filter in one-dimensional unstable linear-Gaussian models.

In the introduction we list some of the seminal data assimilation and ensemble Kalman methodology articles. Here, we give a by no means complete list of some methodological approaches that in some sense can be considered cousins of  

Different particle filtering methods based on evolving particles according to an ordinary differential equation that bridges the prior or predictive distribution with the posterior have been developed, e.g. see \cite{Daum2010,bunch2016}. These so-called Bayesian homotopy methods aim to introduce nonlinearities arising from the Bayes update in a tempered fashion. See also \cite{Reich2022} for an introduction with connections to ensemble Kalman filtering. Applications of the multilevel Monte Carlo method \cite{giles2008} to ensemble Kalman filtering have been considered in, e.g., \cite{Hoel2016,chada2022}. The continuous-time linear-Gaussian case similarly formulated as in this article is considered in detail in \cite{chada2022}. A related extension in multi-index ensemble Kalman filtering was proposed in \cite{Hoel2022}. The problem of unbiased ensemble Kalman filtering has been considered in \cite{Alvarez2022} with an emphasis and analysis closely related to the formulation considered herein. 

The ensemble Kalman filter has been applied to the problem of (log) normalisation constant estimation for continuous-time filtering problems; e.g. see \cite{CrisanDM2022,Ruzayqat2022}. More generally, ensemble Kalman methods for inverse problems have also been considered in the literature \cite{Iglesias2013,chada2018} with some related analysis \cite{Schillings2017a,Schillings2017b}. Particle filtering and ensemble filtering methods have also been applied in optimisation \cite{Zhang2018}. See these references for further details on the respective topics.

It was noted in the preceding section that certain McKean-Vlasov equations can be derived \cite{Yang2016} that superficially resemble the ensemble filters in (\ref{EnKF-sampled-diffusions-nonlinear}), but with gain functions derived as the solution of certain Poisson-type partial differential equations. In general, the (conditional) law of these nonlinear diffusions is equal to the filtering distribution. Applying mean-field particle approximations in the linear-Gaussian case, the filter of \cite{Yang2016} coincides with the {\tt DEnKF}. In the general nonlinear setting, other gain function approximations can be employed such as discussed in \cite{Taghvaei2020,Taghvaei23}.

Finally, we remark in passing that different models involving backward matrix Riccati diffusions arise in linear-quadratic optimal control problems with random coefficients; see e.g. \cite{bismut1976linear,hu2003indefinite,Kohlmann2003}. Another class of random Riccati equations, different from the Riccati (matrix quadratic) diffusion equations studied herein arises in network control and filtering with random observation losses; see e.g. \cite{sinopoli2004kalman,Tanwani2020}. The details of these works are beyond the scope of the forward-in-time Riccati diffusions considered herein. The forward-in-time Riccati diffusion of (\ref{EnKF-1}) is nevertheless of interest on its own, as with $\kappa=0$ it is a prototypical model of a matrix stochastic differential equation with a quadratic drift term. In the deterministic setting, it is worth noting that the contraction of the Riccati equation, e.g. (\ref{nonlinear-KB-Riccati}), (\ref{ref-phi-1-stability}), can naturally be studied with different metrics, e.g. in the Riemannian space of positive definite matrices, see e.g. \cite{Bougerol1993,Lee2008,Bonnabel2013,Levy2016}, and such contraction results may be of interest and/or practical value in the stochastic setting of the Riccati diffusion describing the flow of sample covariances.

\subsection{Some Open Problems}

The stochastic analysis and stability of the {\tt EnKF} models considered in this article are rather well understood in the linear-Gaussian, continuous-time, setting \textit{even with unstable latent signal processes}. However, there are still some worthy gaps. In particular, we may contrast the very strong and complete picture in the one-dimensional setting (following from \cite{2017arXiv171110065B} and presented above) with the more general and relevant multivariate setting. 

For example, the $\mathbb{L}_n$ contraction estimates on the Ricatti diffusion flow describing the sample covariance in (\ref{uniformly-Lipschitz-moments-riccati-scalar}) are only available in the scalar case. We know, e.g. see Theorem \ref{theo-stab-intro}, that the Markov semigroup is contractive and that the law of the sample covariance converges to an invariant measure in the multi-dimensional setting, analogously to the convergence of the deterministic Riccati equation to an equilibrium state. However, multi-variate versions of the contraction estimates (\ref{uniformly-Lipschitz-moments-riccati-scalar}) are unknown.

Moreover, multi-variate contraction estimates for the {\tt EnKF} mean as in (\ref{uniformly-Lipschitz-moments-OU-scalar}) are known only on certain finite time intervals (with a finite sample size) and a general multi-variate filtering contraction estimate as in (\ref{uniformly-Lipschitz-moments-OU-scalar}) under mild assumptions would be of interest. There is also a disparity in available results in the multi-variate setting between the {\tt VEnKF} and {\tt DEnKF}. 

In both (\ref{uniformly-Lipschitz-moments-riccati-scalar}) and (\ref{uniformly-Lipschitz-moments-OU-scalar}) in the scalar case we recover in the limit with the ensemble size the deterministic optimal rates of the classical Kalman-Bucy filter. Extensions of these rate estimates to the multi-variate setting would also be of interest.

Finally in reference to the scalar setting, we can explicitly state in closed form the invariant measure of the sample covariance, e.g. see (\ref{invariant-measure-venkf}) for the {\tt VEnKF} and (\ref{invariant-measure-denkf}) for the {\tt DEnKF}, and Figure \ref{fig:invariantmeasureEnKFvsDEnKF}. It is unlikely that closed-form expressions are possible to derive in the multivariate setting, however, it would be of further interest to confirm various properties such as the heavy-tailed nature of the stationary measure in the case of the {\tt VEnKF}. Such properties as discussed earlier have practical consequences such as the potential to lead to so-called catastrophic divergence and numerical instability. 

While not given explicitly, we may crudely introduce the signal dimension into the constants of the results presented in this work. However, a detailed study of the dimension as it pertains to stochastic fluctuation properties of the relevant sample covariance and subsequently ensemble mean is warranted since dimensionality versus computational expense is a primary driver of these methods in practice. Work in this direction with various covariance inflation mechanisms was considered in \cite{Bishop/DelMoral/Pathiraja:2017} in the linear-Gaussian setting and more generally in \cite{deWiljes2019} for the transport-inspired ensemble Kalman-Bucy filter {\tt DEnTF}. 

The main assumptions employed throughout are Assumptions \ref{mainAssumpObs} and \ref{mainAssumpCon}. The Assumption \ref{mainAssumpCon} can be relaxed to just stabilisability, and is thus completely in line with classical Kalman filtering analysis. It would be of interest to relax Assumption \ref{mainAssumpObs} to observability or even just detectability in line with classical Kalman-Bucy filtering. We have discussed in the preceding sections various issues issues surrounding the need for this stronger observability-type Assumption \ref{mainAssumpObs}, and difficulties with its possible relaxation. 

This article focused on the continuous-time linear-Gaussian model case, and the preceding open problems also fall under this setting. Direct extensions of the results presented in this article to the discrete-time linear-Gaussian model setting would be of interest and some work has been published in this case, particularly in the scalar setting, e.g. see \cite{dm-horton-21} and Theorem \ref{timeuniform-discrete-time} and the earlier discussion. Finally, without considering specifics and re-discussing the relevant literature, extensions of these results to the nonlinear model setting is of great practical interest, particularly under testable and natural model assumptions (e.g. observability-type assumptions) also accommodating unstable or transient latent signal processes (with the latter relaxation being a major driver of the results presented herein).

%

\end{document}